\documentclass[12pt,reqno]{amsart}
\usepackage{a4wide,enumerate,color}
\usepackage[left=2cm, right=2cm, bottom = 3cm]{geometry}
\usepackage{amsmath}
\usepackage{bm}
\usepackage{scrextend}% for labeling environment
\usepackage{graphicx}
\usepackage{tikz}
\usepackage{pgfplots}
\usepackage{caption}
\usepackage{subfig}
\usepackage{enumerate}
\usepackage{xcolor}
\usepackage[colorlinks]{hyperref}
\hypersetup{colorlinks=true,linkcolor=blue,filecolor=blue,citecolor = blue,urlcolor=blue,}

\usepackage{amsthm}
\usepackage{amsfonts}
\usepackage{amssymb}
\usepackage{mathrsfs}
\usepackage{esint}
\colorlet{darkgreen}{green!45!black}
\usepackage{graphicx}
\usepackage{caption}

\usepackage{tikz}
\usepackage{pgfplots}
\usepackage{pgfplotstable}
\usepackage{etex}
\usepackage{subfig}
\usetikzlibrary{shapes}
\usepackage{pifont}
\usepackage{multicol}

\usepackage{mathtools}
\usepackage{esint}

\newtheorem{theo}{Theorem}[section]
\newtheorem{prop}[theo]{Proposition}
\newtheorem{lem}[theo]{Lemma}
\newtheorem{coro}[theo]{Corollary}
\theoremstyle{definition}
\newtheorem{defini}[theo]{Definition}
\newtheorem{Remark}[theo]{Remark}

\numberwithin{equation}{section}

\newcommand{\rom}{\rho_-}
\newcommand{\rop}{\rho_+}
\newcommand{\ros}{\rho_\sharp}
\newcommand{\rot}{\hat\rho_+}
\newcommand{\rob}{\hat\rho_-}
\newcommand{\ov}{\overline}
\newcommand{\dw}{\downarrow}
\newcommand{\be}{\begin{equation}}
\newcommand{\ee}{\end{equation}}
\newcommand{\longto}{\longrightarrow}
\newcommand{\pt}{\partial}
\newcommand{\oo}{\infty}
\newcommand{\MM}{\mathbf{M}}

\newcommand{\R}{\mathbb{R}}

\newcommand{\Ups}{\Upsilon}
\newcommand{\TT}{\mathcal{T}}
\newcommand{\vhi}{\varphi}
\newcommand{\eps}{\varepsilon}
\newcommand{\te}{\theta}
\newcommand{\PI}{{\varPi}}
\newcommand{\om}{\omega}
\newcommand{\Om}{\Omega}
\newcommand{\lt}{\left}
\newcommand{\rt}{\right}
\newcommand{\st}{\stackrel}
\newcommand{\ds}{\displaystyle}
\newcommand{\sub}{\subset}

\newcommand{\loc}{_{\text{loc}}}
\newcommand{\restr}{
	\hskip2.5pt{\vrule height7pt width.5pt depth0pt}
	\hskip-.2pt\vbox{\hrule height.5pt width7pt depth0pt}
	\, }
\newcommand{\ie}{\textit{i.e.~}}

\DeclareMathOperator{\Id}{Id}
\DeclareMathOperator{\Lip}{Lip}

\DeclareMathOperator{\supp}{supp}

\newcommand{\cm}[1]{}

\newcommand{\A}{\mathcal{A}}
\newcommand{\PSI}{{\bs\Psi}}

\newcommand{\Co}{\mathcal{C}}

\renewcommand{\t}{\times}

\newcommand{\up}{\uparrow}
\newcommand{\pf}{\,\!_\#\,}

\newcommand{\Leb}{\mathcal{L}}
\newcommand{\Wass}{\mathbf{W}_2}
\newcommand{\wass}{\mathcal{W}_2}
\newcommand{\Abb}{\mathbb{A}}

\newcommand{\J}{\mathbf{J}}
\newcommand{\I}{\mathbf{I}}

\newcommand{\bs}{\boldsymbol}

\newcommand{\com}[1]{}

\usepackage{caption}
\usepackage{subfig}

\usepackage{array,multirow,makecell}
\usepackage{color}

\begin{document}

\title[Analysis of a free boundary model]{Analysis of a one dimensional energy dissipating free boundary model with nonlinear boundary conditions.\\
Existence of  weak solutions}

\author[B. Merlet]{Beno\^it Merlet}
\address{Univ. Lille, CNRS, UMR 8524, INRIA - Laboratoire Paul Painlev\'e, F-59000 Lille, France.}
\email{benoit.merlet@univ-lille.fr}

\author[J. Venel]{Juliette Venel}
\address{Univ. Polytechnique Hauts-de-France, CERAMATHS, FR CNRS 2037, F-59313 Valenciennes, France.}
\email{juliette.venel@uphf.fr}

\author[A. Zurek]{Antoine Zurek}
\address{Universit\'e de Technologie de Compi\`egne, LMAC, 60200 Compi\`egne, France.}
\email{antoine.zurek@utc.fr}

\date{\today}

\begin{abstract}
This work is part of a general study on the long-term safety of the geological repository of nuclear wastes. A diffusion equation with a moving free boundary in one dimension is introduced and studied. The model describes some mechanisms involved in corrosion processes at the surface of carbon steel canisters in contact with a claystone formation.  The main objective of the paper is to prove the existence of weak solutions to the problem which are maximal in time. For this, a time semidiscrete  minimizing movements scheme based on a Wasserstein-like distance is introduced. The existence of solutions to the scheme is proved. Then, using \textit{a priori} estimates, it is shown that as the time step goes to zero these solutions converge up to extraction towards a maximal weak solution to the free boundary model. 
\end{abstract}

%\paragraph{Keywords:}

%\paragraph{AMS classification:}

\maketitle

\tableofcontents

For ease of reading, we list here some important notations in their order of appearance in the article.
\subsection*{Index of Notation}
%\label{Sion}
\begin{enumerate}[(a)]
% Equation~\eqref{prime} page~\pageref{prime}. 
\item[{\makebox[3.5em][c]{$'$}}] When $w$ is a function of $x$ or of $(x,t)$ and the distributional derivative $\pt_xw$ is a measure we denote by $w'$ the density of the absolutely continuous part of $\pt_x w$.\\
 For differentiable functions of $x$, $w'$ denotes the classical derivative of $w$. 
\item[{\makebox[3.5em][c]{$\dot~$}}]  $\dot w$ denotes the time derivative of a function $w=w(x,t)$ or $w=w(t)$.
\item[{\makebox[3.5em][c]{$\rom,\rop$}}] Equation~\eqref{P.bord1} page~\pageref{P.bord1}. Given data $\rop\ge\rom>0$.
\item[{\makebox[3.5em][c]{$\beta,\theta$}}] Equation~\eqref{def_betatheta} page~\pageref{def_betatheta}, $\beta\in\R$ and $\theta>0$ are defined by  $-\beta\pm\te=\log\rho_\pm$.
\item[{\makebox[3.5em][c]{$f$}}] Equation~\eqref{f} page~\pageref{f}. This function is defined by $f(r):=r(\log r +\beta - 1)$ for $r\ge0$.
\item[{\makebox[3.5em][c]{}}]  We have $f'(r)=\log r +\beta$\ and\ $ f''(r)=1/r$\ for $r>0$. Moreover $\min f=-e^{-\beta}=-\sqrt{\rom\rop}$.
\item[{\makebox[3.5em][c]{$\PSI$}}] Equation~\eqref{PSI} page~\pageref{PSI}. Gibbs energy.
\item[{\makebox[3.5em][c]{$D_T$}}] Equation~\eqref{def_D_T} page~\pageref{def_D_T}. Space-time domain.
\item[{\makebox[3.5em][c]{$\MM$}}] Equation~\eqref{def_MM} page~\pageref{def_MM}, $\MM(\rho):=\int_{\R_+}\rho(x)\, dx$.
\item[{\makebox[3.5em][c]{$\rob,\rot$}}] Equation~\eqref{rhominrhomax} page~\pageref{rhominrhomax}. Lower and upper bounds for the  density of the solutions.
\item[{\makebox[3.5em][c]{$\wass$}}] Equation~\eqref{2.def.wass} page~\pageref{2.def.wass}. The 2-Wasserstein distance.
\item[{\makebox[3.5em][c]{$\mu_{\rho^0,\rho},\,\mu^0$}}] Equation~\eqref{2.mu1} page~\pageref{2.mu1} and equation~\eqref{3.mu1} page~\pageref{3.mu1}.
\item[{\makebox[3.5em][c]{$\Wass$}}] Equation~\eqref{def_W2} page~\pageref{def_W2}. ``Unbalanced'' 2-Wasserstein distance.
\item[{\makebox[3.5em][c]{$\ell_+,\ell_-$}}]  Lengths defined in~\eqref{ell+}, \eqref{ell-} (and~\eqref{ell+ell-0}) page~\pageref{ell+}.
\item[{\makebox[3.5em][c]{$\Abb$}}] Equation~\eqref{SpaceAbb} page~\pageref{SpaceAbb}. Space of solutions for the optimization problem~\eqref{OP}.
\item[{\makebox[3.5em][c]{$x_\rho$}}] Equation~\eqref{xrho} page~\pageref{xrho}, $x_\rho:=\inf\{x\ge0:\supp\rho\sub[0,x]\}$.
\item[{\makebox[3.5em][c]{$\A$}}] Equation~\eqref{SpaceA} page~\pageref{SpaceA}. Precised space of solutions for the optimization problem~\eqref{OP}. 
\item[{\makebox[3.5em][c]{$\textbf{d}_1$}}] Equation~\eqref{2.def_d1} page~\pageref{2.def_d1}. Linear part of the metric distance in $\Abb$.
\item[{\makebox[3.5em][c]{$\textbf{d}_2^{\,2}$}}] Equation~\eqref{2.def_d2} page~\pageref{2.def_d2}. Quadratic part of the metric distance in $\Abb$.
\item[{\makebox[3.5em][c]{$p_\tau$}}] Equation~\eqref{2.def.penalty} page~\pageref{2.def.penalty}. Penalty functional.
\item[{\makebox[3.5em][c]{$\J^\tau_{(\rho^0,X^0)}$}}]  Equations~\eqref{Jtau},~\eqref{Jtau2} page~\pageref{Jtau}.  The functional minimized in the optimization problem~\eqref{OP}. 
\item[{\makebox[3.5em][c]{$\I^\tau_{(\rho^0,X^0)}$}}]  Equation~\eqref{Itau} page~\pageref{Itau}.   Functional minimized by $\rho$, assuming $\MM(\rho)$ is fixed.
\item[{\makebox[3.5em][c]{$\TT$}}] Equation~\eqref{Tmax} page~\pageref{Tmax}. A temporary bound on the maximum existence time of solutions.
\item[{\makebox[3.5em][c]{$\rho^n$}}] Equation~\eqref{JKOscheme1} page~\pageref{JKOscheme1}. Discrete in time density of oxygen.
\item[{\makebox[3.5em][c]{$X^n$}}] Equation~\eqref{JKOscheme1} page~\pageref{JKOscheme1}. Discrete in time width of the oxyde layer.
\item[{\makebox[3.5em][c]{$\rho_\tau$}}] Equation~\eqref{def.rhotauXtau} page~\pageref{def.rhotauXtau}. Piecewise constant interpolation of $\rho^n$.
\item[{\makebox[3.5em][c]{$X_\tau$}}] Equation~\eqref{def.rhotauXtau} page~\pageref{def.rhotauXtau}. Piecewise constant interpolation of the sequence $X^n$.
\item[{\makebox[3.5em][c]{$M^n$}}] Equations~\eqref{Mn} page~\pageref{Mn}.  $M^n=\MM(\rho^n)$.
\item[{\makebox[3.5em][c]{$m^n$}}] Equations~\eqref{mn} page~\pageref{mn}. Discrete time derivative of $M^n$, $m^n=M^{n+1}-M^n$.
\item[{\makebox[3.5em][c]{$\ell_+^n,\ell_-^n$}}]  Lengths defined in~\eqref{ell+n}, \eqref{ell-n} page~\pageref{ell+n}.
\item[{\makebox[3.5em][c]{$\hat T$, $\hat\tau$}}] Positive constants introduced page~\pageref{prop.Ttau} in Proposition~\ref{prop.Ttau}(d).
\item[{\makebox[3.5em][c]{${\tilde X}_\tau$}}] Equation~\eqref{def.tildeXtauMtau} page~\pageref{def.tildeXtauMtau}. Continuous interpolation of $X^n$.
\item[{\makebox[3.5em][c]{$M_\tau$}}] Equations~\eqref{def.tildeXtauMtau} page~\pageref{def.tildeXtauMtau}. Continuous interpolation of $M^n$.
\item[{\makebox[3.5em][c]{$\sigma_{-\tau}$}}] Equation~\eqref{sigmatau} page~\pageref{sigmatau}. Time shift operator.
\item[{\makebox[3.5em][c]{$H^*$}}] Dual space of $H^1(\R_+)$ introduced in Proposition~\ref{prop_apriori}, page~\pageref{prop_apriori}. 
\item[{\makebox[3.5em][c]{$\delta$}}] Equation~\eqref{defdelta} page~\pageref{defdelta}. $\delta=\delta(\tau,\om'):=\tau^{\omega'}$ is the characteristic length of the mollifier $\zeta_\delta$. 
\end{enumerate}

\section{Introduction}
%\label{sec.intro}

This work is motivated by the study of the so-called Diffusion Poisson Coupled Model (DPCM) introduced in~\cite{BBC10} by Bataillon\textit{et al.} This system was developed to model the corrosion of a steel plate in contact with a solution. In particular, it is relevant for describing the corrosion of steel canisters containing nuclear wastes which are confined in a glass matrix and stored at a depth of several hundred meters in a claystone layer. Since this storage method is considered by various countries, its reliability requires investigations. In particular, wastes stay radioactive for several hundred of years and it is important to understand the \emph{long term} behaviour of the system. Our main concern is about corrosion and the quantity of hydrogen molecules released during the process which can lead to safety issues. As it is not possible to perform physical experiments at these time scales, the use of reliable models (such as DPCM mentioned above) allowing \textit{in silico} experiments are required. However, to design accurate numerical methods capable of predicting the values of the relevant physical quantities over a long time, it is necessary to understand the mathematical properties of the model.

Let us briefly explain the main features of the DPCM. This is a one dimensional free boundary system. The space is decomposed in three regions: the oxide layer is in contact on one side with the claystone, viewed as an aqueous solution and on the other side with the metal. The DPCM is a system of drift-diffusion equations describing the evolution inside the oxide layer of charge carriers (electrons, Fe$^{3+}$  and  O$^{2+}$ cations) and coupled with a Poisson equation governing the dynamics of the electrical potential. The positions of the solution/oxide layer and oxide layer/metal interfaces evolve along time according to some given ordinary differential equations. Besides, the electrochemical reactions only occur at these interfaces (\textit{i.e.} there are no reaction terms in the bulk of the three regions). These electrochemical reactions are modelled by some nonlinear Fourier boundary conditions at the interfaces.

Due to the numerous coupling of the equations and their definition on a domain with free boundary the mathematical study of the DPCM is a challenging task. So far, only a few results are available in the literature. In~\cite{CLV12,CLV14} the well-posedness of the system has been proved for a simplified version of the DPCM where the positions of the interfaces are fixed. A finite-volume scheme approximating the solutions to the DPCM has been proposed in~\cite{BBCF12}. The numerical experiments with relevant physical data presented in~\cite{BBC10,BBCF12} suggest the existence of a global solution to the system. In particular, the existence of travelling wave solutions is established: after a transient time both interfaces move at the same speed, the width of the oxide domain remains constant and the charge carriers and the electrical potential admit a stationary profile. The existence of such travelling wave solutions has been proved in~\cite{CG15} for a reduced model which assumes electroneutrality in the oxide layer. Thanks to a computer-assisted proof, the existence of travelling wave solutions for the ``full'' DPCM has been obtained in~\cite{BCZ21}. More recently, in~\cite{CCMRV22} another simplified model (only two species are considered and the interfaces are fixed) is proposed with some changes in the nonlinear boundary conditions that correct a thermodynamical inconsistency of the initial model. This modification makes the mathematical study  more tractable in the case of a fixed domain but the well-posedness of the complete free boundary problem is still open.
 
Up to now, no existence result has been proved for the evolutionary DPCM with free boundaries. 
Some of the main difficulties come from the evolution of positions of the free boundaries which are governed by nonlinear equations.  
%One of the main difficulties for establishing the existence of a global solution is to justify that the width of the oxide layer, where the equations of the systems are defined, stays positive along time (as numerically suggested in~\cite{BBC10,BBCF12}). 
Besides, the structure of the system and in particular its lack of obvious gradient flow structure prevents the derivation of classical \textit{a priori} estimates which would lead (even formally) to a positive lower bound for the width of the oxide layer along time.\footnote{This general idea might possibly work but we have not found the right way to implement it.} To bypass this difficulty we adapt the approach of Portegies and Peletier in~\cite{PP10} which makes use of tools from optimal transport for studying a moving boundary problem. Indeed, in~\cite{PP10}, these authors introduce a one dimensional parabolic free boundary model with two moving interfaces describing the variation of the length of a piece of crystal by dissolution/precipitation. The thermodynamical consistency of the model is deeply connected with its gradient-flow structure with respect to some Wasserstein metric. In particular, the existence of solutions to this system is obtained thanks to a Jordan, Kinderlehrer and Otto (JKO) minimizing scheme (see~\cite{JKO98}). The relevance of this approach in the context of parabolic equations in a fixed domain is well known, see for instance~\cite{A05,JKO98,KMX17,Ot98,Ot01,OTAM} and the idea to recast some free boundary problems in the Wasserstein gradient-flow setting seems promising. 

The main goal of the present work is to show that this idea is effective for the DPCM model. For this, we consider a free boundary model which, compared to the DPCM, is very simple in the bulk of the oxide layer but retains all the difficulties related to the nonlinear boundary conditions and to the equation of motion of the oxide-metal interface.

\subsection{Presentation of the reduced model}

Let us first explain the phenomena that our reduced model has to describe. We only consider the evolution of the density $\rho$ of oxygen atoms inside the oxide layer.\footnote{In the original model of~\cite{BBC10} the equations are written in terms of the density of ``oxygen vacancies'' which corresponds to $1-\rho$ here.} We neglect the other charged species as well as the influence of the electrical potential in this domain so that $\rho$ satisfies a heat equation defined on a moving domain. In the original model, the concentration of oxygen in the oxide layer can not exceed some threshold. Here we will assume that $\rho$ is non negative and bounded by some constant $\rot>0$ at initial time. We will check that these properties propagate along time.

Next, we fix the position of the interface solution/oxide layer at $x=0$, while the interface oxide layer/metal is moving according to an ordinary differential equation. Finally, in order to take into account the chemical reactions at the interfaces we impose boundary conditions which model the exchange of matters at the interfaces of each regions, \textit{i.e.} at the interfaces solution/oxide layer and oxide layer/metal. More precisely, we denote by $X(t)>0$ the position of the moving interface at time $t$. The respective domains of the solution, of the oxide layer and of the metal at time $t$ are $(-\oo,0)$, $[0,X(t)]$ and $(X(t),+\oo)$. The density of oxygen in the oxide layer is denoted $\rho(x,t)$ for $0\le t\le T$, $0\le x\le X(t)$. We consider for $T>0$ the following free boundary model with unknowns~$(\rho,X)$:
\begin{subequations}
\label{P}
\begin{align}
\label{P.a}
\pt_t\rho(x,t) -\pt^2_x\rho(x,t) = 0 &\quad\mbox{for }x\in[0,X(t)],\, t\in[0,T],\\
\label{P.c}
\pt_x\rho(X(t),t) +\dot{X}(t)\,\rho(X(t),t) =0&\quad\mbox{for }t\in[0,T],\\
\label{P.e}
\lambda\,\dot{X}(t) =\rho(X(t),t)-\ros &\quad\mbox{for }t\in[0,T],\\
X(0) = X^0,&\\
\label{P.f1}
\rho(x,0) =\rho^0(x) &\quad\mbox{for }x\ge0,
\end{align}
\end{subequations}
where $\lambda$ and $\ros$ are given positive constants and the given initial data is composed of $X^0>0$ and $\rho^0:[0,X^0]\to(0,+\oo)$.

\begin{figure}[ht]
\centering
\captionsetup{width=1\textwidth}
\includegraphics[]{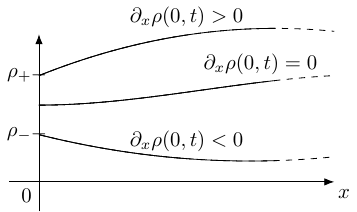}
\caption{The three possible cases of the boundary condition~\eqref{P.bord2}\label{figcbx=0}}
\end{figure}

We also have to specify a boundary condition on $\rho$ at the solution/oxide layer interface at $x=0$. In order to be consistent with the previous simplifications we should take a non-homogeneous \emph{linear} Fourier condition but we would lose an important difficulty in the problem. According to~\eqref{P.a}, $\pt_x\rho(0,t)$ represents the flux of oxygen  at time $t$ from the oxide layer into the solution, this quantity should be a nondecreasing function of $\rho(0,t)$ and more generally, we should have 
\[
\pt_x\rho(0,t)\in\pt F(\rho(0,t)),
\]
where $F$ is some (lower semicontinuous) convex function and $\pt F$ denotes its subderivative. We pick a worst case scenario (from the point of view of regularity) and choose $F$ as the indicatrix of $[\rom,\rop]$ for some 
\[
0<\rom\le\rop,
\]
that is $F(\rho)=0$ if $\rom\le\rho\le\rop$ and $F(\rho)=+\oo$ in the other cases. This leads to the following nonlinear conditions for $ t\in[0,T]$, see Figure~\ref{figcbx=0}. 
\begin{subequations}
\label{P.bord}
\begin{gather}
\label{P.bord1}
\rom\le\rho(0,t)\le\rop,\\
\label{P.bord2}\text{and }\qquad
\begin{cases}
\phantom{aaa}\rho(0,t) =\rop &\quad\mbox{if }\pt_x\rho(0,t)> 0,\\
\phantom{aaa}\rho(0,t) =\rom &\quad\mbox{if }\pt_x\rho(0,t)< 0,\\
\rom\le\rho(0,t)\le\rop &\quad\mbox{if }\pt_x\rho(0,t)= 0.
\end{cases}\phantom{aaa}
\end{gather}
\end{subequations}

The conditions in~\eqref{P.bord} can be seen as a generalization of a Signorini problem where the one sided-constraint on the solution model usually some irreversible phenomena on the boundary. They are used for instance in some unilateral contact problems in elasticity~\cite{BB03,HN07,KO88}, in some continuum mechanics models to describe a semipermeable membrane~\cite{DL76} or in chemistry to model electrochemical reacting interfaces~\cite{GHM99}. Here, the conditions in~\eqref{P.bord} allow the exchange of matters, in both ways, at the interface solution/oxide layer with two distinct thresholds $\rom$ and $\rop$.\\
From a physical point of view, the choice of~\eqref{P.bord} is disputable as the transport of oxygen is a reversible phenomenon. However, there is another phenomenon, neglected in the simplified model considered here, which is irreversible, namely the dissolution of the oxide layer. Indeed, the iron in the aqueous solution, rather than possibly reconstituting the oxide layer, will form oxide complexes, less organized and with a porous structure (rust). For this reason, considering nonlinear monotonic and \emph{non-smooth} boundary conditions like~\eqref{P.bord} anticipates future studies on a full model.\\
In the sequel we will use~\eqref{P.bord2} in the following equivalent form (assuming that~\eqref{P.bord1} holds true),
\be\label{Pbordeta}
\pt_x\rho(0,t)\, (\rho(0,t)-\eta)\ge 0\quad\text{ for every }\rom\le\eta\le\rop.\medskip
\ee

Let us now discuss the model at the interface $x=X(t)$, namely conditions~\eqref{P.c},~\eqref{P.e}. First observe that~\eqref{P.c} expresses the conservation of oxygen  as the position of the oxide/metal interface varies while the rate at which the interface moves is given by~\eqref{P.e}. This does not correspond to a motion of matter but to a change of state. More precisely, a change in the arrangement of iron atoms at the interface (from metal to oxide if $\dot X>0$), hence~\eqref{P.e}, is the rate of a chemical reaction. In the original model of~\cite{BBC10}, this rate is given by a formula of the form\footnote{In the original model, $k_1$ and $m_1$ depend on the electric potential which is assumed to be constant here.}
\[
\dot X= k_1\rho - m_1\lt(1-\rho\rt),
\]
for some $m_1,k_1>0$. With our notation, this corresponds to $m_1=\ros/\lambda$ and $k_1+m_1=1/\lambda$ and we see that in this situation there holds
\[
\ros=\dfrac{m_1}{m_1+k_1}\in(0,1).
\]
According to~\eqref{P.e}, the oxide layer may grow or shrink depending on whether $\rho>\ros$ or $\rho<\ros$ at point $X(t)$, see the discussion of Subsection~\ref{SubsectionSchrink}.

\subsection{Evolution of the Gibbs free energy}
%\label{subsec_Evolution_of_Gibbs_energy}

In~\cite{CCMRV22}, the authors consider another simplification of the DPCM and study its well-posedness. In that paper, the space domain is the interval $(0,1)$ (no free boundary) and the system describes the evolution of the density $u_1$ of ferric cations, the density of electrons $u_2$ and the electric potential $v_0$. The system is rewritten by expressing the fluxes $J_i$ of each chemical species (ferric cations for $i=1$ and electrons for $i=2$) on the form
\[
J_i=-\sigma_i\pt_x\xi_i,
\]
where $\sigma_i$ and $\xi_i$ are respectively the mobility and the electrochemical potential of the species $i$.  The study relies heavily on the evolution of the physical relevant local quantities (concentrations and electrothermodynamical potentials of each species) as well at the evolution of the global Gibbs energy of the system. As already said, we consider here instead the density of oxygen and by contrast, we retain the free boundary at the oxide/metal interface from the complete model of~\cite{BBC10}. The  mathematical treatment of the evolution of this interface is one of the two objectives of the present work (the other one being the treatment of the irreversible boundary condition at $x=0$). However, we follow~\cite{CCMRV22} by interpreting the system~\eqref{P} with a thermodynanymical point of view. More precisely, we consider that the flux of oxygen is given by 
\[
J=-\sigma\pt_x\xi,
\]
with 
\[
\sigma=1/\rho\qquad\text{and}\quad\xi=\log(\rho)+\beta,
\]
for some constant $\beta$. As soon as we consider the evolution of $\rho$ in the bulk of the oxide layer the chemical potential may be defined up to an additive constant and the value of $\beta$ is irrelevant. However, this is no longer true when we consider exchanges and chemical reactions at the boundary. Here, $\beta>0$ (and a new parameter $\te>0$) are defined from $\rom,\rop$  by the following equations
\be\label{def_betatheta}
\lt\{\begin{array}{c}\exp(-\beta-\te)=\rom\\ \exp(-\beta+\te)=\rop
\end{array}\rt\}\ \iff\\
\log\rho_\pm +\beta=\pm\te\ \iff\ 
\lt\{\begin{array}{l}\beta=-\dfrac12\log\lt(\rom\rop\rt)\\ \te=\dfrac12\log\lt({\rop/\rom}\rt)\end{array}\rt\}.
\ee

Now, we define the density of Gibbs energy function in the oxide layer as   
\be\label{f}
f(r):=r(\log r +\beta-1)\qquad\text{for }r\ge0,
\ee
so that 
\[
\xi=f'(\rho)=\log\rho +\beta.
\]
We also define the Gibbs energy density in the metal domain as $-\ros$. Then, if we consider an oxide domain $(0,X)$ and an oxygen density $\rho\in L^1(\R_+,\R_+)$ with $\rho\equiv 0$ on $(X,+\oo)$, the  associated Gibbs energy in a large domain $(0,\Lambda)$ (with $\Lambda>X$) is given by 
\[
\int_0^Xf(\rho(x))\,dx+\int_{X}^{\Lambda}(-\ros)\, dx=\int_0^X f(\rho(x))\,dx+\ros X -\ros\Lambda.
\] 
Finally, we get rid of the last constant term and define the total Gibbs energy of the pair $(\rho,X)$ as
\be\label{PSI}
\PSI(\rho,X):=\int_0^Xf(\rho(x))\,dx+\ros X.
\ee
Let us now study the evolution of the Gibbs energy of a solution~\eqref{P}--\eqref{P.bord}. Let $(\rho,X)$ be a smooth solution in the space-time domain 
\be\label{def_D_T}
D_T:=\{(x,t):0 < t<T,\, 0< x< X(t)\},
\ee
and let us denote
\[
\Psi(t):=\PSI(\rho(t),X(t)),
\]
with the notation $\rho(t):=\rho(\cdot,t) $.
We compute 
\[
\dot\Psi(t)=\int_0^{X(t)}\pt_t\rho(x,t)\lt(\log\rho(x,t)+\beta\rt)\, dx
+f(\rho_X)\dot X(t) +\ros\dot X(t) =:A(t)+B(t)+C(t),
\]
where we use the shorthands $\rho_X$ for $\rho(X(t),t)$. Similarly, in the computations below, we write  $\pt_x\rho_X$ for $\pt_x\rho(X(t),t)$, $\rho_0$ for $\rho(0,t)$ and $\pt_x\rho_0$ for $\pt_x\rho(0,t)$.\\
Using~\eqref{P.a} and integrating by parts, the first term rewrites as 
\[
A(t)
=-\int_0^{X(t)}\dfrac{(\pt_x\rho(x,t))^2}{\rho(x,t)}\, dx-\pt_x\rho_0 f'(\rho_0)+\pt_x\rho_X f'(\rho_X)=:A_0(t)+A_1(t)+A_2(t).
\]
We first notice that the bulk term $A_0(t)$ is nonpositive and  is therefore a dissipation term. The term $A_1(t)$ accounts for the fluxes through the interface $x=0$. We have 
\[
A_1(t)\, = -\pt_x\rho_0\,f'(\rho_0)\st{\eqref{f}}=-\pt_x\rho_0 (\log(\rho_0)+\beta).
\]
Combining the boundary condition~\eqref{P.bord2} with~\eqref{def_betatheta}, we see that either $\pt_x\rho_0=0$ or 
\[
\pm\pt_x\rho_0 >0\implies\rho_0=\rho_\pm\implies\log(\rho_0)+\beta=\pm\te.
\]
In any case, we have
\[
A_1(t)\,=-\te |\pt_x\rho_0 |\le0.
\]
Eventually, we gather the boundary terms at $X(t)$, namely $A_2(t)$, $B(t)$, $C(t)$. We compute
\begin{align*}
R_X(t)&:=A_2(t)+B(t)+C(t)\\
&\st{\phantom{\eqref{P.c}}}=\pt_x\rho_X f'(\rho_X)+f(\rho_X)\dot X(t)+\ros\dot X(t)\\
&\st{\eqref{P.c}}=\dot X(t)\big( -\rho_X f'(\rho_X)+f(\rho_X)+\ros\big)\\
&\st{\phantom{\eqref{P.c}}}=\dot X(t)\lt(-\rho_X+\ros\rt)\\
&\st{\eqref{P.e}}=-\lambda\lt(\dot X(t)\rt)^2\ \le0.
\end{align*}
In summary, we have for $t\in[0,T)$,
\be\label{dissip}
\dot\Psi(t)=-\int_0^{X(t)}\dfrac{\lt(\pt_x\rho(x,t)\rt)^2}{\rho(x,t)}\, dx -\te\left|\pt_x\rho(0,t)\right|-\lambda\lt(\dot X(t)\rt)^2,
\ee
and the three terms in the right hand side are nonpositive and correspond to different dissipation phenomena, respectively: in the bulk of the domain, at the left interface and at the right interface.\medskip

The form of the Lyapunov functional $\PSI$ is typical of Wasserstein gradient flow systems. This suggests using the techniques developed for the latter to analyze our model.  However, we cannot claim that~\eqref{P}--\eqref{P.bord} admits a gradient flow structure even in the generalized sense introduced by Mielke in~\cite{Mie11}. Indeed, in this theory, one has to specify an energy functional (or driving functional) and a \emph{quadratic} dissipation potential. In our case, using~\eqref{dissip}, the energy functional and the dissipation potential are clearly identified. But, due to the linear dissipation term $-\te\left|\pt_x\rho(0,t)\right|$ in~\eqref{dissip}, we cannot recast the dynamics of~\eqref{P}--\eqref{P.bord} in an Hilbertian setting and interpret this system as a generalized gradient flow. Nevertheless, the ideas developed in~\cite{JKO98,Ot01} can be applied in our case in order to prove the existence of weak solutions. In particular, we use a semidiscrete scheme inspired from the JKO-scheme.\medskip

\subsection{Notion of weak solution and main result}
%\label{subsec.weaksol}

In this subsection we define a notion of weak solution for system~\eqref{P}--\eqref{P.bord1}, first ignoring the boundary condition~\eqref{P.bord2} at $x=0$. Then, following the classical approach to deal with Signorini problem~\cite{DL76,KS00}, this later is expressed separately, in a weak form, as a variational inequality.

It is convenient to consider that the concentration of oxygen is defined on $\R_+$ and not only on the oxyde domain $[0,X(t)]$. The physically relevant choice is to set $\rho(x,t)=0$ for $t>X(t)$ (no oxygen in the metal domain). With this convention, for a (weak) solution $(\rho,X)$ of~\eqref{P}, we will have  
\[
\rho\in L^1_c(\R_+\times[0,T),\R_+),
\]
where $L^1_c$ denotes the space of compactly supported integrable functions (see Subsection~\ref{subsec.notation} below). In the sequel, for a pair $(\sigma,Y)$ representing a density $\sigma$ of oxygen in an oxide domain $[0,Y]$, at some fixed time $t$ (that is  $Y=X(t)$ and $\sigma=\rho(t)$), the function $\sigma$ is always continuous on $[0,Y$] and by convention $\sigma(x)=0$ for $x>Y$. In particular, $\sigma$ is left continuous and the value of $\sigma$ at $Y$ is not ambiguous: 
\[
\sigma(Y)=\lim_{\substack{x\to Y\\x<Y}}\sigma (x) =\sigma(Y^-,t).
\]

Let us derive a variational identity satisfied by any (sufficiently smooth) solution  $(\rho,X)$ to~\eqref{P}. Now, we introduce the following notation for $\sigma\in L^1(\R_+)$,
\be\label{def_MM}
\MM(\sigma):=\int_{\R_+}\sigma(x)\,dx.
\ee
Let $(\rho,X)$ be a solution of ~\eqref{P} and let $D_T$ be the space-time domain as in~\eqref{def_D_T}. We denote for $0\le t<T$,
\[
M(t):=\MM(\rho(t))=\int_{\R_+}\rho(x,t)\,dx =\int_0^{X(t)}\rho(x,t)\,dx .
\]
We first observe that for $0\le t<T$, 
\begin{align}
\nonumber
\dot M(t)&=\int_0^{X(t)}\pt_t\rho(x,t)\,dx +\dot X(t)(\rho(X(t),t)\\
\label{dotM}
&\st{\eqref{P.a}}=-\pt_x\rho(0,t)+\lt[\pt_x+\dot X(t)\rt]\rho(X(t),t)\st{\eqref{P.c}}=-\pt_x\rho(0,t).
\end{align}
To deal with weak solutions such that $M$ has $BV$-regularity, in the sequel, we consider $\dot{M}$ as the measure on $[0,T)$  defined by
\be\label{def:dotM}
\langle\dot{M},\psi\rangle:= -\int_{-\oo}^T M(t)\,\dot{\psi}(t)\, dt,\qquad\text{for every }\psi\in \Co_c([0,T)).
\ee

Let $\vhi\in\Co^\oo_c(\R_+\t [0,T))$,  multiplying~\eqref{P.a} by $\vhi$, integrating over $D_T$, 
integrating by parts the first term with respect to time and the second term with respect to space and using the conservative boundary condition~\eqref{P.c}, we get
\begin{multline*}
\iint_{D_T}\lt[-\rho(x,t)\pt_t\vhi(x,t)+\pt_x\rho(x,t)\pt_x\vhi(x,t)\rt]\,dx\,dt\\
-\int_0^{X(0)}\rho(x,0)\vhi(x,0)\,dx
+\int_0^T\pt_x\rho(0,t)\vhi(0,t)\, dt
=0.
\end{multline*}
Substituting the initial condition~\eqref{P.f1} and the boundary condition~\eqref{dotM} leads to
\begin{multline}\label{wf_0}
\iint_{D_T}\lt[-\rho(x,t)\pt_t\vhi(x,t)+\pt_x\rho(x,t)\pt_x\vhi(x,t)\rt]\,dx\,dt\\
-\int_0^{X^0}\rho^0(x)\vhi(x,0)\,dx
-\int_0^T\dot M(t)\vhi(0,t)\, dt
=0.\phantom{aaa}
\end{multline}

The above weak formulation has to be complemented with equation~\eqref{P.e}, that is $\lambda\,\dot{X}(t) =\rho(X(t),t)-\ros$, which prescribes the motion of $X(t)$. It is a standard procedure to check that~\eqref{P.e} and the weak formulation~\eqref{wf_0} are equivalent to~\eqref{P} as soon as $X$ is Lipschitz continuous, $M$ is $BV$ and  $\rho$ has regularity $L^2_tH^2_x\cap H^1_tL^2_x$ in the domain $D_T$.\medskip

We now derive a convenient weak formulation of the boundary conditions~\eqref{P.bord}. Let us assume that $(\rho,X)$ is a smooth solution to~\eqref{P}--\eqref{P.bord} on $(0,T)$  for some $T>0$ such that $\inf_{t\in(0,T)}X(t)>0$.\\
For the first part of the boundary condition (\ie~\eqref{P.bord1}), we consider the weaker formulation 
\be
\label{P.bord1w}
\rom\le\rho(0,t)\le\rop\qquad\text{for almost every }t\in (0,T).
\ee
Then, let $X_*>0$ with $X(t) \ge X_*$ for $t\in[0,T]$ and let $\chi\in\Co^\oo_c(\R_+)$ such that 
\be
\label{X*chi}
0\le\chi\le1,\qquad\chi\equiv1\ \text{ on }[0,X_*/2)\quad\text{ and }\quad\supp\chi\sub [0,3X_*/4).
\ee
We set $w(x,t):=\chi(x)\rho(x,t)$ for $(x,t)\in D_{T}$ and we compute
\be\label{1.ineg.var1}
\pt_t w(x,t) -\pt_x^2 w(x,t) = g(x,t)\quad\mbox{for }(x,t)\in [0,X_*)\t (0,T),
\ee
where the source term $g$ is given by
\[
g(x,t):= -\chi''(x)\,\rho(x,t) - 2\chi'(x)\,\pt_x\rho(x,t)\quad\mbox{for }(x,t)\in [0,X_*)\t (0,T).
\]
Let us also introduce some functions $\phi\in\Co^\oo_c([0,T),\R_+)$ and $\eta\in\Co^\oo_c(\R_+\t [0,T))$ with $\eta(0,t)\in [\rom,\rop]$ for all $ t\in[0,T]$. Thanks to the boundary conditions~\eqref{P.bord} we have (recall~\eqref{Pbordeta}), 
\be
\label{P.bordw}
\pt_x\rho(0,t)\, (\rho(0,t)-\eta(0,t))\ge 0\quad\text{ for every }t\in[0,T).
\ee
In fact, assuming~\eqref{P.bord1w} it is easily seen that~\eqref{P.bord2} holds true in $[0,T)$ if and only if~\eqref{P.bordw} holds true for every $\eta$ such that $\eta(0,t)\in[\rom,\rop]$.

Next, since $\phi\ge0$ and $w\equiv\rho$ in the neighborhood of $x=0$ we have
\be\label{1.ineg.var2}
\phi(t)\pt_x w(0,t)\, (w(0,t)-\eta(0,t))\ge0\quad\text{ for every }t\in[0,T).
\ee
Multiplying~\eqref{1.ineg.var1} by $\phi (w-\eta)$, integrating in space and time, integrating by parts and using inequality~\eqref{1.ineg.var2}, we obtain
\begin{multline}\label{1.ineg.var}
-\int_0^{T}\!\!\int_0^{X_*}\dot{\phi}\left(\dfrac{w^2}{2}-\eta\, w\right)\, dx\,dt +\int_0^{T}\!\!\int_0^{X_*}\phi\, w\,\pt_t\eta\, dx\,dt +\int_0^{T}\!\!\int_0^{X_*}\phi\,\pt_x w\pt_x (w-\eta)\, dx\,dt\\
\le\phi(0)\int_0^{X_*}\left(\dfrac{w^2}{2}- w\eta\right)(x,0)\, dx +\int_0^{T}\!\!\int_0^{X_*}\phi g (w-\eta)\, dx\,dt.\phantom{aaa}
\end{multline}
On the one hand this computation is valid for $\rho$ such that $w=\chi\rho\in H^1((0,T);H^2(\R_+))$ and in this case~\eqref{1.ineg.var} and~\eqref{P.bord1w} imply~\eqref{P.bord2}. On the other hand~\eqref{1.ineg.var} has a meaning as soon as  $w\in L^2((0,T);H^1(0,X^*))$. In this sense,~\eqref{P.bord1w}\&\eqref{1.ineg.var} form a weak formulation of the boundary conditions~\eqref{P.bord}.

\begin{defini}
\label{def.sol.faible}
Let $\ov T\in(0+\oo]$, $X^0 > 0$ and $\rho^0\in L^2_{\rm loc}(\R_+)$. 
Then, we say that $(\rho, X)$ is a weak solution to~\eqref{P}--\eqref{P.bord} if the following conditions hold true.
\begin{enumerate}[(a)]
\item  $X: [0,\ov T)\to(0,+\oo)$ is Lipschitz continuous and $X(0)=X^0$.\medskip
\item Denoting by $\rho'$ the absolutely continuous part of $\pt_x\rho$, the function $\rho:\R_+\t(0,\ov T)\to\R_+$ satisfies:
\begin{itemize}
\item[$\circ$] $\rho\in L^\infty(\R_+\t(0,{\ov T}))$,\smallskip
\item[$\circ$] $\supp\rho=\ov{D_{\ov T}}$,\smallskip
\item[$\circ$] $\rho'\in L^2(D_T)$, for every $T\in(0,\ov T)$.\smallskip
%\footnote{\label{footnoterho'}Where $\rho'$ is the density of the absolutely continuous component of the measure $\partial_x \rho$, see Subsection~\ref{subsec.notation} below.}
\end{itemize}
\item Setting $M(t):=\MM(\rho(t))= \displaystyle \int_0^{X(t)}\rho(x,t)\,dx $ for $t\in(0,\ov T)$, there holds $M\in BV\loc([0,\ov T))$\footnote{We recall that possibly ${\ov T}= +\oo$.}.\smallskip
\item For all $\vhi\in\Co^\oo_c(\R_+\t[0,\ov T))$ 
\begin{multline}
\label{1.eq_limit_rho}
-\int_0^{\ov T}\int_{\R_+}\rho(x,t)\,\pt_t\vhi(x,t) dx\, dt -\int_{\R_+}\rho^0(x)\,\vhi(x,0)\, dx\\
-\langle\dot M,\vhi(0,\cdot)\rangle +\int_0^{\ov T}\int_{\R_+}\rho'(x,t)\,\pt_x\vhi(x,t)\, dx\, dt = 0,\phantom{aaa}
\end{multline}
where the measure $\dot{M}$ is defined in~\eqref{def:dotM} and $\rho'$ is defined at point~(b) above.
\item There holds
\be\label{1.eq_limit_X}
\lambda\,\dot{X}(t)  =\rho(X(t),t)-\ros\qquad\text{for almost every }t\in(0,\ov T).
\ee
\item The inequalities~\eqref{P.bord1w} hold true, that is, 
\[
\rom\le\rho(0,t)\le\rop\qquad\text{for almost every }t\in (0,\ov T).\smallskip
\]
\item Let  $T\in(0,\ov T)$, pick $0<X_*\le\min_{[0,T]}X$ and fix $\chi\in\Co^\oo_c(\R_+)$ satisfying~\eqref{X*chi}. Then, the variational inequality~\eqref{1.ineg.var} holds true for every $\phi\in\Co^\oo_c([0,T),\R_+)$ and  $\eta\in\Co^\oo_c(\R_+\t [0,T))$  with $\eta(0,t)\in [\rom,\rop]$ on $[0,T)$. 
\end{enumerate}
\end{defini}

We are now in position to state our main result.

\begin{theo}\label{th.main}
Let the following assumptions hold:
\begin{itemize}
\item[(H1)]Parameters of the model: let $\rop\ge\rom>0$, $\ros>0$ and $\lambda>0$.

\medskip 

\item[(H2)] Initial data: Let $X^0 > 0$ and $\rho^0\in L^1(\R_+,\R_+)$ such that 
\[
\rom\le\rho^0(0)\le\rop,\qquad\rho^0_{|[0,X^0]}\in \Co^{0,1}([0,X^0])\qquad\text{with}\qquad\rho^0(x) = 0\iff  x > X^0.
\] 
\end{itemize}
Then, there exists $\ov T\in(0,+\oo]$ and a weak solution $(\rho,X)$ to the system~\eqref{P}--\eqref{P.bord} on $(0,\ov T)$ in the sense of Definition~\ref{def.sol.faible}. Besides,
\begin{enumerate}[(a)]
\item If  $\ov T<+\oo$ then $X(t)\to 0$ as $t\up\ov T$.\smallskip
\item Let 
\be\label{rhominrhomax}
 0\ <\ \rob:=\min\lt\{\min_{[0,X^0]}\rho^0,\ros,\rop\rt\}\ \le\ \rot:=\max\lt\{\max_{[0,X^0]}\rho^0,\ros,\rom\rt\},
\ee
then we have $\rob\le\rho\le\rot$  in $D_{\ov T}$.
\item Recall the definition of the Gibbs energy $\PSI$ in~\eqref{PSI} and denote $\psi(t):=\PSI(\rho(t),X(t))$ for $t\in[0,\ov T)$ then $\psi\le\PSI(\rho^0,X^0)$ on $(0,T)$.
\end{enumerate}
\end{theo}

\subsection{Formal discussion about the possible disappearance of the oxide layer}\label{SubsectionSchrink}
The point~(a) states that in the case $\ov T<\oo$, the oxide domain disappears at time $\ov T$. 
Let us explain why we cannot prevent this from happening in general, at least when $\ros>\rop$. For this, we assume $\ros>\rop$ and  we build a ``super-solution'' $(\tilde\rho,\tilde X)$ to~\eqref{P}--\eqref{P.bord}.\\
Let us introduce a constant $\kappa>0$ to be determined later and let us  set 
\be\label{idrksupersolution}
\hat X:=\dfrac{\ros-\rop}{2\kappa}>0,\qquad\hat T:=\dfrac\lambda\kappa\log 2.
\ee
We then define for $0\le t\le\hat T$,
\[
\tilde X(t):=\hat X\lt(2-2^{t/\hat T}\rt),\qquad
\tilde\rho(x,t):=\begin{cases}\rop+\kappa x&\text{for }0\le x\le\tilde X(t),\\
\quad 0&\text{for }x>\tilde X(t).
\end{cases}
\]
We notice that $(\tilde\rho,\tilde X)$ satisfies~\eqref{P.a},~\eqref{P.e} and~\eqref{P.bord} with $\tilde X$ decreasing, $\tilde X(0)=\hat X$ and $\tilde X(\hat T)=0$ so that $\tilde X$ vanishes in finite time $\hat T$. However, $(\tilde\rho,\tilde X)$ fails to be a solution to~\eqref{P}--\eqref{P.bord} with respect to its own initial condition because the balance equation~\eqref{P.c} does not hold. Instead, using~\eqref{idrksupersolution} we compute for $0\le t\le\hat T$,
\begin{multline*}
 Q(t):=\pt_x\tilde\rho(\tilde X(t),t) +\dot{\tilde X}(t)\,\tilde\rho(\tilde X(t),t)\\
=\kappa\lt(1-\dfrac{\hat X}\lambda 2^{t/\hat T}\lt(\ros -\kappa\hat X2^{t/\hat T}\rt)\rt)
\ge\kappa\lt(1-\dfrac{2\hat X\ros}{\lambda}\rt)=\kappa-\dfrac{\ros(\ros-\rop)}\lambda.
\end{multline*}
Choosing $\kappa>0$ large enough, we have $Q>0$ on $(0,\hat T)$. This means that even substituting for~\eqref{P.c} the balance equation with positive source term
\[
\pt_x\tilde\rho(\tilde X(t),t) +\dot{\tilde X}(t)\,\tilde\rho(\tilde X(t),t)= Q(t)>0,
\]
the oxide layer disappears in finite time. Using an \textit{ad-hoc} version of the maximum principle\footnote{Assuming that there is some first time $t^*$ and a point $0\le x^*\le\tilde X(t^*)$, such that $\rho(x^*,t^*)=\tilde\rho(x^*,t^*)$ or $X(t^*)=\tilde X(t^*)$, we obtain a contradiction from~\eqref{P}--\eqref{P.bord} and $Q>0$.} we see that a smooth solution $(\rho,X)$ of~\eqref{P}--\eqref{P.bord} with an initial data $(\rho^0,X^0)$ such that $\rho^0\le\tilde\rho(\cdot,0)$ and $X^0<\tilde X(0)=\hat X$  satisfies $\rho<\tilde\rho$ and $X<\tilde X$ on its interval of existence. In particular $X$ vanishes in a finite time $\ov T<\hat T$.
 
When $\rop>\ros$, we expect that this phenomenon does not occur. Indeed, if $X(t)$ were decreasing to 0 as $t\up {\ov T}$, then we would have $\rho(t,X(t))<\ros<\rop$ and assuming enough regularity  we would also have $\rho(t,0)<\rop$ and the boundary condition~\eqref{P.bord} would lead to $\frac d{dt}\MM(\rho(t))=-\pt_x\rho(0,t)\ge0$  preventing the oxide domain from collapsing.

 To make the above arguments rigorous we should  study the regularity of the solutions. This goes beyond the scope of the paper.

\subsection{A few words about the proof of Theorem~\ref{th.main}}
The existence of a weak maximal solution to~\eqref{P}--\eqref{P.bord} is established by considering a discretization in time of the system through a minimizing movements scheme \textit{\`a la De Giorgi}. The problem which is solved at each discrete time is an optimization problem inspired by the JKO scheme.\footnote{The JKO scheme of~\cite{JKO98}  produces approximate solutions of the Fokker--Planck equation in the whole space.} Here the minimization drives the evolution of the oxygen density in the bulk but also the position of the moving interface $X$ and the Signorini boundary conditions at $x=0$. This is reflected, on the one hand, in the driving energy (the same Gibbs free energy introduced above with an additional  term $p_\tau$ penalizing large fluxes at $x=0$) and on the other hand, in the underlying metric (which contains in particular a Wassertein diistance between measure with possibly different masses and supports).

 The energy, the metric and the basic optimization step are introduced in Section~\ref{sec.intmin}. The rest of this section is devoted to the properties of the  minimizers for one step of the scheme.%: existence, \textit{a priori} estimates and Euler Lagrange equations which turn out to be a discrete version of~\eqref{P}--\eqref{P.bord}.  The compactness of the family of discrete solutions and the passage to the limit in the weak formulations of~\eqref{P}--\eqref{P.bord} are achieved in Section~\ref{Sec.Th2}.  
 
 The idea of using a minimizing movement scheme to treat moving boundaries is taken from~\cite{PP10}. However, to take into account the double Signorini boundary condition at $x=0$ we are led to introduce a linear term in the metric (see~\eqref{2.def_d1}) which breaks the nice Riemannian-like structure. The most difficult part of the proof of Theorem~\ref{th.main} is to establish that the candidate solution obtained in the limit by sending the time step to 0 satisfies the weak form~\eqref{1.ineg.var} of this condition. Indeed, the Signorini inequalities involve the trace of $\pt_x\rho$  for which we do not even have uniform bounds. To improve the regularity of the solution we add a  penalty function $p_\tau$ to the driving energy.  This allows us to obtain a bound on $\Lip(\rho)$ which degenerates as the time step goes to 0 but which is sufficient, when combined to mollification techniques, to get~\eqref{1.ineg.var} at the limit.

\subsection{Notation and conventions}\label{subsec.notation}

We denote by $\Leb$ the Lebesgue measure on $\R$ and we use standard notation for Lebesgue and Sobolev spaces. For any of these spaces, the corresponding subspace of compactly supported functions is denoted with a subscript $c$. For instance $L^1_c(\Om)$ is the space of compactly supported integrable function on $\Om$.

When $\nu$ is a measure on $\mathcal{X}$ and $K\sub\mathcal{X}$ is a $\nu$-measurable set, we denote $\nu\restr K$, the restriction of $\nu$ to $K$, that is $\nu\restr K(\Om):=\nu(K\cap\Om)$.

 The derivatives of locally integrable functions of one or several variables are meant in the sense of distributions whenever classical derivatives do not exist.  The partial derivatives of a function of two variables $(x,t)\mapsto w(x,t)$ are denoted $\pt_xw(x,t)$ and $\pt_tw(x,t)$. We also often use $\dot w$ the time derivatives.   When the distribution $\pt_x w$ is a measure, we denote $w'$ the absolute continuous part of  $\pt_x w$, that is, as a measure, 
\[%\be\label{prime}
\pt_x w = w'\Leb +\nu\qquad\text{ where }\nu=\nu\restr K\text{ for some set }K\text{ with }\Leb(K)=0.  
\]%\ee
For a function $w$ defined on a set $I$ and $J\sub I$, we denote $w_{|J}$ the restriction of $w$ on $J$.

 When $I$ is an interval, the space of functions with bounded variations on $I$ is denoted by $BV(I)$.  We denote $\Co^{0,1}(I)$ the space of Lipschitz continuous functions on $I$. More generally for a positive integer $k$, $\Co^{k,1}$ denotes the spaces of functions $w:I\to\R$ of class $\Co^k$ on $I$ with $w^{(k)}\in \Co^{0,1}(I)$. The Lipschitz constant of $w\in \Co^{0,1}(I)$ is denoted 
\[
\Lip(w):=\sup_{x,y\in I,\ x\ne y,}\dfrac{|w(x)-w(y)|}{|x-y|}.
\]

We deal with functions $\rho\in L^1(\R_+)$ such that, for  some $X>0$, $\rho(x)= 0$ for almost every $x>X$ and $\rho_{|[0,X]}$ admits a continuous representative.  We always consider the representative which is continuous on $[0,X]$ and vanishes on $(X,+\oo)$ and we also denote it $\rho$. This function is defined  everywhere and is left-continuous.

For functions which depend on space and time such as $\rho\in L^1(\R_+\t (0,T))$, we denote $\rho(t):=\rho(\cdot,t)\in L^1(\R_+)$. The function $\rho(t)$ is well defined for almost every $t\in (0,T)$. Notice that we already used this convention, for instance when defining $M(t)=\MM(\rho(t))$. 

We denote the positive part of $x\in\R$ by 
\[
x_+:=\max(x,0).
\]

Let $\mathcal{X}$ and $\mathcal{Y}$ be complete metric spaces and let $\mu$ be a finite Radon measure on $\mathcal{X}$. For $T\in\Co(X,Y)$, we denote $T\pf\mu=\nu$ the pushforward of $\mu$ by $T$. This is the Radon measure  on $\mathcal{Y}$ characterized  by 
\be\label{pfdef}
\int_\mathcal{Y}\psi(y)\, d\nu(y) =\int_\mathcal{X}\psi(T(x))\,d\mu(x)\qquad\text{for every }\psi\in\Co(\mathcal{Y}).
\ee

Finally, we draw the reader's attention to the following inconsistency in the use of the letter $T$. 
\begin{enumerate}[($*$)]
\item In Section~\ref{sec.intmin}, the notation $T$ and its variants $T_+$, $T_-$,\dots denote transport maps.
\item In Section~\ref{Sec.Th2}, $T$, $\ov T$, $\mathcal{T}$ denote different positive times.
\end{enumerate}
This should not be confusing since the continuous time variable does not appear in Section~\ref{sec.intmin} and we do not make explicit use of transport maps in Section~\ref{Sec.Th2}.

\section{The optimisation problem defining one time step of the semidiscrete scheme}\label{sec.intmin}

In this section, we first recall the definition of the 2-Wasserstein distance and some of its basic properties. Then, in Subsection~\ref{sec.not}, we introduce a 2-Wasserstein like distance for specific measures with unequal masses. In Subsection~\ref{subsec_OTP} we describe the optimal transport plans between such measures  Subsection~\ref{sec.min} is devoted to the definition of one step of the minimizing movement scheme. We prove the existence of at least one solution to this optimisation problem in Subsection~\ref{subs.existenceOP} and we show in Subsection~\ref{subs:onestep_timediscretization} that this problem can be seen as a time discretisation of~\eqref{P}. In Subsection~\ref{Sec:comp_inter_fixe}, we study the behavior of the solutions at the fixed interface $x=0$. The remaining of Section~\ref{sec.intmin} is dedicated to the proof of $L^\infty$ and Lipchitz bounds on the solutions to the optimisation problem. For the latter, we use smooth approximations of the initial data $\rho^0$ on $[0,X^0]$. We justify this approximation procedure in Subsection~\ref{subsec.M(rho)prescibed}. The $L^\oo$  and Lipschitz estimates are then established in Subsection~\ref{subsec.LooLipbounds}.\medskip

\subsection{The 2-Wasserstein distance in dimension 1}~
%\label{sec.dist}

 We use basic tools from optimal transport theory in dimension 1 that are described for instance in~\cite{OTAM}.  
Let $I=[a,b]$ and $J=[c,d]$ be two bounded intervals and, for $K=I$ or $K=J$, let $\mathcal{M}_+(K)$  be the set of finite positive Radon measures on $X$. Given $\mu\in\mathcal{M}_+(I)$ and $\nu\in\mathcal{M}_+(J)$ with $\mu(I)=\nu(J) =m$ for some $m > 0$, the squared Wasserstein distance between $\mu$ and $\nu$ is defined as
\be\label{2.def.wass}
\wass^2(\mu,\nu):=\inf_{\gamma\in\PI(\mu,\nu)}\int_{I\t J}(x-y)^2\, d\gamma(x,y),
\ee
where $\PI(\mu,\nu)$ denotes the set of transport plans from $\mu$ to $\nu$, that is,
\[
\PI(\mu,\nu):=\ds\lt\lbrace\gamma\in\mathcal{M}_+(I\t J)\,\,:\,\,\gamma(I\t J) = m,\,\,\pi_x\pf\gamma =\mu,\,\,\pi_y\pf\gamma =\nu\rt\rbrace.
\]
Here $\pi_x$ and $\pi_y$ denote the projections into the first and second coordinates ($\pi_x(x,y)=x,\ \pi_y(x,y)=y$ and (recall~\eqref{pfdef})  the above pushforwards of $\gamma$  by these mappings are the Radon measures defined for $\vhi\in\Co(I)$, $\psi\in\Co(J)$ by
\[
\int_I\vhi\, d(\pi_x\pf\gamma)=\int_{I\times J}\vhi(\pi_x(x,y)) d\gamma(x,y),\qquad\int_J\psi\, d(\pi_y\pf\gamma)=\int_{I\times J}\psi(\pi_y(x,y)) d\gamma(x,y).
\]
For the sake of completeness, let us state below some classical results in optimal transport theory.

\begin{theo}{\cite[Chapter 1]{OTAM}}\label{th.rappel}
 Let $\mu$, $\nu$ as above. 
\begin{enumerate}[(a)]
\item The optimization problem~\eqref{2.def.wass} admits a dual formulation:
\[
\wass^2(\mu,\nu) =\sup_{(\vhi,\psi)\in\mathcal{S}}\int_I\vhi\, d\mu +\int_J\psi\,d\nu,
\]
 where 
\[
\mathcal{S}=\lt\{(\vhi,\psi)\in\Co(\R)\t\Co(\R):\vhi(x)+\psi(y)\le (x-y)^2\text{ for every }(x,y)\in\R^2\rt\}.
\]
 
\item There exists a unique optimal transport plan $\gamma\in\PI(\mu,\nu)$ associated with the optimization problem~\eqref{2.def.wass}.
 This is the unique monotonic transport plan of $\PI(\mu,\nu)$ in the sense that $(x_2-x_1)(y_2-y_1)\ge 0$ for any  $(x_1,y_1),(x_2,y_2)\in\supp\gamma$. 
\smallskip
\item If moreover, $\mu$ is atomless then this optimal transport plan $\gamma$ is induced by a nondecreasing map $T$ such that $\gamma = (\mathrm{id},T)\pf\mu$. Recalling the notation $I=[a,b]$, $J=[c,d]$, we have in particular, 
\[
\mu([a,x))=\nu([c,T(x))\qquad\text{ for }x\in (a,b]. 
\]
In this case, there exists a unique (up to an additive constant) Lipschitz function $\Phi$, called Kantorovich potential, such that there holds
\[
\Phi'(x) = x - T(x)\quad\mbox{for almost every }x\in I.
\]
As a consequence,
\be\label{Lip(Phi)}
\Lip(\Phi)\le 2\max_{x\in I,\, y\in J}|x-y| =2\max(|d-a|,|b-c|).
\ee
\item Assuming that $\mu=p\Leb$ with $p\le C$ and that $\nu$ decomposes as $\nu=q\Leb+\tilde\nu$ with $\tilde\nu\perp\Leb$ and $q\ge\kappa>0$,  then $T$ is Lipschitz continuous with $0\le T'\le C/\kappa$ and there holds
\[
p(x)=T'(x) q(T(x))\qquad\text{ for almost every }x\in(a,b). 
\]  
\end{enumerate}
\end{theo}

\subsection{An unbalanced optimal transport problem}\label{sec.not}~

Let $\rho^0,\rho\in L^1(\R_+,\R_+)$ be compactly supported. We wish to define a Wasserstein-like distance between $\rho^0$ and $\rho$. As we do not assume that $\MM(\rho^0)=\MM(\rho)$, we cannot transport the measure $\rho^0\Leb$ on $\rho\Leb$ and
\[
\wass(\rho^0\Leb,\rho\Leb)=+\oo\qquad\text{whenever }\MM(\rho^0)\ne\MM(\rho).
\]  
Here, we start by balancing the masses by adding a Dirac measure at $x=0$.  We set:
\be
\label{2.mu1}
\mu_{\rho^0,\rho}:=\rho^0\Leb\restr\R_+ + (\MM(\rho)-\MM(\rho^0))_+\,\delta_0,
\ee
where $x_+:=\max(x,0)$ and $\delta_0$ denotes the Dirac measure at  $x=0$. 
Symmetrically, we set
\[
\mu_{\rho,\rho^0}=\rho\Leb\restr\R_+ + (\MM(\rho^0)-\MM(\rho))_+\,\delta_0.
\]
For shortness, we write
\be
\label{3.mu1}
\mu^0:=\mu_{\rho^0,\rho} ,\quad\qquad\mu:=\mu_{\rho,\rho^0}.
\ee
These measures are compactly supported and satisfy $\mu^0(\R_+)=\mu(\R_+)$. We  define 
\be\label{def_W2}
\Wass(\rho^0,\rho):=\wass\lt(\mu^0,\mu\rt).
\ee
We see that $\Wass(\rho^0,\rho)<\oo$ for any pair of compactly supported functions  $\rho^0,\rho\in L^1(\R_+,\R_+)$.
\begin{Remark}\label{rem_dualWass}
Using the dual formulation in Theorem~\ref{th.rappel} (a), we see that 
\[
\Wass^2 (\rho^0,\rho)=\sup_{(\vhi,\psi)\in\mathcal{S}}B((\rho^0,\rho),(\vhi,\psi)),
\]
where we introduce the bilinear functional,
\[
 B((\rho^0,\rho),(\vhi,\psi)):=\int_{\R_+}\rho^0\vhi\,+\int_{\R_+}\rho\psi\, +\lt[\MM(\rho)-\MM(\rho^0)\rt]_+\vhi(0)-\lt[\MM(\rho^0)-\MM(\rho)\rt]_+\psi(0).
\]
Let us consider two sequences $(\rho^0_j)$ and $(\rho_j)$ of nonnegative functions in $L^1(\R_+)$ all supported in $[0,\Lambda]$ for some $\Lambda>0$ and such that,
\[
\rho^0_j\longto\rho^0,\qquad\rho_j\longto\rho\qquad\text{weakly-$\star$ in }\mathcal{M}(\R_+)\text{ as }j\up\oo,
\]
that is, 
\[
\quad\int_{\R_+}\rho^0_j\phi\,\ \st{j\up\oo}\longto\ \int_{\R_+}\rho^0\phi\,,\qquad\int_{\R_+}\rho_j\phi\,\ \st{j\up\oo}\longto\ \int_{\R_+}\rho\phi\,\qquad\text{for every }\phi\in\Co(\R).
\]
In particular, by taking  $\phi\equiv 1$, we have  $\MM(\rho^0_j)\to\MM(\rho^0)$ and $\MM(\rho_j)\to\MM(\rho)$.\\
For such sequences and for $(\vhi,\psi)\in\mathcal{S}$, there holds
\[
\lim_{j\up\oo}B((\rho^0_j,\rho_j),(\vhi,\psi)) =  B((\rho^0,\rho),(\vhi,\psi)). 
\]
Taking the supremum over $(\vhi,\psi)\in\mathcal{S}$ we obtain 
\[
\Wass(\rho^0,\rho)\le\lim_{j\up\oo}\Wass(\rho^0_j,\rho_j).
\]
In other words, the distance $\Wass$ is lower semicontinuous with respect to the weak-$\star$  convergence in $L^1(\R_+)$ in the subset of nonnegative and uniformly compactly supported functions.
\end{Remark}

\subsection{Optimal transport plans}

\label{subsec_OTP}
\begin{figure}[ht]
\centering
\begin{minipage}{0.48\textwidth}
\centering
\captionsetup{width=.9\textwidth}
\includegraphics[]{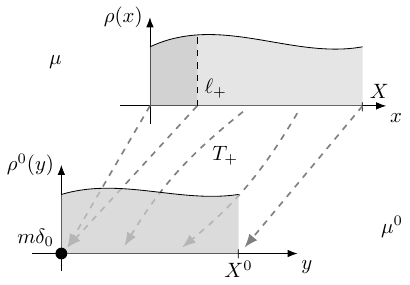}
\caption{The transport map $T_+$ in the case $0\le\MM(\rho)-\MM(\rho^0)=:m$.\label{figT+}}
\end{minipage}\hfill
\begin{minipage}{0.48\textwidth}
\centering
\captionsetup{width=.9\textwidth}
\includegraphics[]{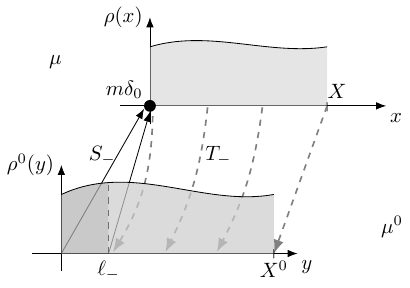}
\caption{The transport maps $T_-$ and $S_-$ in the case $0<\MM(\rho^0)-\MM(\rho)=:m$.\label{figT-}}
\end{minipage}

	\vspace{4pt}
\end{figure}

~\\
Let us now describe the optimal transport plan from $\mu$ to  $\mu^0$. For this we denote 
\[
X^0:=\inf\{x>0:\rho^0\equiv 0\text{ on }(x,+\oo)\},\qquad X:=\inf\{x>0:\rho\equiv 0\text{ on }(x,+\oo)\}.
\]
According to Theorem~\ref{th.rappel}, the measure $\mu$ is sent monotonically on the measure $\mu^0$. We have two situations, depending on the sign of $\MM(\rho)-\MM(\rho^0)$, that we describe below. The corresponding constructions are also illustrated in Figures~\ref{figT+} and~\ref{figT-}. 

\smallskip
\noindent\textbf{Case $\MM(\rho)\ge\MM(\rho^0)$}.  Let us denote $\ell_+\ge0$ the smallest number such that 
\be\label{ell+}
\int_0^{\ell_+}\rho\, =\MM(\rho)-\MM(\rho^0)=:m.
\ee
The measure $\rho\Leb\restr (0,\ell_+)$ is sent to the Dirac mass $m\delta_0$ and the remaining part $\rho\Leb\restr (\ell_+,X)$ is sent monotonically on $\rho^0\Leb\restr (0,X^0)$. More precisely, we set $T_+(x)=0$ for $0\le x\le\ell_+$ and for $\ell_+\le x\le  X$, we define  $T_+(x)\ge0$ as the smallest number such that 
\be\label{T+}
\int_{\ell_+}^x\rho\,=\int_0^{T_+(x)}\rho^0\,.
\ee 
By construction, $T_+(X)=X^0$ and $T_+\pf\rho\,\Leb\restr (\ell_+,X) =\rho^0\,\Leb\restr (0,X^0)$. The optimal transport plan $\gamma_+\in\PI(\mu,\mu^0)$ is given by
\[
\gamma_+:=(\Id, T_+)\pf (\rho\,\Leb\restr (0,X)).
\]
The  corresponding Wasserstein energy is 
\[
\Wass^2 (\rho,\rho^0)  =\int_0^{\ell_+}x^2\,\rho(x)\, dx +\int_{\ell_+}^X ( x -T_+(x))^2\,\rho(x)\, dx =\int_{0}^X ( x -T_+(x))^2\,\rho(x)\, dx.
\]

\medskip

\noindent\textbf{Case $\MM(\rho)<\MM(\rho^0)$.} In this case we define  $\ell_- > 0$ by
\be\label{ell-}
\int_0^{\ell_-}\rho^0\, =\MM(\rho^0)-\MM(\rho),
\ee
and  we define  $T_-(x)\ge\ell_-$ for $x\in [0,X]$ as the smallest number such that 
\be\label{T-}
\int_0^x\rho\,=\int_{\ell_-}^{T_-(x)}\rho^0.
\ee
We have by construction, $T_-(X)=X^0$ and
\[
T_-(0)=\ell_-,\qquad T_-\pf\rho\Leb\restr(0,X) =\rho^0\Leb\restr (\ell_-,X^0).
\]
The unique optimal transport plan $\gamma_-\in\PI(\mu,\mu^0)$ is 
\[
\gamma_-:= (S_-,\Id)\pf\rho^0\,\Leb\restr (0,\ell_-)+ (\Id, T_-)\pf\rho\,\Leb\restr\R_+, 
\]
where  $S_-$ is the constant map $S_-\equiv 0$ on $[0,\ell_-]$. In other words, the ``misssing mass'' $\MM(\rho)-\MM(\rho^0)$ is added to $\rho\Leb$ as a Dirac mass at $0$ and sent to $\rho^0\Leb\restr(0,\ell^-)$. The measure $\rho\Leb$ is then distributed monotonically on $\rho^0\Leb\restr (\ell_-,X^0)$. The Wasserstein energy is
\[
\Wass^2(\rho,\rho^0) =\int_0^{X}(x-T_-(x))^2\,\rho(x)\, dx +\int_0^{\ell_-}y^2\,\rho^0(y)\, dy.\medskip
\]

In the sequel, we set by convention 
\be\label{ell+ell-0}
\begin{cases}
\ell_-=0&\text{ when }\MM(\rho)\ge\MM(\rho^0),\\
\ell_+=0&\text{ when }\MM(\rho)\le\MM(\rho^0).
\end{cases}
\ee
\medskip

\subsection{Definition of one step of the minimizing movement scheme}~

\label{sec.min}

In the rest of Section~\ref{sec.intmin}, the time step $\tau>0$ is fixed as well with the coefficients $\rot\ge\rom>0$, $\ros>0$ and $\lambda>0$ as in assumption~(H1) of Theorem~\ref{th.main}.

Let us first introduce the set of candidates for the minimization problem. 
\be
\label{SpaceAbb}
\Abb:=\lt\{ (\rho,X)\in L^1_c(\R_+)\t (0,+\oo):\rho\ge0\ \text{ and }\ \rho\equiv 0\text{ on }(X,+\oo)\rt\}.
\ee
For $\rho\in L_c^1(\R_+)$, we denote
\be\label{xrho}
x_\rho:=\inf\lt\{x\ge 0:\rho\equiv 0\text{ almost everywhere on }(x,+\oo)\rt\}=\inf\lt\{x\ge 0:\supp\rho\sub [0,x]\rt\}.
\ee
We then define the following subset of  $\Abb$,
\be\label{SpaceA}
\A:=\lt\{ (\rho,X)\in\Abb: X=x_\rho,\ \rho_{|[0,X]}\in \Co^{0,1}([0,X]),\ \rho>0\text{ on }[0,X]\rt\}. 
\ee
Remark that for $(\rho,X)\in\A$, $\rho$ does not vanish identically (indeed $x_\rho=X>0$). Notice also that $X$ is uniquely determined by $\rho$ so that we could get rid of the explicit mention of $X$. However it is convenient to see $\A$ as a subset of $\Abb$.\medskip

\noindent{\it Energy functional.} We repeat the definition~\eqref{PSI} of the functional $\PSI:\Abb\to\R\cup\{+\oo\}$:
\[
\PSI(\rho,X):=\int_{\R_+}f(\rho(x))\, dx  +\ros X, 
\]
where (recall~\eqref{f}), $f(r)=r(\log r +\beta -1)$.\medskip

\noindent{\it Distance on $L^1_c(\R_+,\R_+)\t\R_+$.} For every $(\rho,X)$, $(\rho^0,X^0)\in L^1_c(\R_+,\R_+)\t\R_+$ we define the following ``distances'' $\textbf{d}_1$, $\textbf{d}_2$ on $\Abb$,
\begin{align}
\label{2.def_d1}
\textbf{d}_1\lt((\rho, X),(\rho^0,X^0)\rt)&:=\te\left|\MM(\rho)-\MM(\rho^0)\right|,\\
\label{2.def_d2}
\textbf{d}_2^{\,2}\lt((\rho, X),(\rho^0,X^0)\rt)&:=\Wass^2\left(\rho,\rho^0\right) +\lambda\left(X-X^0\right)^2.\medskip
\end{align}
We notice that $\textbf{d}_1$ is built on a $L^1$-distance and that $\textbf{d}_2$ is modelled on a $L^2$-distance (the Wasserstein distance $\wass$ is thought of as a Riemannian metric on sets of positive measures).\smallskip

\noindent{\it Penalty term.} Let $\om\in(0,1)$ to be fixed later. We define $p_\tau:\R\to\R_+$  by
\be\label{2.def.penalty}
p_\tau(m):=\dfrac1{2\tau^{1-\om}}\lt(|m|-\tau^{1-\om}\rt)_{\!+}^{\,2}\qquad\text{for }m\in\R,
\ee
where we recall that $\tau>0$ is the time step. The function $p_\tau$ is used as a technical penalization term to ensure some Lipschitz bounds on the solutions of the minimization problem~\eqref{OP} below (see Proposition~\ref{prop.Lipbounds}). The exponent $\om\in(0,1)$ will be fixed at the end of the proof (see the condition~\eqref{cond.vartheta} of Proposition~\ref{prop.inegvar}).\medskip

\noindent{\it Minimization problem.} We define the functional $\J^\tau_{(\rho^0,X^0)}:\,\Abb\to\R\cup\lt\{ +\infty\rt\}$ by
\begin{multline}
\label{Jtau}
\J^\tau_{(\rho^0,X^0)}(\rho, X) 
:=\dfrac1{2\tau}\textbf{d}
_2^{\,2}\lt((\rho,X),(\rho^0,X^0)\rt)
 +\textbf{d}_1\lt((\rho,X),(\rho^0,X^0)\rt)\\
+ p_\tau\left(\MM(\rho)-\MM(\rho^0)\right)+\PSI(\rho,X).
\end{multline}
Expanding all terms, the functional decomposes into a bulk term (depending on the function $\rho$), a left boundary term (depending on $\MM(\rho)$) and a right boundary term (depending on $X$).
\begin{align}
\nonumber
\J^\tau_{(\rho^0,X^0)}(\rho, X) =&\dfrac1{2\tau}\Wass^2\left(\rho,\rho^0\right)+\int_{\R_+}f(\rho(x))\, dx\\
\label{Jtau2}
&+\te\left|\MM(\rho)-\MM(\rho^0)\right|+ p_\tau\left(\MM(\rho)-\MM(\rho^0)\right)\\
\nonumber
&+\dfrac1{2\tau}\lambda\left(X-X^0\right)^2+\ros X.
\end{align}
Eventually, we introduce the optimization problem 
\be
\label{OP}
\min_{(\rho,X)\in\Abb}\,\J^\tau_{(\rho^0,X^0)}(\rho,X).
\ee
\medskip

\subsection{Existence of a minimizer and first properties}~
\label{subs.existenceOP}

Throughout the paper, we assume that $(\rho^0,X^0)\in\A$ and  $\rom\le\rho^0(0)\le\rop$ so that the assumption~(H2) of Theorem~\ref{th.main} holds true. Using the Direct Method of the Calculus of Variation, we establish in Theorem~\ref{theo.ExistenceJKO.1} that~\eqref{OP} admits at least one solution when the time step $\tau>0$ is small enough. Then, in Theorem~\ref{theo.ExistenceJKO.2} we establish that such minimizer $(\rho,X)$ belongs to $\A$ with $\log\rho$ of class $\Co^{1,1}$ on $[0,X]$, we describe the boundary condition satisfied by $\rho$ at $x=0$ and we give a bound on $|\MM(\rho)-\MM(\rho^0)|$. 

\begin{theo}
\label{theo.ExistenceJKO.1}
Let us denote 
\[
R:=\ros-\sqrt{\rop\rom}\qquad\text{and}\qquad P(\rho^0,X^0):=\int_0^{X^0}f(\rho^0)\,+\sqrt{\rop\rom}\,\,X^0=\PSI(\rho^0,X^0)-RX^0.
\]
Provided that 
\be\label{X0large}
 X^0>\dfrac{\tau R}\lambda +\sqrt{\dfrac{\tau^2 R^2}{\lambda^2}+\dfrac{2\tau}\lambda P(\rho^0,X^0)}\,,
\ee
the minimization problem~\eqref{OP} admits at least one solution $(\rho,X)\in\Abb$. 
\end{theo}

\begin{Remark}\label{RemX0large}~
Since $X^0>0$, the condition~\eqref{X0large} holds true for $\tau>0$ small enough (depending on $(\rho^0,X^0)$).\medskip
\end{Remark}

\begin{proof}[Proof of Theorem~\ref{theo.ExistenceJKO.1}]~
Let $(\rho_j,X_j)\in\Abb$ be a minimizing sequence for~\eqref{OP} indexed by $j\ge1$. Considering the candidate $(\rho,X)=(\rho^0,X^0)$, we may assume that 
\[
\J^\tau_{(\rho^0,X^0)}(\rho_j,X_j)\le\J^\tau_{(\rho^0,X^0)}(\rho^0,X^0) =\PSI(\rho^0,X^0).
\]
 In particular, 
\[
\PSI(\rho_j,X_j) +\dfrac\lambda{2\tau}(X^0-X_j)^2\le\PSI(\rho^0,X^0)\qquad\text{for }j\ge1.
\]
Using the definition~\eqref{PSI} of $\PSI$ and rearranging the terms we get
\[
\dfrac\lambda{2\tau}(X^0-X_j)^2 +\ros (X_j-X^0)\le\int_0^{X^0}f(\rho^0)\,  -\int_0^{X_j}f(\rho_j)\,.
\]
We then use $\min_{\R_+}f=-\sqrt{\rop\rom}$ to obtain 
\[
\dfrac\lambda{2\tau}(X^0-X_j)^2 +(\ros-\sqrt{\rop\rom}) (X_j-X^0)\le\int_0^{X^0}f(\rho^0)\, +\sqrt{\rop\rom}X^0.\\
\]
Denoting $x:=X_j-X^0$ and defining, as in the statement of the proposition,
\[
R:=\ros-\sqrt{\rop\rom}\quad\text{ and }\quad  P:=P(\rho^0,X^0):=\int_0^{X^0}f(\rho^0)\,+\sqrt{\rop\rom}X^0,
\]
the last inequality rewrites as
\be\label{proofprop.ExistenceJKO_0}
x^2 +2\dfrac{\tau R}\lambda x-\dfrac{2\tau P}\lambda\le0. 
\ee
Notice that 
\[
P=\int_0^{X^0}f(\rho^0)\, +\sqrt{\rop\rom}\,X^0\ge (\min f +\sqrt{\rop\rom}) X^0=0.
\]
Computing the roots of the polynomial function of $x$ in the left hand side of~\eqref{proofprop.ExistenceJKO_0}, we  have
\be\label{proofprop.ExistenceJKO_1}
X_j=X^0+x\ge X^0 -\dfrac{\tau R}\lambda  -\sqrt{\dfrac{\tau^2R^2}{\lambda^2}+\dfrac{2\tau P}\lambda}\,\st{\eqref{X0large}}\ge\,\eps,
\ee
for some $\eps>0$. 

Now, we write
\[
\PSI(\rho^0,X^0)\ge\J^\tau_{(\rho^0,X^0)}(\rho_j,X_j)
\ge\int_0^{X_j}f(\rho_j)\, +\ros X_j+\dfrac{\lambda(X^0-X_j)^2}{2\tau}.
\]
By convexity, of $f$, the first term of the right hand side satisfies
\[
\int_0^{X_j}f(\rho_j)\,\ge X_j f(\MM(\rho_j)/X_j),
\]
and writing $Y_j:=\MM(\rho_j)/X_j$, we get 
\be\label{ineqXjYj}
\PSI(\rho^0,X^0)\ge f(Y_j)  X_j +\ros X_j  +\dfrac{\lambda(X^0-X_j)^2}{2\tau}.
\ee
The term $\ros X_j  +\lambda(X^0-X_j)^2/(2\tau)$ is bounded from below and since $X_j\ge\eps>0$ by~\eqref{proofprop.ExistenceJKO_1} and $\lim_{+\oo}f=+\oo$, \eqref{ineqXjYj} implies that the sequence $(Y_j)$ is bounded. By continuity of $f$, the sequence $(f(Y_j))$ is  bounded and  using  again~\eqref{ineqXjYj} with this information, we deduce that $X_j$ is bounded. Now, as $\MM(\rho_j)=Y_j/X_j\le Y_j/\eps$, we get that $\MM(\rho_j))=\|\rho_j\|_{L^1}$ is bounded.\\
 Consequently, there exists $X\ge\eps$ such that, up to extraction, $X_j\to X$ as $j\up\oo$. Let us fix $\Lambda>0$ such that $X_j\le\Lambda$  for $j\ge1$ (hence $\rho_j\equiv 0$ on $(X_j,+\oo)$).  We see that the sequence  $(\int_0^\Lambda\rho_j\log\rho_j\,)$ is uniformly bounded  and we deduce that the sequence $(\rho_j)$ is equi-integrable. By Dunford--Pettis theorem, there exists $\rho\in L^1_c(\R_+)$ with $(\rho,X)\in\Abb$ such that, up to extraction, $\rho_j\to\rho$ weakly in $L^1(\R_+)$. By Remark~\ref{rem_dualWass}, we have 
\[
\MM(\rho_j)\st{j\up\oo}\longto\ \MM(\rho)\qquad\text{and}\qquad\Wass(\rho^0,\rho)\le\liminf_{j\up\oo}\Wass(\rho^0,\rho_j).
\]
Eventually, by convexity and lower semicontinuity of $f$ the weak convergence $\rho_j\to\rho$  yields
\[
\int f(\rho)\,\le\liminf_{j\up\oo}\int f(\rho_j)\,.  
\]
We conclude that  
\[
\J^\tau_{(\rho^0,X^0)}(\rho,X)\le\liminf_{j\up\oo}\J^\tau_{(\rho^0,X^0)}(\rho_j,X_j),
\]
so that $(\rho,X)\in\Abb$ is a minimizer of~\eqref{OP}.
\end{proof}

\begin{theo}[Existence of a minimizer in $\A$]
\label{theo.ExistenceJKO.2}
Let $(\rho,X)\in\Abb$ be a minimizer of~\eqref{OP}.  We have $(\rho,X)\in\A$ with $\rho$ of class $\Co^{1,1}$ on $[0,X]$ and of class $\Co^\oo$ on $[0,\ell_+]$ (if $\ell_+>0$).\\ 
Moreover, the transport map $T$ (defined as in~\eqref{notationT} below) is of class $\Co^{0,1}$ on $[0,X]$ (with $T\equiv 0$ in $[\ell_+,X]$).
\end{theo}

\begin{proof}Let $(\rho,X)$ be a minimizer of~\eqref{OP}. We start by noticing that $X>0$ being fixed, $\rho$ solves a one step of a classical JKO scheme. Let us be more specific. With the convention~\eqref{ell+ell-0} on $\ell_-$ and $\ell_+$, we denote 
\begin{align*}
\tilde\mu^0 &:=\rho^0\Leb\restr (\ell_-,X^0) + [\MM(\rho)-\MM(\rho^0)]_+\delta_0\\
 &=
\begin{cases}
\rho^0\Leb\restr (0,X^0) + [\MM(\rho)-\MM(\rho^0)]\delta_0&\text{ if }\MM(\rho)\ge\MM(\rho^0),\\
\rho^0\Leb\restr (\ell_-,X^0) &\text{ if }\MM(\rho)<\MM(\rho^0).
\end{cases}
\end{align*}
By construction $\supp\tilde\mu^0=[\ell_-,X^0]$, $\tilde\mu^0(\R_+)=\MM(\rho)$ and $\rho\Leb\restr(0,X)$  minimizes the functional 
\be
\label{JKO(nu)}
\nu\mapsto\dfrac1{2\tau}\wass^2(\tilde\mu^0 ,\nu) +\int_0^X f(v)\,,
\ee
in the set of measures
\[
\mathcal{M}_\rho=\lt\{\nu\in\mathcal{M}_+([0,X]):\nu([0,X])=\MM(\rho),\ \nu = v\Leb\,\text{ for some }v\in L^1(\R_+)\rt\}.
\]
In other words, $\rho$ is the solution of one step of a classical JKO scheme in $[0,X]$ with initial data $\tilde\mu^0$ and time step $\tau>0$. This is a well known situation, for which we refer to~\cite[Chapter~7.4.1]{OTAM}. Recaling the notation $T_+$ or $T_-$ of~\eqref{T+},~\eqref{T-} for the optimal transport map we set, for $x\in[0,X]$, 
\be\label{notationT}
T(x):=\begin{cases}
T_+(x)&\text{ if }\MM(\rho)\ge\MM(\rho^0),\\
T_-(x)&\text{ if }\MM(\rho)<\MM(\rho^0).
\end{cases}
\ee
In both cases, $\tilde\mu^0 = T\pf\rho$. 
The function $\rho$ being a minimizer of~\eqref{JKO(nu)} in $\mathcal{M}_\rho$, we conclude from~\cite[Chapter~7.4.1]{OTAM}, that this minimizer is unique and Lipschitz continuous. Besides, using the notation of Theorem~\ref{th.rappel}, the Kantorovitch potential $\Phi$ associated with the transport $T$ is Lipschitz continuous  and we have the relation
\be\label{Phi/tau+logrho}
\dfrac{\Phi(x)}{\tau}+\log\rho(x) = C\qquad\mbox{ for }x\in [0,X],
\ee
for some constant $C\in\R$. In particular, $\rho>0$ in $[0,X]$ which proves that $x_\rho=X$. We conclude that $(\rho,X)\in\A$.

Next, we note from point~(d) of Theorem~\ref{th.rappel}, we have 
\be\label{T'(x)=}
T'(x)=\rho(x)/\rho^0(T(x)),
\ee
for almost every $x\in [\ell_+,X]$ and from point~(c) of the same theorem we have  
\be\label{PhiT}
\Phi'(x)= x-T(x)\qquad\mbox{for }x\in [0,X].
\ee
Since $T'$ is bounded on $[\ell_+,X]$ by~\eqref{T'(x)=} and vanishes on $[0,\ell_+]$, we deduce that $T\in \Co^{0,1}([0,X])$. In particular, the relation~\eqref{PhiT} implies that $\Phi\in \Co^{1,1}([0,X])$. Thus, by~\eqref{Phi/tau+logrho} it holds 
\[
\log\rho\,\in \Co^{1,1}([0,X]).
\]
Hence, we conclude from~\eqref{T'(x)=} and assumption~(H2) of Theorem~\ref{th.main}, that $T\in \Co^{1,1}([\ell_+,X])$. Eventually, if $\ell_+>0$, since $T$ vanishes in $[0,\ell_+]$, the relation~\eqref{Phi/tau+logrho} implies
\[
\rho(x) = C\exp\lt(-\frac{x^2}{2\tau}\rt)\qquad\mbox{for }x\in [0,\ell_+],
\]
so that $\rho\in \Co^\oo([0,\ell_+])$. This ends the proof of the theorem.
\end{proof}\smallskip

\subsection{The optimization problem~\eqref{OP} as a time discretisation of~\eqref{P}}\label{subs:onestep_timediscretization}

\begin{theo}\label{theo.EL}
Let $(\rho,X)\in\A$  be a minimizer of the optimization problem~\eqref{OP}.Then, there holds
\be\label{estim_grad_L2}
\int_0^{X}\left|\rho'(x)\right|^2\, dx\le\|\rho\|_{L^\oo(\R_+)}\,\dfrac{\Wass^2(\rho,\rho^{0})}{\tau^2}.
\ee
Next, for every $\psi\in\Co^\oo_c(\R_+)$ 
\[
\int_{\R_+}\dfrac{\rho(x)-\rho^0(x)}{\tau}\,\psi(x)\, dx-\dfrac{\MM(\rho)-\MM(\rho^0)}{\tau}\psi(0)+\int_0^{X}\rho'(x)\,\psi'(x)\, dx = Q_\tau(\psi),
\]
where the right hand side satisfies
\be
\label{3.cont_EL_reste}
\lt|Q_\tau(\psi)\rt|\le\dfrac{\|\psi''\|_{L^{\oo}(\R_+)}}{2\tau}\,\Wass^2(\rho,\rho^0) +\dfrac{\|\psi'\|_{L^\oo(\R_+)}}{\tau}\,\lt(\int_0^{\ell_-}y\,\rho^0(y)\, dy\rt).
\ee
%We recall the convention $\ell_-=0$ when $\MM(\rho)\ge\MM(\rho^0)$.
\end{theo}

\begin{proof}We use the notation of the proof of Theorem~\ref{theo.ExistenceJKO.2}. In particular, $T$ is defined by~\eqref{notationT}. From~\eqref{Phi/tau+logrho} and~\eqref{PhiT}, we have $(\log\rho)'=(T(x)-x)/\tau$ on $[0,X]$, 
hence
\be\label{rho'}
\rho'(x)=\dfrac{T(x)-x}\tau\rho(x)\qquad\text{for }x\in[0,X].
\ee

\noindent
\textit{Step 1.}
Using~\eqref{rho'}, squaring and integrating on $[0,X]$, we have,
\[
\int_0^X(\rho')^2(x)\, dx=\dfrac1{\tau^2}\int_0^X(T(x)-x)^2\rho^2(x)\,dx\le\dfrac{\|\rho\|_{L^\oo(\R_+)}}{\tau^2}\int_0^X(T(x)-x)^2\rho(x)\,dx.
\]
This proves~\eqref{estim_grad_L2}.\medskip

\noindent
\textit{Step 2.} Let $\psi\in\Co^\oo_c(\R_+)$. Multiplying the last identity by $\psi'$ and integrating on $[0,X]$, we obtain, 
\[
\int_0^X\psi'(x)\,\dfrac{x-T(x)}{\tau}\,\rho(x)\, dx +\int_0^{X}\rho'(x)\,\psi'(x)\,dx  = 0.
\]
We then only have to show that
\be
\label{3.EL3}
\int_0^X\psi'(x)\,\dfrac{x-T(x)}{\tau}\,\rho(x)\, dx =\int_{\R_+}\dfrac{\rho(x)-\rho^0(x)}{\tau}\psi(x)\, dx -\dfrac{\MM(\rho)-\MM(\rho^0)}{\tau}\psi(0) + Q_{\tau}(\psi),
\ee
for some $Q_{\tau}(\psi)$ satisfying the estimate~\eqref{3.cont_EL_reste}. We consider two cases.
\medskip

\noindent
\textit{Case 1, $\MM(\rho)\ge\MM(\rho^0)$.} Writing the Taylor expansion,
\[
\psi(T_+(x))=\psi(x)+\psi'(x)(T_+(x)-x) +\int_x^{T_+(x)}(T_+(x)-y)\,\psi''(y)\,dy,
\]
we get 
\be\label{technical1}
\psi'(x) (x-T_+(x)) =\psi(x)-\psi(T_+(x))+R_1(x)\qquad\text{with}\quad|R_1(x)|\le\|\psi''\|_{L^{\oo}(\R_+)}\dfrac{(T_+(x)-x)^2}2.
\ee
Multiplying by $\rho(x)/\tau$ and integrating, we obtain,
\be\label{3.EL4}
\int_0^X\psi'(x)\,\dfrac{x-T_+(x)}{\tau}\,\rho(x)\, dx =\dfrac{1}{\tau}\int_{0}^{X}\rho(x)\,\psi(x)\, dx -\dfrac{1}{\tau}\int_{0}^X\rho(x)\,\psi(T_+(x))\, dx + Q_\tau(\psi),
\ee
where $Q_\tau(\psi)$ satisfies
\[
\lt|Q_\tau(\psi)\rt|\le\dfrac{\|\psi''\|_{L^{\oo}(\R_+)}}{2\tau}\Wass^2(\rho,\rho^0) .
\]
Recalling the definition~\eqref{T+} of the mapping $T_+$, we get
\begin{align*}
\int_0^X\rho(x)\,\psi(T_+(x))\, dx &=\int_0^{\ell_+}\psi(0)\,\rho(x)\,dx +\int_{\ell_+}^X\psi(T_+(x))\,\rho(x)\,dx\\
&=\lt[\MM(\rho)-\MM(\rho^0)\rt]\psi(0) +\int_0^{X^0}\rho^0(y)\,\psi(y)\, dy.
\end{align*}
Putting this in~\eqref{3.EL4}, we conclude that~\eqref{3.EL3} holds in this case.
\smallskip

\noindent\textit{Case 2, $\MM(\rho)<\MM(\rho^0)$.} Applying~\eqref{technical1} with $T_-$ instead of $T_+$, we obtain 
\be\label{3.EL45}
\int_0^X\psi'(x)\,\dfrac{x-T_-(x)}{\tau}\,\rho(x)\, dx  =\dfrac{1}{\tau}\int_0^X\rho(x)\,\psi(x)\, dx -\dfrac{1}{\tau}\int_0^X\rho(x)\,\psi(T_-(x))\,dx + Q^a_{\tau}(\psi),
\ee
with
\be\label{3.EL5}
\lt|Q_{\tau}^{a}(\psi)\rt|\le\dfrac{\|\psi''\|_{L^{\oo}(\R_+)}}{2\tau}\Wass^2(\rho,\rho^0).
\ee
Using the definitions~\eqref{ell-} and~\eqref{T-} of $\ell_-$ and $T_-$, we compute 
\begin{align*}
\int_0^X\rho(x)\,\psi(T_-(x))\,dx  &=\int_{\ell_-}^{X^0}\rho^0(y)\,\psi(y)\, dy\\
&=\int_0^{X^0}\rho^0(y)\,\psi(y)\, dy -\lt[\MM(\rho^0)-\MM(\rho)\rt]\psi(0)  -\int_0^{\ell_-}\rho^0(y)\lt(\psi(y)-\psi(0)\rt)\, dy.
\end{align*}
Putting this in~\eqref{3.EL45}, we obtain
\begin{multline}
\label{3.EL6}
\int_0^X\psi'(x)\,\dfrac{x-T_-(x)}{\tau}\,\rho(x)\, dx =\int_{\R_+}\psi(x)\,\dfrac{\rho(x)-\rho^0(x)}{\tau}\, dx
 -\dfrac{\MM(\rho)-\MM(\rho^0)}\tau\psi(0)\\
+ Q_{\tau}^{a}(\psi)+ Q_{\tau}^{b}(\psi),
\end{multline}
where the last term  is
\[
Q_{\tau}^{b}(\psi):=\dfrac1\tau\int_0^{\ell_-}\rho^0(y)\, (\psi(y)-\psi(0))\, dy\quad\text{ for which }\quad\lt|Q_{\tau}^{b}(\psi)\rt|\le\dfrac{\|\psi'\|_{L^{\oo}(\R_+)}}{\tau}\lt(\int_0^{\ell_-}y\rho^0(y)\, dy\rt). 
\]
Eventually,~\eqref{3.EL3} follows from  the last estimate,~\eqref{3.EL5} and~\eqref{3.EL6}. This completes the proof of the theorem.
\end{proof}

\subsection{Behavior of the minimizers at the fixed interface}
\label{Sec:comp_inter_fixe}~

\begin{figure}[ht]
\centering
\captionsetup{width=1\textwidth}
\includegraphics[]{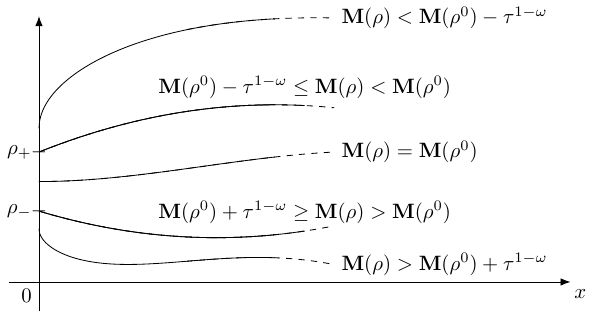}
\caption{Possible behaviors of $\rho$ at $x=0$ according to~\eqref{bcx=0.bis}. Compare with Figure~\ref{figcbx=0}.\label{figcbx=0.bis}}
\end{figure}

The next result describes the boundary condition satisfied by $\rho$ at $x=0$, see Figure~\ref{figcbx=0.bis}. 
\begin{prop}
\label{prop.BCx=0}
Let $(\rho,X)\in\A$ be a minimizer of~\eqref{OP}, then we have the following alternatives.
\be\label{bcx=0.bis}
\begin{cases}
\log\rho(0)=\log\rop + p'_\tau(\MM(\rho^0)-\MM(\rho))&\text{if }\MM(\rho^0)>\MM(\rho),\smallskip\\
\qquad\rom\le\rho(0)\le\rop&\text{if }\MM(\rho)=\MM(\rho^0),\smallskip\\
\log\rho(0)=\log\rom - p'_\tau(\MM(\rho)-\MM(\rho^0))&\text{if }\MM(\rho)>\MM(\rho^0).
\end{cases}
\ee
Recalling the definition~\eqref{2.def.penalty} of the penalty function $p_\tau$, we notice that for $M>0$, $p'_\tau(M)= (M\tau^{\om-1}- 1)_{+}\ge0$.
\end{prop}

\begin{proof}[Proof of Proposition~\ref{prop.BCx=0}]

We split the proof in four steps. In each step, we introduce some values $m_\eps$ and some candidates $(\rho_\eps,X)\in\Abb$ depending on $\eps>0$. Please note that these objects are not the same from one step to the next.\medskip

\noindent 
\textit{Step 1.} We first consider the case $\MM(\rho)\le\MM(\rho^0)$ and we set  $M:=\MM(\rho)-\MM(\rho^0)\le 0$. We use the notation $\ell_-$ from~\eqref{ell-}, {\it i.e.}
\[
\int_0^{\ell-}\rho^0\, = -M,
\]
and we define $\mu^0  =\rho^0\,\Leb $ and $\mu =\rho\,\Leb + (-M)\delta_0$. By Theorem~\ref{th.rappel}, there exists a nondecreasing Lipschitz mapping $S_0:[0,X^0]\to[0,X]$ such that $\mu=S_0\pf\mu^0$.  In particular,
$S_0=0$ on $[0,\ell_-]$ and $S_0$ is increasing on $[\ell_-,X^0]$ with $S_0([\ell_-,X^0])=[0,X]$.

Let $\eps\in(0,1)$ small enough so that $y_\eps:=\ell_-+\sqrt\eps<X^0$ and let us define 
\[
 x_\eps:=S_0(y_\eps)\qquad\text{and}\qquad m_\eps:=\int_{\ell_-}^{\ell_-+\eps}\rho^0\,>0.
\]
We notice that since $\rho^0,\rho\in\A$ and $S_0$ is Lipschitz continuous, there exists some $c>0$ such that
\[
c\sqrt\eps\le x_\eps\le c^{-1}\sqrt\eps\quad\text{ and}\quad c\eps\le m_\eps\le c^{-1}\eps.
\]
Let us set for $x\ge0$, 
\[
\rho_\eps(x):=
\begin{cases}
\rho(x) -\dfrac{m_\eps}{x_\eps}&\text{for }x\in [0,x_\eps],\\
\ \quad\rho(x)&\text{for }x\in (x_\eps, X].\\
\ \ \quad0&\text{for }x>X.
\end{cases}
\] 
The pair $(\rho_\eps,X)\in\Abb$ is a candidate for the optimization problem~\eqref{OP} which satisfies  $\MM(\rho_\eps)=\MM(\rho)-m_\eps$.  Since $\rho$ is continuous and positive on $[0,X]$ and $m_\eps/x_\eps=O(\sqrt\eps)$, we have $(\rho_\eps,X)\in\Abb$ for $\eps>0$ small enough.  For such $\eps$ we have by optimality of $(\rho,X)$,
\[
\J^\tau_{(\rho^0,X^0)}(\rho,X)\le\J^\tau_{(\rho^0,X^0)}(\rho_\eps,X).
\]
This rewrites as 
\begin{multline}\label{proof.propBCx=0.1}
q(\eps):=\int_0^{x_\eps}f(\rho(x))\, dx -\int_0^{x_\eps}f\lt(\rho(x) -\dfrac{m_\eps}{x_\eps}\rt)\, dx -\theta m_\eps +p_\tau(-M)-p_\tau(-M+m_\eps)\\
\le\dfrac1{2\tau}\lt(\Wass^2\lt(\rho_\eps,\rho^0\rt) -\Wass^2\lt(\rho,\rho^0\rt)\rt).
\end{multline}
Sending $\eps$ to 0 we have for the  left hand side of~\eqref{proof.propBCx=0.1},
\be\label{proof.propBCx=0.2}
q(\eps)\ \st{\eps\dw0}\sim\  
\lt[f'(\rho(0)) -\theta-p_\tau'(-M)\rt]m_\eps=\lt[\log\rho(0)-\log\rop -p_\tau'(\MM(\rho^0)-\MM(\rho))\rt]m_\eps.
\ee
Here, we have used the identities $f'(r)=\log r +\beta$ and $-\beta+\theta=\log\rop$ from the definitions of $\theta$, $\beta$ and $f$ in~\eqref{def_betatheta} and~\eqref{f}. Recalling $c\eps\le m_\eps\le c^{-1}\eps$, it is enough to establish that 
\be\label{proof.propBCx=0.3}
\Wass^2\lt(\rho_\eps,\rho^0\rt) -\Wass^2\lt(\rho,\rho^0\rt)= o(\eps),
\ee
to deduce from~\eqref{proof.propBCx=0.1} and~\eqref{proof.propBCx=0.2} that 
\be\label{proof.propBCx=0.4}
\log\rho(0)\le\log\rop +p_\tau'(\MM(\rho^0)-\MM(\rho)).
\ee 
Let us prove~\eqref{proof.propBCx=0.3}. The optimal transport map sending $\mu^0$ to $\mu_\eps:=\rho_\eps\Leb +(-M+m_\eps)\delta_0$ is the unique non-decreasing  mapping $S_\eps:[0,X^0]\to[0,X]$ such that $S_\eps\pf\mu^0=\mu_\eps$. We have by Theorem~\ref{th.rappel},
\be\label{proof.propBCx=0.5}
\begin{cases}
S_\eps(y)=0&\text{for }0\le y\le\ell_-+\eps,\smallskip\\
S_\eps(x)=S_0(x)&\text{for }y_\eps =\ell_-+\sqrt{\eps}\le y\le X^0,\smallskip\\
\ds\int_{\ell_-+\eps}^y\rho^0(y)\, dy =\ds\int_0^{S_\eps(y)}\lt(\rho(x)-\dfrac{m_\eps}{x_\eps}\rt)\, dx&\text{for }0\le\ell_-+\eps\le y\le y_\eps.
\end{cases}
\ee
Notice that for $\eps>0$ small enough, $S_\eps$ is Lipschitz continuous  with a uniform Lipschitz constant.

 Let us now estimate $|S_\eps-S_0|$. We first observe that 
\[
S_\eps=S_0\quad\text{on }[0,\ell_-]\cup[\ell_-+\sqrt\eps ,X^0]. 
\]
Besides, since $S_0(\ell_-)=0$ and $S_\eps=0$ on $[0,\ell_-+\eps]$, we have
\[
|S_\eps-S_0|\le\Lip(S_0)\eps\quad\text{ on }[\ell_-,\ell_-+\eps].
\]
For the remaining values, we deduce from the last case of~\eqref{proof.propBCx=0.5} that  for $y\in[\ell_-+\eps ,\ell_-+\sqrt\eps ]$,
\[
\lt[\rho(S_\eps(y)) -\dfrac{m_\eps}{x_\eps}\rt]S_\eps'(y)=\rho^0(y)=\rho(S_0(y))S_0'(y).
\]
Rearranging the terms we obtain,
\[
\rho(S_0(y))\lt[S_\eps'(y)-S_0'(y)\rt]= (\rho(S_0(y))-\rho(S_\eps(y)))S_\eps'(y)   +\dfrac{m_\eps}{x_\eps}S_\eps'(y).
\]
Denoting $\eta:=\min_{[0,X]}\rho>0$, this leads to 
\[
\lt|(S_\eps-S_0)'\rt(y)|\le\dfrac{\Lip(S_\eps)}\eta\lt[\dfrac{m_\eps}{x_\eps}+\Lip(\rho_{\,|[0,X]}) |S_\eps-S_0|(y)\rt].
\]
Integrating the above inequality and using the estimates $|S_\eps-S_0|(\ell_-+\eps)\le\Lip(S_0)\eps$ and $m_\eps/x_\eps\le c^{-2}\sqrt\eps$, we obtain that there exists $C\ge0$ such that for $\eps>0$ small enough, 
\[
|S_\eps-S_0|(y)\le C\lt(\eps +\int_{\ell_-+\eps}^y|S_\eps-S_0|(z)\, dz\rt)\qquad\text{for }y\in\lt[\ell_-+\eps ,\ell_-+\sqrt\eps\rt].
\]
By Gronwall lemma we deduce that 
\[
|S_\eps-S_0|\le C'\eps\qquad\text{on }\lt[\ell_-+\eps ,\ell_-+\sqrt\eps\rt],
\]
for some $C'\ge0$. Increasing  $C'$ if necessary, we have $C'\ge\Lip(S_0)$ and the inequality holds true also on $[\ell_-,\ell_-+\eps)$. Now, we can estimate  
\begin{align*}
\lt|\Wass^2\lt(\rho_\eps,\rho^0\rt) -\Wass^2\lt(\rho,\rho^0\rt)\rt|
 &=\lt|\int_{\ell_-}^{\ell_-+\sqrt\eps}\lt(\lt[y-S_\eps(y)\rt]^2-\lt[y-S_0(y)\rt]^2\rt)\rho^0(y)\, dy\rt|\\
 &\le\int_{\ell_-}^{\ell_-+\sqrt\eps}\lt(S_\eps(y)+2y+S_0(y)\rt)\,|S_\eps(y)-S_0(y)|\rho^0(y)\, dy\\
 &\le 2\max(X^0,X)\lt(\max_{[0,X^0]}\rho^0\rt)\int_{\ell_-}^{\ell_-+\sqrt\eps}C'\eps\, dy\\
 &= 2C'\max(X^0,X)\lt(\max_{[0,X^0]}\rho^0\rt)\eps^{3/2}=o(\eps).
\end{align*}
This establishes~\eqref{proof.propBCx=0.3} and thus~\eqref{proof.propBCx=0.4}. Let us notice here that in the limit case $M=0$, we have $p'_\tau(M)=0$ and~\eqref{proof.propBCx=0.4} leads to 
\be\label{proof.propBCx=0.55}
\MM(\rho)=\MM(\rho^0)\ \implies\ \log\rho(0)\le\log\rop,
\ee 
which is the second inequality of the case $\MM(\rho)=\MM(\rho^0)$ of the proposition.\medskip

\noindent
\textit{Step 2.} Let us now assume that $M=\MM(\rho)-\MM(\rho^0)<0$ and let us establish the reverse inequality of~\eqref{proof.propBCx=0.4}. Instead of sending more mass from $\mu^0$ to the point $x=0$, we send less. More precisely, for $0<\eps<1$ such that  $\eps<\ell_-$  we define  $y_\eps:=\ell_-+\sqrt\eps<X^0$ and  
\[ 
m_\eps:=-\int_{\ell_--\eps}^{\ell_-}\rho^0\,<0.
\]
Then, denoting, $x_\eps=S_0(y_\eps)$ as above we introduce a new candidate $(\rho_\eps,X)\in\Abb$ for the optimization problem~\eqref{OP}, 
\[
\rho_\eps(x):=
\begin{cases}
\rho(x) -  m_\eps/ x_\eps&\text{for }x\in [0,x_\eps],\\
~\quad\rho(x)&\text{for }x\in (x_\eps, X].
\end{cases}
\] 
We have $\MM(\rho_\eps)=\MM(\rho)- m_\eps$ and repeating the computations of Step~1, we obtain~\eqref{proof.propBCx=0.2} and~\eqref{proof.propBCx=0.3}. As in Step~1, the right hand side of~\eqref{proof.propBCx=0.1} is negligible compared to $|m_\eps| \sim \rho^0(\ell_-)\eps$. However, now $m_\eps<0$ and when we send $\eps$ to 0, we obtain the reverse inequality of~\eqref{proof.propBCx=0.4}. This proves the case $\MM(\rho)<\MM(\rho^0)$ of the proposition.
\medskip

\noindent
\textit{Step 3.} We now assume that $M:=\MM(\rho)-\MM(\rho^0)\ge 0$. For $\eps>0$ such that $\eps<X$, we define the candidate $(\rho_\eps,X)\in\Abb$ for~\eqref{OP} by 
\be\label{proof.propBCx=0.6}
\rho_\eps(x)=\begin{cases}
\rho(x)+\eps&\text{for }0\le x\le\eps,\smallskip\\
 \quad\rho(x)&\text{for }x>\eps.
\end{cases}
\ee
We have $\MM(\rho_\eps)=\MM(\rho)+\eps^2$. By optimality of $(\rho,X)$ we have $\J^\tau_{(\rho^0,X^0)}(\rho,X)\le\J^\tau_{(\rho^0,X^0)}(\rho_\eps,X)$ which reads, 
\begin{multline}\label{proof.propBCx=0.7}
q(\eps):=\int_0^{\eps}f(\rho(x))\, dx -\int_0^{\eps}f(\rho(x) +\eps)\, dx -\theta\eps^2 +p_\tau(M)-p_\tau(M+\eps^2)\\
\le\dfrac1{2\tau}\lt(\Wass^2\lt(\rho_\eps,\rho^0\rt) -\Wass^2\lt(\rho,\rho^0\rt)\rt).
\end{multline}
Sending $\eps$ to 0 we have for the  left hand side of~\eqref{proof.propBCx=0.7},
\be\label{proof.propBCx=0.8}
q(\eps)\ \st{\eps\dw0}\sim\ \lt[-f'(\rho(0))-\theta-p_\tau'(M)\rt]\eps^2=-\lt[\log\rho(0)-\log\rom +p_\tau'(\MM(\rho)-\MM(\rho^0))\rt]\eps^2,
\ee
where we have used $f'(r)=\log r +\beta$ and $-\beta-\theta=\log\rom$ from the definitions~\eqref{def_betatheta} and~\eqref{f}. Next, we estimate the right hand side of~\eqref{proof.propBCx=0.7}. By definition, $\mu^0=\rho^0\Leb + M\delta_0$ and 
\[
\Wass^2\lt(\rho,\rho^0\rt)=\wass^2\lt(\rho\Leb,\mu^0\rt).
\]
Similarly, denoting $\mu^0_\eps:=\rho^0\Leb+(M+\eps^2)\delta_0$, we have 
\[
\Wass^2\lt(\rho_\eps,\rho^0\rt)=\wass^2\lt(\rho_\eps\Leb,\mu^0_\eps\rt).
\]
Let $\gamma$ be the optimal transport plan from $\mu$ to $\mu^0$. Setting
\[
\gamma_\eps:=\gamma +\eps\lt(\Leb\restr(0,\eps)\rt)\otimes\delta_0,
\] 
we have $\gamma_\eps\in\PI(\rho_\eps\Leb,\mu^0_\eps)$ and 
\begin{align*}
\Wass^2\lt(\rho_\eps,\rho^0\rt)
&=\wass^2\lt(\rho_\eps\Leb,\mu^0_\eps\rt)\\
&\le\iint_{[0,X]\times[0,X^0]}|x-y|^2\, d\gamma_\eps(x,y)\\
&=\iint_{[0,X]\times[0,X^0]}|x-y|^2\, d\gamma(x,y)+\eps\int_0^\eps x^2\, dx\\
&=\wass^2\lt(\rho\Leb,\mu^0\rt) +\eps^4/3=\Wass^2\lt(\rho,\rho^0\rt) +\eps^4/3.
\end{align*}
Using this inequality in the right hand side of~\eqref{proof.propBCx=0.7}, dividing by $\eps^2$ and sending $\eps$ to 0 we deduce from~\eqref{proof.propBCx=0.8} the inequality
\be\label{proof.propBCx=0.9}
\log\rho(0)\ge\log\rom -p_\tau'(\MM(\rho)-\MM(\rho^0)).
\ee
In particular, 
\[
\MM(\rho)=\MM(\rho^0)\ \implies\ \log\rho(0)\ge\log\rom.
\]
Together with~\eqref{proof.propBCx=0.55}, this proves the proposition in the case $\MM(\rho)=\MM(\rho^0)$.
\medskip

\noindent
\textit{Step 4.} Eventually, we assume that $M=\MM(\rho)-\MM(\rho^0)>0$ and we establish the reverse inequality of~\eqref{proof.propBCx=0.9}. 
For this, contrarily to~\eqref{proof.propBCx=0.6} we remove a small mass from $\mu=\rho\Leb$ near 0. Let $\eps>0$ such that $\eps^2<M$, $\eps<X$ and $\eps<\min_{[0,X]}\rho$, we define a candidate $(\rho_\eps,X)\in\Abb$ for~\eqref{OP} by 
\[
\rho_\eps(x)=\begin{cases}
\rho(x)-\eps&\text{for }0\le x\le\eps,\smallskip\\
\quad \rho(x)&\text{for }x>\eps.
\end{cases}
\]
We have $\MM(\rho_\eps)=\MM(\rho)-\eps^2$. Following the reasoning of the previous step, the optimality of $(\rho,X)$ now leads to
\be\label{proof.propBCx=0.10}
\log\rho(0)-\log\rom +p_\tau'(\MM(\rho)-\MM(\rho^0))\le\dfrac1{2\tau\eps^2}\lt(\Wass^2\lt(\rho_\eps,\rho^0\rt) -\Wass^2\lt(\rho,\rho^0\rt)\rt)+o(1).
\ee
In order to establish that the right hand side is negligible we proceed as in the previous step. We define the transport plan from $\mu_\eps:=\rho_\eps\Leb$ to $\mu_0-\eps^2\delta_0$ as
\[
\gamma_\eps:=\gamma-\eps\lt(\Leb\restr(0,\eps)\rt)\otimes\delta_0,
\] 
which belongs to  $\PI(\rho_\eps\Leb,\mu^0_\eps)$ (recall that $\eps^2<M$ and $\eps<\min_{[0,X]}\rho$ so that $\gamma_\eps$ is a positive measure).
The rest of the proof is exactly as in Step~3. We obtain that 
\[
\Wass^2\lt(\rho_\eps,\rho^0\rt) -\Wass^2\lt(\rho,\rho^0\rt)\le\eps^4/3,
\]
which with~\eqref{proof.propBCx=0.10}, this leads to the reverse inequality of~\eqref{proof.propBCx=0.9} and settles the case $\MM(\rho)>\MM(\rho^0)$. 
\end{proof}
\medskip

\subsection{The optimisation problem with a prescribed mass $\MM(\rho)$}~
\label{subsec.M(rho)prescibed}

In Subsection~\ref{subsec.LooLipbounds}, in order to establish some bounds on $u$ and $\Lip(u)$ for $u=(\log\rho)_{|[0,X]}$, we will work with smooth approximations of $\rho^0$ on $[0,X^0]$. The goal of the current subsection is to justify this approximation procedure. 

We study here a simplification of the optimisation problem~\eqref{OP} where the total mass $\MM(\rho)$ of the candidates is given.  
For $M>0$, we define 
\[
\Abb_M:=\lt\{(\rho,X)\in\Abb:\MM(\rho)=M\rt\}\qquad\text{and}\qquad\A_M:=\{(\rho,X)\in\A:\MM(\rho)=M\}.
\]
Given $(\rho^0,X^0)\in\A$, we consider the optimisation problem: 
\be\label{OP2}
\min_{(\rho,X)\in\Abb_M}\,\J^\tau_{(\rho^0,X^0)}(\rho,X).
\ee
In this setting, the terms $\theta |\MM(\rho)-\MM(\rho^0)|$ and $p_\tau(\MM(\rho)-\MM(\rho^0))$ in  the definition of the functional  $\J^\tau_{(\rho^0,X^0)}$ (see~\eqref{Jtau2}) are irrelevant. Removing those terms, we define 
\be\label{Itau}
\I^\tau_{(\rho^0,X^0)}(\rho,X):=\dfrac{\Wass^2\left(\rho,\rho^0\right) +\lambda\left(X-X^0\right)^2 }{2\tau}+\PSI(\rho,X), 
\ee
so that, for $(\rho,X)\in\Abb_M$, there holds
\[
\J^\tau_{(\rho^0,X^0)}(\rho,X)=\I^\tau_{(\rho^0,X^0)}(\rho,X) +\theta |M-\MM(\rho^0)| + p_\tau(M-\MM(\rho^0)).
\]
As a consequence, the optimization problem~\eqref{OP2} is equivalent to
\be
\label{OP3}
\min_{(\rho,X)\in\Abb_M}\,\I^\tau_{(\rho^0,X^0)}(\rho,X).
\ee
\begin{prop}
\label{propOP3.1}
Let $(\rho^0,X^0)\in\A$ and $M>0$.  The optimization problem~\eqref{OP3} has a unique solution $(\rho,X)\in\Abb_M$.  We denote it
\[
\Ups^\tau((\rho^0,X^0),M):=(\rho,X).
\]
Besides, $(\rho,X)\in\A_M$ and  $\rho$ is of class $\Co^{1,1}$ on $[0,X]$.
\end{prop}

\begin{Remark}\label{rem.rho=Upsilon.1}
By uniqueness, if $(\rho,X)\in\A$ is a minimizer of the original problem~\eqref{OP} then 
\[(\rho,X)=\Ups^\tau((\rho^0,X^0),\MM(\rho)).\] 
However, we do not claim that~\eqref{OP} admits a unique minimizer. We only have that, if $(\rho_0,X_0)$, $(\rho_1,X_1)$ are minimizers of~\eqref{OP}, there holds $\MM(\rho_1)=\MM(\rho_0)\implies (\rho_1,X_1)=(\rho_0,X_0)$. 
\end{Remark}\medskip

\begin{proof}[Proof of Proposition~\ref{propOP3.1}]~

\noindent
\textit{Step 1. Existence and regularity.} The method is similar to the one used in Subsection~\ref{subs.existenceOP} for the existence of a solution to the original optimization problem~\eqref{OP}.  First, we introduce a reference candidate $(\hat\rho,X^0)\in\A_M$ defined by  $\hat\rho=(M/X^0){\bf{1}}_{[0,X^0]}$. Next, 
as in the proof of Theorem~\ref{theo.ExistenceJKO.1}, we consider a minimizing sequence $(\rho_j,X_j)\in\Abb_M$ for~\eqref{OP3}. By optimality, there holds,
\[
\limsup_{j\up\oo}\I^\tau_{(\rho^0,X^0)}(\rho_j,X_j)\le\I^\tau_{(\rho^0,X^0)}(\hat\rho,X^0) =: Q(\rho^0,X^0,M).
\]
 Hence, for $j\ge 1$,
\[
\PSI(\rho_j,X_j) +\dfrac\lambda{2\tau}(X^0-X_j)^2\le  Q(\rho^0,X^0,M) +o(1).
\]
This leads to
\be\label{proof.propO3.1}
\dfrac\lambda{2\tau}(X_j-X^0)^2 +\ros\,X_j\le Q(\rho^0,X^0,M)-\int_0^{X_j}f(\rho_j) +o(1).
\ee
Now, since $r\mapsto f(r)=r(\log r+\beta-1)$ is convex and $\MM(\rho_j)=M$ we have by Jensen's inequality,
\[
\int_0^{X_j}f(\rho_j)\ge X_jf\lt(\dfrac M{X_j}\rt)=\ M\dfrac{f\lt(M/X_j\rt)}{M/X_j}\ \st{X_j\dw0}\longto\ +\oo,
\]
where we used $f(r)/r\to+\oo$ as $r\up\oo$. Therefore, sending $X_j$ to 0 in~\eqref{proof.propO3.1} yields the contradiction 
\[
\dfrac\lambda{2\tau}(X^0)^2\le -\oo.
\] 
We deduce that there exists $\eps>0$ such that $X_j\ge\eps$ for every $j\ge1$. Now, we can follow the proof of Theorem~\ref{theo.ExistenceJKO.1} after~\eqref{proofprop.ExistenceJKO_1} and conclude that the optimization problem~\eqref{OP3} admits at least one solution $(\rho,X)\in\Abb_M$.  Then, proceeding exactly as in the proof of Theorem~\ref{theo.ExistenceJKO.2}, we obtain that $(\rho,X)\in\A_M$ and  that $\log\rho$ is of class $\Co^{1,1}$ on $[0,X]$.\medskip

\noindent
\textit{Step 2. Preliminaries about convexity along generalized geodesics.} We still have to establish that the minimizer of the optimization problem~\eqref{OP3} is unique. This is done in Step~3 using  the notion of generalized geodesics (see~\cite[Chapter~7.3.3]{OTAM}). For the convenience of the reader, we give here the definitions and prove all the necessary results in our particular setting.

Let  $\nu$, $\nu_0$, $\nu_1$ be  three positives measures on $\R$. We make the following assumptions: 
\begin{enumerate}[(h1)]
\item $0<\nu(\R)=\nu_0(\R)=\nu_1(\R)<\oo$.
\item $\nu$, $\nu_0$, $\nu_1$ are absolutely continuous with respect to the Lebesgue measure. We write,
\[
\nu=\sigma\Leb,\qquad\nu_0=\sigma_0\Leb,\qquad\nu_1=\sigma_1\Leb.
\]
\item The supports of $\sigma$, $\sigma_0$ and $\sigma_1$ are bounded intervals of the following form:
\[
\supp\sigma=[a,b]\qquad\text{and}\qquad\supp\sigma_i=[0,d_i]\quad\text{for } i=0,1.
\]
\item There exists $\kappa\in(0,1)$ such that for $\rho\in\{\sigma,\sigma_0,\sigma_1\}$,
\[
\kappa\le\rho\le 1/\kappa\quad\text{ almost everywhere on }\supp\rho.
\]
\end{enumerate}
For $i=0,1$, let $S_i:[a,b]\to\R$ be the unique optimal transport map from $\nu$ to $\nu_i$. It is nondecreasing and characterized by 
\be\label{proof.propOP3.1.1}
\int_a^y\sigma(z)\, dz =\int_0^{S_i(y)}\sigma_i(x)\, dx\qquad\text{for }y\in [a,b].
\ee
In particular, $S_i\pf\nu=\nu_i$. We also observe that differentiating~\eqref{proof.propOP3.1.1}, we get 
\be\label{proof.propOP3.1.1.5}
\sigma(y)=S_i'(y)\sigma_i(S_i(y))\quad\text{for  almost every }y\in[a,b].
\ee
We deduce from~(h4) that $S_i$ is increasing and Lipschitz continuous with 
\be\label{proof.propOP3.1.2}
S_i(0)=a\qquad\text{and}\qquad\kappa^2\le S_i'\le\kappa^{-2}\qquad\text{almost everywhere on }[a,b].
\ee
 Next, we define for $t\in[0,1]$, 
\be\label{proof.propOP3.1.3}
S_t:=(1-t)S_0+tS_1\qquad\text{and}\qquad\nu_t:=S_t\pf\nu.
\ee
By construction, for $t=0$ and $t=1$, $\nu_t$ coincides with $\nu_0$ and $\nu_1$. Moreover, for $t\in[0,1]$, $\nu_t$ is a positive measure. The path $t\in[0,1]\mapsto\nu_t$ is called a generalized geodesic from $\nu_0$ to $\nu_1$. We also see that $\nu_t(\R)=\nu(\R)$ and from~\eqref{proof.propOP3.1.2} we get that $S_t$ is Lipschitz continuous with 
\[
\kappa^2\le S_t'\le\kappa^{-2}\qquad\text{almost everywhere on }[a,b].
\]
We deduce that $\nu_t$ is absolutely continuous with respect to the Lebesgue measure and writing $\nu_t=\sigma_t\Leb$, we have 
\be\label{proof.propOP3.1.4}
\supp\sigma_t=[0,d_t]\quad\text{with } d_t:=(1-t)d_0+td_1,
\ee
and 
\be\label{proof.propOP3.1.5}
\kappa^3\le\sigma_t\le\kappa^{-3}\qquad\text{almost everywhere on }[0,d_t].
\ee

Let us now establish the convexity of $t\in [0,1]\mapsto\mathcal{F}(\nu_t)$ for three different types of functionals.\medskip

\noindent
\textit{Step 2.1.}
Let $\hat x,\hat q_0,\hat q_1\in\R$ with $\hat q_1\ge 0$. Recalling~\eqref{proof.propOP3.1.4}, we see that 
\be\label{proof.propOP3.1.5.5}
t\in[0,1]\mapsto\hat q_1(d_t-\hat x)^2 +\hat q_0  d_t\qquad\text{ is convex}.
\ee

\noindent
\textit{Step 2.2.}
Let $Z:[a,b]\to\R$ be a measurable bounded function. We see from~\eqref{proof.propOP3.1.3} that 
\be\label{proof.propOP3.1.6}
E_Z: t\in[0,1]\mapsto\int_a^b |S_t(y) - Z(y)|^2\sigma(y)\, dy=\int |S_t - Z|^2\,d\nu\qquad\text{ is convex}.
\ee
We use this in Step 3, with $Z(y)=y$ and $Z(y)=y_+$. In the former case, $E_Z(t)=\wass^2(\nu,\nu_t)$.
\smallskip

\noindent
\textit{Step 2.3. Generalized strict convexity of the entropy.}
Let us define for $t\in[0,1]$
\[
\mathcal{F}(\sigma_t):=\int_0^{d_t}f(\sigma_t(x))\, dx,
\]
where we recall that $f(r)=r(\log r +\beta -1)$. We claim that  
\be\label{proof.propOP3.1.7}
t\mapsto\mathcal{F}(\sigma_t)\text{ is convex on }[0,1]\text{ and strictly convex if }\nu_1\ne\nu_0.
\ee
Here $\nu_1$ is a translate of $\nu_0$ means that $\sigma_1(x)=\sigma_0(x-\hat x)$ for some $\hat x\in\R$. 

Let us establish the claim. For $r>0$ we denote $g(r):=f(r)/r=\log r +\beta -1$. By~\eqref{proof.propOP3.1.5}, we have  $\sigma_t\ge\kappa^3>0$ on $[0,d_t]$ and we can write 
\[
\mathcal{F}(\sigma_t)=\int_0^{d_t}g(\sigma_t(x))\sigma_t(x)\, dx =\int_a^b g(\sigma_t(S_t(y)))\sigma(y)\, dy.
\]
As in~\eqref{proof.propOP3.1.1.5}, we have $\sigma_t(S_t(y))=\sigma(y)/S_t'(y)$ for almost every $y\in[a,b]$, which leads to  
\[
\mathcal{F}(\sigma_t)=\int_a^b g\lt(\dfrac{\sigma(y)}{(1-t)S_0'(y)+tS_1'(y)}\rt)\sigma(y)\, dy.
\]
Now for every $c>0$, the mapping $r\in(0,+\oo)\mapsto g(c/r)=-\log(r)+\log(c) +\beta-1$ is strictly convex and we see that $t\mapsto\mathcal{F}(\sigma_t)$ is convex. If moreover, $\nu_1\ne\nu_0$ then $S_0\ne S_1$ which implies that $S_0'\ne S_1'$ on a subset of $[a,b]$ with positive measure (recall that $S_0(0)=S_1(0)=a$). In such case, the strict convexity of $r\mapsto g(c/r)$ implies that  $t\mapsto\mathcal{F}(\sigma_t)$ is also strictly convex. This proves the claim~\eqref{proof.propOP3.1.7}.
\medskip

\noindent
\textit{Step 3. Uniqueness.} Let $(\rho_0,X_0)$, $(\rho_1,X_1)\in\A_M$ be two minimizers of the optimization problem~\eqref{OP3} and let us establish that $(\rho_0,X_0)=(\rho_1,X_1)$. In the sequel, note carefully the distinction between $\rho^0$ and $\rho_0$. We treat successively the cases $M\le\MM(\rho^0)$ and $M>\MM(\rho^0)$.\medskip

\noindent
\textit{Step 3.1, $M\le\MM(\rho^0)$.} Let $\ell_-\in[0,X^0)$ such that 
\[
\int_0^{\ell_-}\rho^0=\MM(\rho^0)-M. 
\] 
Denoting $\nu=\rho^0\Leb\restr[\ell_-,X^0]$ and $\nu_i=\rho_i\Leb$ for $i=0,1$, we have, with the notation of Step~2, 
\[
\sigma=\rho^0\Leb\restr[\ell_-,X^0],\qquad\sigma_i=\rho_i\quad\text{for } i=0,1,
\]
and 
\[
[a,b]=[\ell_-,X^0],\qquad [0,d_i]=[0,X_i]\quad\text{for } i=0,1.
\]
 Besides the measures $\nu$, $\nu_1$, $\nu_2$ satisfy the assumptions~(h1)--(h4). Now, we define $\nu_t$,  $d_t$ and $\sigma_t$ for $t\in [0,1]$ as in Step~2. We see that $(\rho_t,X_t):=(\sigma_t,d_t)\in\Abb_M$ for $t\in[0,1]$.  Next, recalling the definition~\eqref{Itau} of  $\I^\tau_{(\rho^0,X^0)}$,
we have for $t\in[0,1]$,
\begin{multline*}
\I^\tau_{(\rho^0,X^0)}(\rho_t,X_t)=\dfrac1{2\tau}\wass^2\lt((\MM(\rho^0)-M)\delta_0,\rho^0\Leb\restr[0,\ell_-]\rt)\\
+\dfrac1{2\tau}\wass^2\lt(\nu,\nu_t\rt)+\lt[\dfrac{\lambda}{2\tau}\lt(d_t-b\rt)^2 -\ros d_t\rt] +\int_0^{d_t}f(\sigma_t(x))\, dx.
\end{multline*}
The first term is constant. For the second term, denoting $Z(y)=y$, we have 
\[
\wass^2\lt(\nu,\nu_t\rt)=E_Z(t),
\]
and by~\eqref{proof.propOP3.1.6}, this  is a convex function of $t$. The third term is also convex by~\eqref{proof.propOP3.1.5.5} and by~\eqref{proof.propOP3.1.7}, so is the last term. Moreover, this last term is strictly convex if $\nu_1\ne\nu_0$. Since $\rho_0$ and $\rho_1$ are optimal this is impossible. Hence, $\nu_1$ is a translate of $\nu_0$ and we have in fact $\nu_1=\nu_0$, that is $(\rho_0,X_0)=(\rho_1,X_1)$. This establishes the uniqueness of the minimizer of~\eqref{OP3} in this case.
\smallskip

\noindent
\textit{Step 3.2, $M>\MM(\rho^0)$.} Let us set $m:=M-\MM(\rho^0)$. We start with a preliminary remark. We have for $(\rho,X)\in\Abb_M$,
\[
\Wass^2(\rho^0,\rho)=\wass^2(m\delta_0+\rho^0\Leb,\rho\Leb).
\]
In order to apply the results of Step~2 we replace the measure $\mu^0=m\delta_0+\rho^0\Leb$  by a measure which is absolutely continuous with respect to the Lebesgue measure, namely,
\[
\nu:={\bf 1}_{[-m,0]}+\rho^0\Leb.
\]
Denoting $S:[-m,X^0]\to\R$ the nondecreasing transport map from $\nu$ to $\rho$. Writing $\nu=\sigma\Leb$, we observe that 
\be\label{proof.propOP3.1.8}
\Wass^2(\rho,\rho^0) =\int_{-m}^0 (S(y))^2\sigma(y)\, dy +\int_0^{X^0}(S(y)-y)^2\sigma(y)\, dy =\int_{-m}^{X^0}(S(y)-y_+)^2\sigma(y)\, dy,
\ee
with $y_+ = \max(y,0)$. Now, we set  $\nu_i=\rho_i\Leb$ for $i=0,1$ and we use the notation of Step~2: we also introduce the densities $\sigma_0$, $\sigma_1$, the supports $[a,b]=[-m,X^0]$, and $[0,d_i]=[0,X_i]$ for $i=0,1$ and the maps $S_0$, $S_1$. Next, we define $S_t$, $\sigma_t$ and  $[0,d_t]=[0,(1-t)X_0+tX_1]$. We see that $(\rho_t,X_t):=(\sigma_t,d_t)\in\Abb_m$ for $t\in[0,1]$. Now, using~\eqref{proof.propOP3.1.8}, we compute
\begin{align*}
\I^\tau_{(\rho^0,X^0)}(\rho_t,X_t)
 &=\dfrac1{2\tau}\Wass^2\left(\rho_t,\rho^0\right) +\lt[\dfrac{\lambda}{2\tau}\lt(d_t-b\rt)^2 -\ros d_t\rt] +\int_0^{d_t}f(\sigma_t(x))\, dx\\
 &=\dfrac1{2\tau}E_Z(t) +\lt[\dfrac{\lambda}{2\tau}\lt(d_t-b\rt)^2 -\ros d_t\rt] +\int_0^{d_t}f(\sigma_t(x))\, dx,
\end{align*}
where $Z(y):=y_+$ for $y\in\R$.  Exactly as in the previous case, $\I^\tau_{(\rho^0,X^0)}(\rho_t,X_t)$ is a convex function of $t$ and it is strictly convex if $\nu_1\ne\nu_0$. As this would contradict the optimality of $(\rho_0,X_0)$, $(\rho_1,X_1)$, we deduce again that $\nu_1=\nu_0$, that is $(\rho_0,X_0)=(\rho_1,X_1)$. As a conclusion, the minimizer is unique in all the cases. This ends the proof of  Proposition~\ref{propOP3.1}.
\end{proof}

The next result gives the boundary condition satisfied by $\Ups^\tau((\rho^0,X^0),M)$ at $x=X$.
\begin{prop}
\label{propOP3.2}
Let $(\rho^0,X^0)\in\A$, let $M>0$ and $(\rho,X):=\Ups^\tau((\rho^0,X^0),M)$. There holds,
\be\label{CBenX}
\lambda\dfrac{X-X^0}\tau=\rho(X)-\ros.
\ee
\end{prop}
\begin{Remark}%\label{rem.rho=Upsilon.2}
By Remark~\ref{rem.rho=Upsilon.1} any minimizer $(\rho,X)\in\A$ for the problem~\eqref{OP} satisfies~\eqref{CBenX}.
\end{Remark}

\begin{proof}
We consider here variations of $X$ (and of $\rho$ in the neighborhood of $X$). Let $\eps\in(0,1)$ with $\eps<X^2$. We consider as a competitor for the optimisation problem~\eqref{OP3} the pair, $(\rho_{\pm\eps},X_{\pm\eps})$ defined by $X_{\pm\eps}:=X\pm\eps$ and 
\[
\rho_{\pm\eps}(x):=
\begin{cases}
~\quad\rho(x)&\text{for }0\le x\le X-\sqrt\eps,\\
c_{\pm\eps}\rho(A_{\pm\eps}(x))&\text{for }X-\sqrt\eps  < x\le X\pm\eps,\\
~\quad\ 0&\text{for }x>X\pm\eps.
\end{cases}
\]
Here $A_{\pm\eps}$ is the increasing affine transformation sending $(X-\sqrt\eps,X\pm\eps]$ onto $(X-\sqrt\eps,X]$ and $c_{\pm\eps}$ is the factor ensuring $\int\rho_{\pm\eps}=\int\rho$. We compute,
\[
c_{\pm\eps}=\lt(1\pm\sqrt\eps\rt)^{-1}\qquad\text{and}\qquad A_{\pm\eps}(x))=X-\sqrt\eps +c_{\pm\eps}\lt(x-(X-\sqrt\eps)\rt).
\]
By optimality of $\rho$, we have 
\begin{multline*}
\dfrac1{2\tau}\Wass^2\left(\rho,\rho^0\right) +\lambda\dfrac{\left(X-X^0\right)^2}{2\tau}+\int_{X-\sqrt{\eps}}^Xf(\rho)\, +\ros X\\
\le\dfrac1{2\tau}\Wass^2\left(\rho_{\pm\eps},\rho^0\right) +\lambda\dfrac{\left(X_{\pm\eps}-X^0\right)^2}{2\tau}+\int_{X-\sqrt{\eps}}^{X\pm\eps}f(\rho_{\pm\eps})\, +\ros (X\pm\eps).
\end{multline*}
Substituting the values of $(\rho_{\pm\eps},X_{\pm\eps})$ and simplifying we get 
\begin{multline}\label{proofExistenceJKOb1}
-\lambda\lt(\pm\eps\dfrac{X-X^0}\tau\rt)\mp\eps\ros +\int_{X-\sqrt\eps}^X\lt(f(\rho(x)) -\dfrac1{c_{\pm\eps}}f(c_{\pm\eps}\rho(x))\rt)\, dx\\
\le\dfrac{\lambda\eps^2}{2\tau}+\dfrac1{2\tau}\lt(\Wass^2\left(\rho_{\pm\eps},\rho^0\right)-\Wass^2\left(\rho,\rho^0\right)\rt).
\end{multline}
Recalling that $\rho$ is continuous and positive on $[0,X]$, using the definitions of $c_{\pm\eps}$ and the identity $rf'(r)-f(r)=r$, we compute
\begin{align}
\nonumber
\int_{X-\sqrt\eps}^X\lt(f(\rho(x)) -\dfrac1{c_{\pm\eps}}f(c_{\pm\eps}\rho(x))\rt)\, dx
&=\pm\sqrt\eps\int_{X-\sqrt\eps}^X -f(\rho(x))+\rho(x) f'(\rho(x)) +O\lt(\sqrt\eps\rt)\, dx\\
\nonumber
&=\pm\sqrt\eps\int_{X-\sqrt\eps}^X\rho(x)\, dx+O\lt(\eps^{3/2}\rt)\\
\label{proofExistenceJKOb2}
&=\pm\eps\rho(X) +o\lt(\eps\rt).
\end{align}
Next, we pass from the measure $\rho\Leb$ to $\rho_{\pm\eps}\Leb$ by moving a mass of the order of $\sqrt{\eps}$ to a distance of the order of $\eps$. It follows that 
\[
\Wass^2\left(\rho_{\pm\eps},\rho^0\right)-\Wass^2\left(\rho,\rho^0\right) = O\lt(\eps^{3/2}\rt).
\]
Putting this identity and~\eqref{proofExistenceJKOb2} in~\eqref{proofExistenceJKOb1} and then dividing by $\pm\eps$ we get 
\[
\pm\lt(-\lambda\dfrac{X-X^0}\tau +\rho(X)-\ros\rt)\le o(1).
\]
Eventually, sending $\eps$ to 0, we get~\eqref{CBenX}.
\end{proof}

\medskip
The two last results of this subsection justify the smoothing procedure used in the next one.
\begin{prop}\label{propOP3.3a}
Let $(\rho^0,X^0)\in\A$ with $\rho^0\in \Co^{0,1}([0,X^0])$, let $M>0$ and 
\[
(\rho,X):=\Ups^\tau((\rho^0,X^0),M). 
\]
Then $\rho$ is of class $\Co^{1,1}$ on $[0,X]$ and of class $\Co^\oo$ on $[0,\ell_+]$.\\ 
Moreover, if $\rho^0$ is of class $\Co^{k,1}$ on $[0,X^0]$, for some $k\ge 0$, then:
\begin{enumerate}
\item[$\circ$] $\rho$ is of class $\Co^{k+2,1}$ on $[\ell_+,X]$,
\item[$\circ$] The transport map $T$ is of class  $\Co^{k+1,1}$ on $[\ell_+,X]$. 
\end{enumerate}
\end{prop}\smallskip

\begin{proof}We already know from Proposition~\ref{propOP3.1} that $\rho\in\A$, so that $\rho$ is of class $\Co^{1,1}$ on $[0,X]$. Moreover, arguing as in the proof of Theorem~\ref{theo.ExistenceJKO.2} we deduce that $\rho\in \Co^\oo([0,\ell_+])$ (if $\ell_+>0$) and $T$ is of class $\Co^{1,1}$ on $[\ell_+,X]$. Let us now assume that $\rho^0$ is of class $\Co^{k,1}$ on $[0,X^0]$ for some $k\ge 0$. We consider successively the two cases $M\le\MM(\rho^0)$ and $M<\MM(\rho_0)$.\medskip

\noindent
\textit{Step 1, $M\le\MM(\rho^0)$.} In this case, we have $\ell_+=0$ by convention and  $\ell_-\in[0,X^0)$ is defined by the relation 
\[
\int_0^{\ell_-}\rho^0(y)\, dy =\MM(\rho^0) - M.
\] 
As in the proof of Theorem~\ref{theo.ExistenceJKO.2}, denoting $T=T_-\in \Co^{0,1}([0,X], [\ell_-,X^0])$ the optimal mapping from $\rho\Leb$ to $\rho^0\Leb\restr[\ell_-,X^0]$, we have 
\be
\label{rho=rho0T'}
\rho(x)=\rho^0(T(x))T'(x)\qquad\text{for almost every }x\in[0,X].
\ee
Moreover, there exists a Kantorovitch potential $\Phi\in \Co^{0,1}([0,X],\R)$ such that 
\be\label{logrhoPhi}
(\log\rho)' (x)= -\dfrac{\Phi'(x)}\tau=\dfrac{T(x)-x}\tau\quad\text{for almost every }x\in[0,X].
\ee
Let us first consider the case $k=0$. We already know that $T$ is of class $\Co^{1,1}$ on $[0,X]$ and with ~\eqref{logrhoPhi} we see that $\log\rho$ and $\Phi$ are of class $\Co^{2,1}$ in $[0,X]$. In the case $k=0$, we can not improve the regularity of $T$ by~\eqref{rho=rho0T'} since we only assume that $\rho^0$ is of class $\Co^{0,1}$.\\
Now if $k\ge 1$,~\eqref{rho=rho0T'} implies that $T$ is of class $\Co^{2,1}$ on $[0,X]$ and then by~\eqref{logrhoPhi} that $\log\rho$ and $\Phi$ are of class $\Co^{3,1}$ on $[0,X]$. Repeating the argument, we get by induction that $T$ is of class $\Co^{k+1,1}$ on $[0,X]$ and that $\Phi$ and $\log\rho$ are of class $\Co^{k+2,1}$ on $[0,X]$.\smallskip

\noindent
\textit{Step 2, $M>\MM(\rho^0)$.} In this case, we have $\ell_-=0$ and $\ell_+>0$ is defined by 
\[
\int_0^{\ell_+}\rho(x)\, dx=M-\MM(\rho^0)=:m.
\]
As in the proof of Theorem~\ref{theo.ExistenceJKO.2}, denoting $T=T_+\in \Co^{0,1}([0,X], [0,X^0])$ the optimal mapping from $\mu$ to $\mu_0$, we have $T\pf\mu=\mu_0$ with $T(x)=0$ on $[0,\ell_+]$. Moreover~\eqref{logrhoPhi} holds true. This first leads to  
\be\label{exp-x2surtau}
\rho(x)=\rho(0)\exp\lt(-\dfrac{x^2}{2\tau}\rt)\qquad\text{ on }[0,\ell_+].
\ee
In particular, $\rho$ is smooth on $[0,\ell_+]$. Eventually, the identity $\rho(x)=\rho^0(T(x))T'(x)$ also holds in $[\ell_+,X]$ and proceeding as in Step~1, we get that $T$ is of class $\Co^{k+1,1}$ on $[\ell_+,X]$ and that $\Phi$ and $\log\rho$ are of class $\Co^{k+2,1}$ on $[\ell_+,X]$.
\end{proof}

\begin{prop}\label{propOP3.3b}
Let $(\rho^0,X^0)\in\A$ and let $(\rho^0_j,X^0)$ be a sequence of elements of $\A$ with $\MM(\rho^0_j)=\MM(\rho^0)$ for $j\ge 1$ and such that 
\[
\rho^0_j\to\rho^0\text{ in }\Co([0,X^0])\quad\text{ as }j\up\oo.
\] 
Let $M>0$, denoting  $(\rho_j,X_j):=\Ups^\tau((\rho^0_j,X^0),M)$ and $(\rho,X):=\Ups^\tau((\rho^0,X^0),M)$, we have 
\be\label{propOP3.3_1}
X_j\ \st{j\up\oo}\longto\ X,\qquad\rho_j\ \st{j\up\oo}\longto\ \rho\text{ in }L^1(\R_+)\qquad\text{and the sequence }(\Lip(\rho_j))\text{ is bounded}.
\ee  
\end{prop}\smallskip

\begin{proof}

Let $(\rho^0_j,X^0)\in\A$ and $(\rho^0_j,X_j)$ as in the statement of the proposition.\smallskip

\noindent 
\textit{Step 1.} We first establish that the sequence $(X_j)$ is uniformly bounded. Let us set $\ov M:=\max(M,\MM(\rho^0))$ and $\ov X:=\max(X^0,M)$ and let us define 
\[
\mathcal{S}:=\lt\{\nu\in\mathcal{M}_+([0,\ov X]),\nu([0,\ov X])=\ov M\rt\}.
\]
 The Wasserstein distance $\wass$ is continuous on $\mathcal{S}\times\mathcal{S}$ with respect to the weak-$\star$ convergence of measures. We have $\hat
\rho\in\mathcal{S}$ and by assumption, the measure  
\[
\mu^0_j:=(M-\MM(\rho^0))_+\delta_0 +\rho^0_j\Leb\in\mathcal{S},
\]
weakly (and in fact strongly) converges to $\mu_0:=(M-\MM(\rho^0))_+\delta_0 +\rho^0\Leb$. Defining $(\hat\rho,\hat X):=({\bf 1}_{[0,M]},M)$ as in the proof of Theorem~\ref{theo.ExistenceJKO.1}, it is then easy to see that
\[
C_j:=\I^\tau_{(\rho^0_j,X^0)}(\hat\rho,\hat X)\ \st{i\up\oo}\longto\ \I^\tau_{(\rho^0,X^0)}(\hat\rho,\hat X).
\]
In particular the sequence $(C_j)$ is bounded. Now, by optimality of $(\rho^j,X_j)$ and with the same computations as the ones following~\eqref{ineqXjYj}, we have
\begin{align*}
\oo> C_j&\ge\I^\tau_{(\rho^0_j,X^0)}(\rho_j,X_j)\\ 
 &\ge X_j f(M/X_j) +\ros X_j  +\lambda\,\dfrac{(X^0-X_j)^2}{2\tau}\\
 & = M(\log M -\log X_j +\beta-1)+\ros X_j  +\lambda\,\dfrac{(X^0-X_j)^2}{2\tau}\\
 &\ge  M(\log  M +\beta)   +(\ros-M)  X_j +\lambda\,\dfrac{(X^0-X_j)^2}{2\tau},
\end{align*}
where for the last inequality we have used the positivity of $M$ and the inequality $x-\log(x)-1\ge0$ for every $x\ge0$. Recalling that $M$ is prescribed, we deduce that $(X_j)$ is bounded.\smallskip

\noindent 
\textit{Step 2.} Let us denote $u_j:=\log((\rho_j)_{[0,X_j]})$ and let us show that $(\Lip(u_j))$ is bounded. With obvious notation, there exists Kantorovitch potentials $\Phi_j:[0,X_j]\to\R$ such that 
\[
u_j '=-\dfrac{\Phi_j'}\tau\qquad\text{on }[0,X_j].
\]
By~\eqref{Lip(Phi)} in Theorem~\ref{th.rappel}, we also have 
\[
\Lip(\Phi_j)\le 2\max(X_j,X^0),
\]
and since $(X_j)$ is bounded we conclude that $(\Lip(u_j))$ is also bounded.\smallskip

\noindent 
\textit{Step 3.} We now prove that the sequence $(X_j)$ is bounded from below by a positive constant.   For this, we repeat the method used at the beginning of the proof of Proposition~\ref{propOP3.1}. By~\eqref{proof.propO3.1} applied to $(\rho^0_j,X^0)$ we have
\[
\dfrac\lambda{2\tau}(X_j-X^0)^2 +\ros X_j\le\I^\tau_{(\rho^0_j,X^0)}(\hat\rho,\hat X)-\int_0^{X_j}f(\rho_j).
\]
Continuing as in the proof of Proposition~\ref{propOP3.1}, we get 
\[
\dfrac\lambda{2\tau}(X_j-X^0)^2 +\ros X_j\le\I^\tau_{(\rho^0_j,X^0)}(\hat\rho,\hat X) -M\dfrac{f\lt(M/X_j\rt)}{M/X_j}\ \st{X_j\dw0}\longto\ -\oo.
\]
 Hence there exists $\eps>0$ such that $X_j\ge\eps$ for every $j\ge1$.\medskip
\smallskip

\noindent 
\textit{Step 4.} We deduce from  Steps~1--3 and the fact that $\int_0^{X_j}\rho_j=\int_0^{X_j}\exp(u_j)=M$ that  $(\Lip(u_j))$ is bounded and  then that  $(u_j)$ is bounded. Hence, $(\Lip({\rho_j}_{|[0,X_j]}))$ is bounded and there exist constants $0<c\le C$ such that $0<c\le\rho_j\le C$ on $[0,X_j]$ for every $j$.   Up to extraction, there exists $(\rho,X)\in\A_M$ such that~\eqref{propOP3.3_1} holds true.  

It remains to show that 
$(\rho,X)=\Ups^\tau((\rho^0,X^0),M)=:(\ov\rho,\ov X)$. First, from~\eqref{propOP3.3_1} and the weak continuity of $\wass$ in $\mathcal{S}\times\mathcal{S}$, we have on the one hand
\[
\I^\tau_{(\rho^0,X^0)}(\rho,X)=\lim_{j\up\oo}\I^\tau_{(\rho^0_j,X^0)}(\rho_j,X_j).
\]
On the other hand, for $j\ge1$,
\[
\I^\tau_{(\rho^0_j,X^0)}(\rho_j,X_j)\le\I^\tau_{(\rho^0_j,X^0)}(\ov\rho,\ov X)\ \st{j\up\oo}\longto\ \I^\tau_{(\rho^0,X^0)}(\ov\rho,\ov X).
\]
We obtain that $\I^\tau_{(\rho^0,X^0)}(\rho,X)\le\I^\tau_{(\rho^0,X^0)}(\ov\rho,\ov X)$. Eventually, by Proposition~\ref{propOP3.1}, $(\ov\rho,\ov X)$ is the unique minimizer of $\I^\tau_{(\rho^0,X^0)}$ in $\Abb_M$ and we conclude that $(\rho,X)=(\ov\rho,\ov X)$. This ends the proof of the proposition.
\end{proof}

\subsection{$L^\oo$ and Lipschitz bounds}
\label{subsec.LooLipbounds}
~

This subsection is dedicated to the proof of some upper and lower bounds on $\rho_{|[0,X]}$ independent of $\tau$ and a bound on $\Lip(\rho_{|[0,X]})$ that degenerates in $\tau^{-\om}$ as $\tau$ goes to $0$.  These bounds, stated in the next two propositions, will propagate along the iterations of the minimizing movement scheme. The proof of these results are inspired by~\cite{FS21}.

\begin{prop}\label{prop.Loobounds}
Let $(\rho^0,X^0)\in\A$  and let $(\rho,X)\in\A$ be a minimizer for~\eqref{OP}. Recalling the conditions on $\rob$ and $\rot$ of~\eqref{rhominrhomax}, we have the following uniform bounds,
\be 
\label{Loobounds}
\rob\le\rho\le\rot\qquad\text{ on }[0,X].
\ee 
\end{prop}\medskip

\begin{proof}Let $(\rho^0,X^0)\in\A$  and $(\rho,X)\in\A$ as in the statement of the proposition. We set $M:=\MM(\rho)$ so that $(\rho,X)=\Ups^\tau((\rho^0,X^0),M)$.\\ 
Again, we consider the two cases $M\le\MM(\rho^0)$ and $M>\MM(\rho_0)$ separately.\medskip

\noindent
\textit{Step 1, $M\le\MM(\rho^0)$.} As in the proof of Proposition~\ref{propOP3.3a}, we define  $\ell_-\in[0,X^0)$ by
\[
\int_0^{\ell_-}\rho^0(y)\, dy =\MM(\rho^0) - M.
\] 
Then, denoting $T\in \Co^{1,1}([0,X], [\ell_-,X^0])$ %\footnote{See Proposition~\ref{propOP3.3a}.} 
the optimal mapping from $\rho\Leb$ to $\rho^0\Leb\restr[\ell_-,X^0]$, we have
\be
\label{rho=rho0T'j}
\rho(x)=\rho^0(T(x))T'(x)\qquad\text{for almost every }x\in[0,X].
\ee
Moreover, there exists a Kantorovitch potential $\Phi\in \Co^{2,1}([0,X],\R)$ such that 
\be\label{logrhoPhij}
(\log\rho)'(x) = -\dfrac{\Phi'(x)}\tau=\dfrac{T(x)-x}\tau\quad\text{on }[0,X].
\ee
Let $x^*\in[0,X]$ be an extremal point of $\rho$ on $[0,X]$.  We split the study in three cases.\smallskip

\noindent
\textit{Case 1.} If $x^*\in(0,X)$ and $x^*$ is a minimal point of $\pm\rho$, we have 
\[
(\log\rho)'(x^*)=0\qquad\text{and}\qquad\pm(\log\rho)''(x^*)\ge0.
\]
Differentiating~\eqref{logrhoPhij}, we get 
\[
\pm T'(x^*)=\pm1\pm\tau (\log\rho)''(x^*)\ge\pm1.
\]
By~\eqref{rho=rho0T'j}, we have $\rho(x^*)=T'(x^*)\rho^0(T(x^*))$ and we deduce from~\eqref{hyprho0j},
\be\label{proofLoo.1}
\begin{cases}
\min_{[0,X]}\rho =\rho(x^*)\ge\rho^0(T(x^*))\ge\rob &\text{if }x^*\text{ is a minimal point of }\rho\text{ in }(0,X),\\
\max_{[0,X]}\rho =\rho(x^*)\le\rho^0(T(x^*))\le\rot &\text{if }x^*\text{ is a maximal point of }\rho\text{ in }(0,X).
\end{cases}
\ee
\smallskip

\noindent
\textit{Case 2.} If $x^*=0$, we have by~\eqref{logrhoPhij}, 
\[
(\log\rho)'(0) =\dfrac{T(0)-0}\tau =\dfrac{\ell_-}\tau\ge 0.
\]
We consider two subcases:\\
$\circ$ If $M=\MM(\rho^0)$ then $\ell_-=0$, so $(\log\rho)'(0)=0$. Hence if the point $0$ is a minimal point of $\pm\rho$ we have $\pm(\log\rho)''(0)\ge0$ and, we deduce, as in Case~1, that~\eqref{proofLoo.1} holds true with $x^*=0$.\\
$\circ$ If on the contrary there holds $\ell_->0$ then $(\log\rho)'(0)>0$ and $0$ can only be a minimal point. We conclude that 
\[
 M=\MM(\rho^0)\implies
\begin{cases}
\min_{[0,X]}\rho\ge\rob&\text{if }0\text{ is a minimal point of }\rho\text{ in }[0,X],\\
\max_{[0,X]}\rho\le\rot &\text{if }0\text{ is a maximal point of }\rho\text{ in }[0,X].
\end{cases}
\]
\[
 M<\MM(\rho^0)\implies\  0\text{ is a minimal point of }\rho\text{ in }[0,X]\text{ and }\min_{[0,X]}\rho=\rho(0).
\]

\smallskip

\noindent
\textit{Case 3.} If $x^*=X$, we have by~\eqref{logrhoPhij} and the boundary condition~\eqref{CBenX} of Proposition~\ref{propOP3.2},
\[
(\log\rho)'(X) =\dfrac{T(X)-X}\tau =\dfrac{X^0-X}\tau=\dfrac{\ros-\rho(X)}\lambda.
\]
If $X$ is a minimal point of $\pm\rho$ then $\pm (\log\rho)'(X)\le0$, that is $\pm(\ros-\rho(X))\le 0$. We conclude that, 
\be\label{proofLoo.3}
\begin{cases}
\min_{[0,X]}\rho=\rho(X)\ge\ros\ge\rob &\text{if }X\text{ is a minimal point of }\rho\text{ in }[0,X],\\
\max_{[0,X]}\rho=\rho(X)\le\ros\le\rot &\text{if }X\text{ is a maximal point of }\rho\text{ in }[0,X].
\end{cases}
\ee\smallskip

We deduce from~\eqref{proofLoo.1}--\eqref{proofLoo.3} that if $M\le\MM(\rho^0)$ we have for  $j\ge1$, 
\[
\max_{[0,X]}\rho\le\rot\qquad\text{ and }\qquad
\begin{cases}
\min_{[0,X]}\rho\ge\rob&\text{ if }M=\MM(\rho),\\
\min_{[0,X]}\rho\ge\min(\rho(0),\rob)&\text{ if }M<\MM(\rho).
\end{cases}
\]
We see that~\eqref{Loobounds} holds true in the case $M=\MM(\rho^0)$.  Eventually, for the case $M<\MM(\rho^0)$, the boundary condition of Proposition~\ref{prop.BCx=0} (first case) yields $\rho(0)\ge\rop\ge\rob$. It follows that~\eqref{Loobounds} also holds true in the case $M<\MM(\rho^0)$.\medskip

\noindent
\textit{Step 2, $M>\MM(\rho^0)$.} As in previous proofs, we define $\ell_+>0$ by 
\[
\int_0^{\ell_+}\rho(x)\, dx=M-\MM(\rho^0)=:m.
\]
Setting $\mu^0:=m\delta_0+\rho^0\Leb$, the optimal mapping $T\in \Co^{1,1}([\ell_+,X],[0,X^0])$ between $\rho\Leb$ and $\mu^0$ satisfies
\[
\rho(x)=\rho^0(T(x))T'(x)\qquad\text{for almost every }x\in[\ell_+,X].
\]
Moreover, 
\be\label{logrhoPhij2}
(\log\rho)'(x) =\dfrac{T(x)-x}\tau\quad\text{on }[0,X].
\ee
Now let $x^*\in[0,X]$ be an extremal point of $\rho$ on $[0,X]$. Again, we consider successively three subcases.\smallskip

\noindent
\textit{Case 1.} If $x^*\in(\ell_+,X)$, the reasoning of Case~1 of the previous step ($M\le\MM(\rho^0)$) applies. We have similarly to~\eqref{proofLoo.1}, 
\be\label{proofLoo.4}
\begin{cases}
\min_{[0,X]}\rho\ge\rob &\text{if }x^*\text{ is a minimal point of }\rho\text{ in }(\ell_{j+},X),\\
\max_{[0,X]}\rho\le\rot &\text{if }x^*\text{ is a maximal point of }\rho\text{ in }(\ell_{j+},X).
\end{cases}
\ee
\smallskip

\noindent
\textit{Case 2.} We assume that $x^*\in[0,\ell_+]$. On this interval, we recall that we have (see~\eqref{exp-x2surtau}), 
\[
\rho(x)=\rho(0)\exp\lt(-\dfrac{x^2}{2\tau}\rt)\qquad\text{ on }[0,\ell_+].
\]
We see that $\rho'<0$ on $[0,\ell_+]$ and since $\rho$ is of class $\Co^1$ on $[0,X]$, $x^*$ is not a minimal point of $\rho$. We get 
that 
\be\label{proofLoo.5}
\text{If }x^*\in[0,\ell_+]\text{ then }x^*=0\text{ is a maximal point of }\rho\text{ on }[0,X]\text{ and }\max_{[0,X]}\rho =\rho(0).
\ee

\noindent
\textit{Case 3.} If $x^*=X$ we have~\eqref{proofLoo.3} exactly as in the case $M\le\MM(\rho^0)$.\smallskip

We deduce from~\eqref{proofLoo.3},~\eqref{proofLoo.4},~\eqref{proofLoo.5} that if $M>\MM(\rho^0)$, there holds for  $j\ge1$, 
\[
\min_{[0,X]}\rho\ge\rob ,\qquad\max_{[0,X]}\rho\le\max(\rho(0),\rot).
\]
Using the third case of Proposition~\ref{prop.BCx=0} we have $\rho(0)\le\rom\le\rot$. We conclude that~\eqref{Loobounds} also holds true in the case $M>\MM(\rho^0)$. We have covered every case.
\end{proof}

\begin{coro}
\label{coroborneM}
Let $(\rho^0,X^0)\in\A$  and let $(\rho,X)\in\A$ be a solution to the optimization problem~\eqref{OP}. There holds,
\be\label{bornMdefK}
\lt|\MM(\rho)-\MM(\rho^0)\rt|\le K\tau^{1-\om}\qquad\text{where}\quad K:=1 +\max\lt(\log\dfrac{\rot}{\rop},\log\dfrac{\rom}{\rob}\rt).
\ee
Besides, defining $\ell_+$ and $\ell_-$ as in~\eqref{ell+},~\eqref{ell-} and~\eqref{ell+ell-0}, we have 
\be\label{borneell+ell-}
0\le\ell_+,\ell_-\le A\tau^{1-\om}\qquad\text{where}\quad A:= K/\rob.
\ee
\end{coro}

\begin{proof}First, if $\MM(\rho)=\MM(\rho^0)$ then $\ell_+=\ell_-=0$ and there is nothing to prove. Next, if $\MM(\rho^0)>\MM(\rho)$, we have, by Proposition~\ref{prop.BCx=0}, 
\[
\log\rho(0)=\log\rop + p'_\tau(\MM(\rho^0)-\MM(\rho)).
\] 
By the previous proposition, $\rho(0)\le\rot$, hence 
\[
p'_\tau(\MM(\rho^0)-\MM(\rho))\le\log(\rot/\rop).
\]
Recalling that $p'_\tau(m)=(m\tau^{\om-1}- 1)_{+}\ge0$ for $m\ge0$ by definition~\eqref{2.def.penalty}. We get 
\[
\MM(\rho^0)-\MM(\rho)\le ( 1 +\log(\rot/\rop))\tau^{1-\om}.
\] 
This implies~\eqref{bornMdefK}. Finally, we have $\ell_+=0$ and writing 
\[
m=\int_0^{\ell_-}\rho^0(y)\, dy\ge\rob\ell_-,
\]
the estimate~\eqref{borneell+ell-} follows from~\eqref{bornMdefK}.

The inequalities in the case $\MM(\rho^0)<\MM(\rho)$ are proved in the same way using the third case of Proposition~\ref{prop.BCx=0} and $
\MM(\rho)-\MM(\rho^0)=\int_0^{\ell_+}\rho(x)\, dx\ge\rob\ell_+$.
\end{proof}

\begin{prop}\label{prop.Lipbounds}
Let $(\rho^0,X^0)\in\A$  and let $(\rho,X)\in\A$ be a minimizer for~\eqref{OP}. 
Denoting  $u^0:=(\log\rho^0)_{|[0,X^0]}\in \Co^{0,1}([0,X])$ and  $u:=(\log\rho)_{|[0,X]}\in \Co^{1,1}([0,X])$. We have 
\be\label{Lipbounds}
\Lip(u)\le\max\lt(A\tau^{-\om},B\rt),
\ee
where $A$ is the constant introduced in Corollary~\ref{coroborneM} and 
\[
B:=\max\lt(\Lip(u^0),\dfrac{\rot-\ros}\lambda,\dfrac{\ros-\rob}\lambda\rt).
\]
\end{prop}\medskip

\begin{proof} We consider a sequence $(\rho^0_j,X^0)\in\A$ satisfying the assumptions of Proposition~\ref{propOP3.3b} and such that furthermore $\rho^0_j$ is of class $\Co^{1,1}$ on $[0,X^0]$ for every $j\ge1$ and
\be\label{hyprho0j}
\min_{x\in [0,X^0]}\,\rho^0\le\rho^0_j(x)\le\max_{x\in [0,X^0]}\,\rho^0\qquad\mbox{for }x\in [0,X^0].
\ee
We assume furthermore that for $j\ge1$,\footnote{To obtain such sequence we may proceed as follows. We extend $\rho^0_{|[0,X^0]}$ on $\R$ as an even and $2X^0$-periodic function $\check\rho^0$. We convolve $\check\rho^0$ by a sequence $\zeta_j$ of even, smooth approximation of unity. Eventually, we restrict $\check\rho^0\ast\zeta_j$ to $[0,X^0]$ and we extend it by 0.}
\be\label{proofpropLip0}
\MM(\rho^0_j)=\MM(\rho^0)\qquad\text{and}\qquad\Lip(\rho^0_j)\le\Lip(\rho^0).
\ee
Setting  $M:=\MM(\rho)$, and denoting $(\rho_j,X_j):=\Ups^\tau((\rho^0_j,X^0),M)$, we get from Proposition~\ref{propOP3.3a}, that (with obvious notation) $\rho_j$ is of class $\Co^{1,1}$ on $[0,X_j]$ and of class $\Co^{3,1}$ on the intervals $[0,\ell_{j+}]$ and $[\ell_{j+},X_j]$. Eventually, the sequence $(\rho_j,X_j)$ converges to $(\rho,X)$ in the sense of~\eqref{propOP3.3_1}, that is,
\be\label{proofpropLip0.5}
\qquad X_j\ \st{j\up\oo}\longto\ X\qquad\text{and}\qquad\rho_j\ \st{j\up\oo}\longto\ \rho\text{ in }L^1(\R_+).\medskip
\ee

For $j\ge1$, we denote 
\[
u_j:=\log({\rho_j}_{|[0,X_j]})\in \Co^{1,1}([0,X_j]).
\]
 In the sequel, we study the extrema of $u_j'$ to obtain a bound on $\Lip(u_j)$ and then send $j$ to $+\oo$. The proof follows the same steps and subcases as the proof of Proposition~\ref{prop.Loobounds}.\medskip

\noindent
\textit{Step 1, $M\le\MM(\rho^0)$.} As in the proof of Proposition~\ref{propOP3.3a}, we define  $\ell_-\in[0,X^0)$ by
\[
\int_0^{\ell_-}\rho^0(y)\, dy =\MM(\rho^0) - M.
\] 
Then, denoting $T_j\in \Co^{2,1}([0,X_j], [\ell_-,X^0])$ %\footnote{See Proposition~\ref{propOP3.3a}.} 
the optimal mapping from $\rho_j\Leb$ to $\rho^0_j\Leb\restr[\ell_-,X^0]$, we have 
\be
\label{rho=rho0T'j3}
\rho_j(x)=\rho^0_j(T_j(x))T_j'(x)\qquad\text{for almost every }x\in[0,X_j].
\ee
Moreover, there holds
\be\label{logrhoPhij3}
u_j'(x)=(\log\rho_j)'(x) =\dfrac{T_j(x)-x}\tau\quad\text{on }[0,X_j].
\ee
 Let $x^*\in[0,X_j]$ be a maximal point of $|u_j'|$ on $[0,X_j]$.  We consider three cases.\smallskip

\noindent
\textit{Case 1.} If $x^*\in(0,X_j)$ we have  $(\log\rho_j)''(x^*)=0$ and $\pm(\log\rho_j)'''(x^*)\ge0$ if $x^*$ is a minimal point of $\pm u_j'$. 
Differentiating~\eqref{logrhoPhij3}, we get $T_j'(x^*)=1$ and differentiating again we obtain $\pm T_j''(x^*)\ge0$. Then, using~\eqref{rho=rho0T'j3}, we compute for $x\in (0,X_j)$, 
\[
u_j'(x)=\dfrac{\rho'_j(x)}{\rho_j(x)}= T'_j(x)\lt[\dfrac{{\rho^0_j}'}{\rho^0_j}\rt](T_j(x)) +\dfrac{T_j''(x)}{T'_j(x)}= T'_j(x) {u^0_j}'(T_j(x)) +\dfrac{T_j''(x)}{T'_j(x)},
\]
where  $u^0_j:=\log\rho^0_j$.  At $x=x^*$, using $T_j'(x^*)=1$ and $\pm T_j''(x^*)\ge0$, this leads to
\[
\pm u_j'(x^*)\ge\pm (u_j^0)'(T_j(x^*)).
\]
Since, $x^*$ is a maximal point of $|u_j'|$, we have $\mp u_j'(x^*)= |u'_j(x^*)|\ge 0$. We conclude that 
\be\label{proofLip.1}
\Lip(u_j)=|u'_j(x^*)|\le\Lip (u_j^0).%\quad\text{ if }x^*\text{ is a maximal point of }|u'_j|\text{ in }(0,X_j).
\ee

\noindent
\textit{Case 2.} If $x^*=0$, we have by~\eqref{logrhoPhij3}, 
\[%\be\label{proofLip.2}
\Lip(u_j)=u_j'(0)=(\log\rho_j)'(0) =\dfrac{\ell_-}\tau\ge 0.
\]%\ee\smallskip

\noindent
\textit{Case 3.} If $x^*=X_j$, we have by~\eqref{logrhoPhij3} and the boundary condition~\eqref{CBenX} of Proposition~\ref{propOP3.2},
\be\label{proofLip.3}
u_j'(X_j) =(\log\rho_j)'(X_j) =\dfrac{T_j(X_j)-X_j}\tau =\dfrac{X^0-X_j}\tau=\dfrac{\ros-\rho_j(X_j)}\lambda.
\ee\smallskip
We deduce from~\eqref{proofLip.1}--\eqref{proofLip.3} that if $M\le\MM(\rho^0)$ we have for  $j\ge1$, 
\[
\Lip(u_j)\le\max\lt(\Lip (u_j^0),\dfrac{\ell_-}\tau,\dfrac{|\ros-\rho_j(X_j)|}\lambda\rt).
\]
Sending $j$ to $+\oo$, we have by~\eqref{proofpropLip0} and~\eqref{proofpropLip0.5},
\[
\Lip(u)\le\max\lt(\Lip (u^0),\dfrac{\ell_-}\tau,\dfrac{|\ros-\rho(X)|}\lambda\rt).
\]
Eventually, by Proposition~\ref{prop.Loobounds} (and $\rob\le\ros\le\rot$), we have $|\ros-\rho(X)|\le\max(\rot-\ros,\ros-\rob)$. Then, using the estimate~$\ell_-\le A\tau^{1-\om}$ of Corollary~\ref{coroborneM}, we conclude that~\eqref{Lipbounds} holds in the case $M\le\MM(\rho^0)$.\medskip

\noindent
\textit{Step 2, $M>\MM(\rho^0)$.} For $j\ge1$ we define $\ell_{j+}>0$ by 
\[
\int_0^{\ell_{j+}}\rho_j(x)\, dx=M-\MM(\rho^0)=:m.
\]
We set $\mu^0_j:=m\delta_0+\rho^0_j\Leb$ and denote  $T_j\in \Co^{1,1}([0,X_j],[0,X^0])\cap \Co^{2,1}([\ell_{j+},X_j],[0,X^0])$ the optimal mapping between $\rho_j\Leb$ and $\mu_j^0$. There hold
\[
\rho_j(x)=\rho^0_j(T_j(x))T_j'(x)\qquad\text{for almost every }x\in[\ell_{j+},X_j]
\]
and~\eqref{logrhoPhij2}.

Let $x^*\in[0,X_j]$ be a maximal point of $|u_j'|$ on $[0,X_j]$. We consider again three subcases.\smallskip

\noindent
\textit{Case 1.} If $x^*\in(\ell_{j+},X_j)$, arguing as in the first step, we obtain that 
\be\label{proofLip.35}
\Lip(u_j)\le\Lip (u_j^0).%\quad\text{ if }x^*\text{ is a maximal point of }|u'_j|\text{ in }(\ell_{j+},X_j).
\ee
\smallskip

\noindent
\textit{Case 2.} We assume that $x^*\in[0,\ell_{j+}]$. On this interval, by~\eqref{logrhoPhij2}, $u_j'(x)= - x/\tau$ and consequently,
\be\label{proofLip.4}
0\ge u_j'(x^*)\ge-\dfrac{\ell_{j+}}\tau.
\ee\smallskip

\noindent
\textit{Case 3.} If $x^*=X_j$ we obtain~\eqref{proofLip.3} exactly as in the Case~3 of Step~1.\medskip

We deduce from~\eqref{proofLip.35},~\eqref{proofLip.4} and~\eqref{proofLip.3} that in the case $M<\MM(\rho^0)$ there holds for $j\ge1$, 
\[
\Lip(u_j)\le\max\lt(\Lip (u_j^0),\dfrac{\ell_{j+}}\tau,\dfrac{|\ros-\rho_j(X_j)|}\lambda\rt).
\]
Sending $j$ to $+\oo$ we notice that $\ell_{j+}\to\ell_+$ and using $\ell_+\le A\tau^{-\om}$, we conclude as in Step~1 that~\eqref{Lipbounds} holds true. This ends the proof of the proposition.
\end{proof}

\section{Existence of weak solutions to~\eqref{P}--\eqref{P.bord} }
\label{Sec.Th2}

In this section we prove Theorem~\ref{th.main}.  We assume that the hypotheses (H1)--(H2) of the theorem hold. Consequently,
\begin{enumerate}[($\circ$)]
\item  The parameters $\rop\ge\rom>0$, $\ros>0$ and $\lambda>0$ are given.
\item $\theta>0$ and $\beta\ge0$ are defined by~\eqref{def_betatheta}. 
\item The initial data $(\rho^0,X^0)$ is also fixed and satisfies (H2).%, equivalently $(\rho^0,X^0)\in\A$.
\end{enumerate}

\subsection{Definition of a JKO minimizing scheme}\label{sec.JKOscheme}~

It is convenient  to introduce here an arbitrary and temporary maximal time  
\be\label{Tmax}
\TT>0,
\ee
for the existence of the solutions to~\eqref{P}--\eqref{P.bord}.  $\TT$ will be sent to infinity at the end of the proof. 

In this subsection, we fix a time step $0<\tau\le 1$. Let us try to  build a finite sequence $(\rho^n,X^n)\in\A$ recursively, by 
\be
\label{JKOscheme0}
(\rho^{n+1},X^{n+1})\in\mathrm{argmin}\lt\{\J^\tau_{(\rho^n,X^n)}(\rho,X): (\rho,X)\in\A\rt\}.
\ee
The pair $(\rho^n,X^n)$ is thought of as the solution at time $n\tau$ of a time discretization of~\eqref{P}--\eqref{P.bord}.\\
By Theorems~\ref{theo.ExistenceJKO.1} and~\ref{theo.ExistenceJKO.2}, the existence of $(\rho^{n+1},X^{n+1})$ is ensured (see condition~\eqref{X0large}) if
\be\label{Xnlarge}
 X^n>\dfrac{\tau R}\lambda +\sqrt{\dfrac{\tau^2 R^2}{\lambda^2}+\dfrac{2\tau}\lambda P(\rho^n,X^n)}, 
\ee
 where we recall:
\[
R=\ros-\sqrt{\rop\rom}\qquad\text{and}\qquad P(\rho^n,X^n)=\PSI(\rho^n,X^n)-RX^n.\medskip
\]

The finite sequence $(\rho^n,X^n)$ indexed by $0\le n\le N_\tau$ is then built recursively as follows:
\be\label{JKOscheme1}
\begin{cases}
 (\rho^{n+1},X^{n+1})\text{ is defined by}~\eqref{JKOscheme0}&\text{as soon as~\eqref{Xnlarge} holds true and }(n+1)\tau\le\TT,\\
\quad N_\tau:=n\quad\text{ (we stop)}&\text{whenever~\eqref{Xnlarge} fails or }(n+1)\tau>\TT.
\end{cases}
\ee

 Observe that by Remark~\ref{RemX0large}, $N_\tau\ge 1$ for $\tau>0$ small enough.\medskip
 
Then, we set $T_\tau:=N_\tau\tau$ and we define a function $t\in [0,T_\tau]\mapsto (\rho_\tau,X_\tau)\in\A$ by
\be\label{def.rhotauXtau}
\begin{array}{rl}
&(\rho_\tau, X_\tau)(0):=(\rho^0,X^0)\smallskip\\
\text{and }\qquad&(\rho_\tau, X_\tau)(t):=(\rho^{n+1},X^{n+1})\quad\text{ for }t\in(n\tau,(n+1)\tau],\ 0\le n\le N_\tau-1.
\end{array}
\ee

The rest of the Section is organized as follows. In Subsection~\ref{subsec.apriori} we establish some estimates on the sequence $(\rho^n,X^n)$, which are uniform with respect to $\tau$. They allow us, in Section~\ref{Sec:Compacite}, to deduce some compactness properties on $(\rho_\tau,X_\tau)$. Sending $\tau$ to $0$, we obtain, up to extraction, a limit function $t\in[0,T)\mapsto (\rho(t),X(t))$. We show in Subsections~\ref{subsec.proof.wealsol} and~\ref{subsec.var} that this latter is a weak solution to~\eqref{P}--\eqref{P.bord} in the sense of Definition~\ref{def.sol.faible}.

Let us first introduce some notation.  For $0\le n\le N_\tau$, we set 
\be\label{Mn}
M^n=\MM(\rho^n),
\ee 
and for $0\le n\le N_\tau-1$, we define,
\be\label{mn}
m^n:=M^{n+1}-M^n,
\ee
and $\ell^n_+,\,\ell_-^n\ge0$ by 
\begin{align}
\label{ell+n}
\int_0^{\ell^n_+}\rho^{n+1}\, = m^n\text{ when }m^n>0,\qquad\ell_+^n=0\text{ when }m^n\le0,\\
\label{ell-n}
\int_0^{\ell^n_-}\rho^n\, = -m^n\text{ when }m^n<0,\qquad\ell^n_-=0\text{ when }m^n\ge0.
\end{align}
Moreover, for a $L^1\loc$  function $w$ depending on $x\ge0$ and possibly of $t$, if the distribution $\pt_x w$ is a measure, we recall that $w'$ denotes the absolute continuous part of $\pt_x w$. In particular, for $(\rho,X)\in\A$, we have the following identity, in the space  $\mathcal{M}(\R_+)$,
\[
\pt_x\rho =\rho'\Leb\restr{\R_+}-\rho(X)\delta_X. 
\]

\subsection{Uniform estimates}\label{subsec.apriori}~

Let us first collect the properties satisfied by each $(\rho^n,X^n)$ which are direct consequences of the results of the previous section.

\begin{prop}\label{prop.rhonXn}
We have for $0\le n\le N_\tau-1$:
\begin{enumerate}[(a)]
\item
For every $\psi\in\Co^\oo_c(\R_+)$, there holds, 
\[
\int_{\R_+}\dfrac{\rho^{n+1}(x)-\rho^n(x)}{\tau}\,\psi(x)\, dx-\dfrac{M^{n+1}-M^n}{\tau}\psi(0)+\int_0^{X^{n+1}}(\rho^{n+1})'(x)\,\psi'(x)\, dx := Q^n_\tau(\psi),
\]
where
\[
\lt|Q^n_\tau(\psi)\rt|\le\dfrac{\|\psi''\|_{L^\oo(\R_+)}}{2\tau}\,\Wass^2(\rho^{n+1},\rho^n) +\dfrac{\|\psi'\|_{L^\oo(\R_+)}}{\tau}\,\lt(\int_0^{\ell^n_-}y\,\rho^n(y)\, dy\rt).
\] 
\item At the interface $x=0$, the following discrete version of the boundary condition~\eqref{P.bord} holds,
\[
\begin{cases}
\log\rho^{n+1}(0)=\log\rop + p'_\tau(-m^n)&\text{if }m^n<0,\smallskip\\
\qquad\rom\le\rho^{n+1}(0)\le\rop&\text{if }m^n=0,\smallskip\\
\log\rho^{n+1}(0)=\log\rom - p'_\tau(m^n)&\text{if }m^n>0.
\end{cases}
\]
\item At the interface $x=X^{n+1}$, the following discrete version of~\eqref{P.e} holds, 
\[
\lambda\dfrac{X^{n+1}-X^n}\tau=\rho^{n+1}(X^{n+1})-\ros.
\]
\item We have $\rob\le\rho^{n+1}\le\rot$ on $[0,X^{n+1}]$.
\item
Denoting
\[
u^0:=\log({\rho^0}_{|[0,X^0]})\qquad\text{and}\qquad u^{n+1}:=\log({\rho^{n+1}}_{|[0,X^{n+1}]}),
\]
and recalling the definitions:
\[
K=1 +\max\lt(\log\dfrac{\rot}{\rop},\log\dfrac{\rom}{\rob}\rt),\qquad A=\dfrac K{\rob},\qquad B=\max\lt(\Lip(u^0),\dfrac{\rot-\ros}\lambda,\dfrac{\ros-\rob}\lambda\rt),
\]
there hold,
\[
\Lip(u^{n+1})\le\max\lt(A\tau^{-\om},B\rt),\qquad  |m^n|\le K\tau^{1-\om}\qquad\text{and}\qquad 0\le\ell^n_+,\ell^n_-\le A\tau^{1-\om}.
\]
\end{enumerate}
\end{prop}

\begin{proof}
All the points follow by applying inductively some results of Section~2. Namely, (a) is Theorem~\ref{theo.EL},~(b) is Proposition~\ref{prop.BCx=0} and (c) is Proposition~\ref{propOP3.2}. Next, we get point~(d) by iterating the lower and upper bounds of Proposition~\ref{prop.Loobounds}. Indeed, setting (as in~\eqref{rhominrhomax}),
\[
\rob^{n}:=\min\lt\{\min_{[0,X^n]}\rho^n,\ros,\rop\rt\}\ \le\ \rot^n: =\max\lt\{\max_{[0,X^n]}\rho^n,\ros,\rom\rt\},
\]
we show by induction that,
\[
\rob=\rob^0\le\cdots\le\rob^{n-1}\le\rob^{n}\le\min_{[0,X^{n+1}]}\rho^{n+1}\le\max_{[0,X^{n+1}]}\rho^{n+1}\le\rot^{n}\le\rot^{n-1}\le\cdots\le\rot^0=\rot.
\]
The first inequality of point~(e) comes by iterating Proposition~\ref{prop.Lipbounds}. We have by induction, that
\[
B_n:=\max\lt(\Lip(u^n),\dfrac{\rot^n-\ros}\lambda,\dfrac{\ros-\rob^n}\lambda\rt)\le\max\lt(A\tau^{-\om},B\rt).
\]
The second and third inequalities of~(e) follow from Corollary~\ref{coroborneM}.
\end{proof}

We now establish some properties of the sequences $(\rho^n)$, $(X^n)$ and $(M^n)$.

\begin{prop}\label{prop.Ttau}
There exist positive constants $C_c$, $\hat \TT$, $\hat\tau$ and $C_e$ only depending on $(\rho^0,X^0)$, $\TT$ and on the parameters of the model such that:
\begin{enumerate}[(a)]
\item The sequence $(\PSI(\rho^n,X^n))$ is nonincreasing,
\item There holds, for $0\le n\le N_\tau-1$,
\[
\dfrac{\lt|X^{n+1}-X^n\rt|}\tau\le\max\lt(\dfrac{\rot-\ros}\lambda,\dfrac{\ros-\rob}\lambda\rt)=:C_b.
\]
\item There holds
\[
T_\tau  = N_\tau\,\tau\le\TT -\tau\implies X^{N_\tau}\le C_c\sqrt\tau.
\]
\item There holds
\[ 
0<\tau\le\hat\tau\implies T_\tau\ge\hat T.
\]
\item There holds
\[
\PSI(\rho^n,X^n)\ge C_e\qquad\text{ for}\quad 0\le n\le N_\tau. 
\]
\item There holds
\be\label{e.a.p}
\sum_{n=0}^{N_\tau-1}\left(\dfrac{\Wass^2(\rho^{n+1},\rho^n)}{2\tau}+\lambda\dfrac{(X^{n+1}-X^n)^2}{2\tau}+\te |m^n| + p_{\tau}(m^n)\right)
\le\PSI(\rho^0,X^0)-C_e.
\ee
\end{enumerate}
\end{prop}

\medskip

\begin{proof}~

\noindent 
\textit{(a)} Let $0\le n\le N_\tau-1$. By optimality of $(\rho^{n+1},X^{n+1})$, we have $\J^\tau_{(\rho^{n},X^{n})}(\rho^{n+1},X^{n+1})\le\J^\tau_{(\rho^{n},X^{n})}(\rho^{n},X^{n})$. Using~\eqref{Jtau2}, this is equivalent to 
\be
\label{proof.prop.Ttau.1}
\dfrac{\Wass^2(\rho^{n+1},\rho^n)}{2\tau}+\lambda\dfrac{(X^{n+1}-X^n)^2}{2\tau}+\te |m^n| + p_{\tau}(m^n)
\le\PSI(\rho^n,X^n)-\PSI(\rho^{n+1},X^{n+1}).
\ee
 Observing that the left hand side is nonnegative, we obtain that the sequence $(\PSI(\rho^n,X^n))$, indexed by $0\le n\le N_\tau$, is nonincreasing. 
\medskip

\noindent 
\textit{(b)} This is a direct application of points (c) and~(d) of Proposition~\ref{prop.rhonXn}.
\medskip

\noindent 
\textit{(c)} Let us assume that $T_\tau\le\TT -\tau$ so that~\eqref{Xnlarge} fails for $n=N_\tau$. We have,
\begin{align}
\nonumber
X^{N_\tau}
\le\dfrac{\tau R}\lambda +\sqrt{\dfrac{\tau^2 R^2}{\lambda^2}+\dfrac{2\tau}\lambda P(\rho^{N_\tau},X^{N_\tau})}&=\dfrac{\tau R}\lambda +\sqrt{\dfrac{\tau^2 R^2}{\lambda^2}-\dfrac{2\tau RX^{N_\tau}}\lambda +\dfrac{2\tau}\lambda\PSI(\rho^{N_\tau},X^{N_\tau})}\\
\nonumber
&\st{(a)}\le\dfrac{\tau R}\lambda +\sqrt{\dfrac{\tau^2 R^2}{\lambda^2}-\dfrac{2\tau RX^{N_\tau}}\lambda +\dfrac{2\tau}\lambda\PSI(\rho^0,X^0)}\\
\label{proof.prop.Ttau.2}
&\le\dfrac{2\sqrt\tau |R|}\lambda +\sqrt{\dfrac{2\tau |R|}\lambda}\sqrt{X^N_\tau}+\sqrt{\dfrac{2\tau}\lambda\lt|\PSI(\rho^0,X^0)\rt|}.
\end{align}
To get the last inequality, we have used $\tau\le\sqrt\tau$ (since $0<\tau\le1$) and the fact that the function $r\in\R_+\mapsto\sqrt r$ is increasing and subadditive.\\
 Eventually, by~(b) we have $X^{N_\tau}\le X^0+C_b\TT$. Using this to estimate $\sqrt{X^{N_\tau}}$ in the right hand side of~\eqref{proof.prop.Ttau.2} ends the proof of point~(c). 
\medskip

\noindent 
\textit{(d)} Let $0<\tau\le\min(1,\TT/2)$.   If $T_\tau\ge\TT-\tau$ then $T_\tau\ge\TT/2$ and  we are done as soon as $\hat T\le\TT/2.$ We now assume $T_\tau<\TT-\tau$. By~(c), $X^{N^\tau}\le C_c\sqrt\tau$ and with~(b), we have 
\[
C_b T_\tau=C_b N_\tau\tau\ge X^0-X^{N_\tau}\ge X^0-C_c\sqrt{\tau}.
\]
Assuming $\sqrt\tau\le X_0/(2C_c)$, we get $C_b T_\tau\ge X^0/2$, that is $T_\tau\ge X^0/(2C_b)$. We conclude that~(d) holds true with
\[
\hat T:=\min(\TT/2,X^0/(2C_b))\quad\text{ and }\quad\hat\tau:=\min(1,\TT/2,(X_0/(2C_c))^2).\medskip
\]

\noindent 
\textit{(e)} Let $0\le n\le N_\tau$, we have, using $\min_{\R_+}f=-\sqrt{\rop\rom}$, 
\[
\PSI(\rho^n,X^n)
=\int_0^{X^n}f(\rho^n(x))\,dx +\ros X^n
\ge\lt(\ros-\sqrt{\rop\rom}\rt)X^n\st{(b)}\ge  -\lt(\sqrt{\rop\rom}-\ros\rt)_+ (X^0+C_b\TT).
\]
This proves (e) with $C_e:=-\lt(\sqrt{\rop\rom}-\ros\rt)_+ (X^0+C_b\TT)$.
\medskip

\noindent 
\textit{(f)} Summing~\eqref{proof.prop.Ttau.1} for $n=0,\dots,N_\tau-1$, we obtain that the left hand side of~\eqref{e.a.p} is bounded by 
$\PSI(\rho^0,X^0)-\PSI(\rho^{N_\tau},X^{N_\tau})\le\PSI(\rho^0,X^0)-C_e$. This proves~(f). 
\end{proof}

Now, in addition to the functions $\rho_\tau$ and $X_\tau$ introduced in~\eqref{def.rhotauXtau} we define the functions $\tilde{X}_\tau$ and $M_\tau$ as follows. 
We set for $0\le n\le N_\tau-1$ and $t\in [n\tau, (n+1)\tau]$,
\be\label{def.tildeXtauMtau}
\begin{array}{rl}
\tilde{X}_\tau(t) &:=\dfrac{t-n\tau}{\tau}X^{n+1}+\dfrac{(n+1)\tau-t}{\tau}\, X^n,\medskip\\
M_\tau(t) &:=\dfrac{t-n\tau}{\tau}M^{n+1}+\dfrac{(n+1)\tau-t}{\tau}\, M^n.
\end{array}
\ee
In other words, $\tilde{X}_\tau$ and $M_\tau$ are the continuous piecewise affine interpolations of the sequences $(X^n)$ and $(M^n)$ with respect to the grid $\{t_n=n\tau\}_{0\le n\le N_\tau}$.\medskip

Finally, for $I$ a real interval and a mapping  $w:t \in I \mapsto w(t)$, we introduce for $\tau>0$ the shift operator $\sigma_\tau$ defined by
\be\label{sigmatau}
\sigma_{-\tau}w(t)= w(t-\tau)\quad\text{ for }t\in I+\tau.
\ee

\begin{prop}\label{prop_apriori}
Let  $\tau\in (0,\hat\tau]$ where the constant $\hat\tau$ is the one of Proposition~\ref{prop.Ttau}(d). There exists a constant $C> 0$ only depending on $(\rho^0,X^0)$, $\TT$ and on the parameters of the model such that:
\be
\label{bounds.XtauMtau}
\tilde X_\tau(0)=X_0,\qquad\Lip(\tilde X_\tau)\le C_b,
\qquad\|M_\tau\|_{W^{1,1}(0,T_\tau)}\le C,\medskip
\ee
\be
\label{bounds.rhotau}
\rob\le\rho_\tau(t)\le\rot\text{ in }[0,X^{n+1}],\quad\supp\rho_\tau(t) = [0,X^{n+1}]\quad\mbox{ for }n\tau< t\le (n+1)\tau,\  n< N_\tau,
\ee
\begin{align}
\label{rho_H1}
\int_0^{T_\tau}\lt\|\rho'_\tau(t)\rt\|_{L^2(\R_+)}^2\,dt\le C,\\
\label{time_translate}
\dfrac1\tau\int_\tau^{T_\tau}\lt\|\rho_\tau(t)-\sigma_{-\tau}\rho_\tau(t)\rt\|_{H^*}\, dt\le C,
\end{align}
where $H^*$ denotes the dual space of $H^1(\R_+)$. Notice that in~\eqref{bounds.XtauMtau}, $C_b=C_b(\ros,\rot,\rob,\lambda)$ is the constant of Proposition~\ref{prop.Ttau}(b). Eventually, for $t \in [0,T_\tau] $,
\be\label{bound.rho'oo}
\|\rho_\tau'(t)\|_{L^\oo(\R_+)}\le C\tau^{-\om}. 
\ee
\end{prop}

\begin{proof} In the proof $C$ is a positive constant  that may change from line to line and that only depends on $(\rho^0,X^0)$, $\TT$ and on the parameters of the model.\medskip

\noindent
\textit{Step 1. Proof of~\eqref{bounds.XtauMtau}--\eqref{rho_H1}.} 
First the identity $\tilde X_\tau(0)=X_0$ holds by definition. Next, by construction, the functions ${\tilde X}_\tau$ and $M_\tau$ are continuous on $[0,T_\tau]$ and for $0\le n\le N_\tau-1$, we have $\dot{\tilde X}_\tau=(X^{n+1}-X^n)/\tau$ and $\dot{M}_\tau=m^n$ on $(n\tau,(n+1)\tau)$. With this in mind, we see that the first estimate of~\eqref{bounds.XtauMtau} follows from Proposition~\ref{prop.Ttau}(b). Besides, by Proposition~\ref{prop.Ttau}(f), we have
\[
\te\sum_0^{N_\tau-1}|m^n|\le\PSI(\rho^0,X^0)-C_e\le C.
\]
Since $\te>0$, this yields the second estimate of~\eqref{bounds.XtauMtau}.\\
Next,~\eqref{bounds.rhotau} is a consequence of Proposition~\ref{prop.rhonXn}(d) and of the definition of $X^{n+1}=x_{\rho^{n+1}}$.\\
Then, by construction, there holds
\[
\int_0^{T_\tau}\lt\|\rho'_\tau(t)\rt\|_{L^2(\R_+)}^2\,dt =\sum_{n=0}^{N\tau -1}\tau\int_0^{X^{n+1}}\lt| (\rho^{n+1})'\rt|^2 (x)\, dx .
\]
Applying~\eqref{estim_grad_L2} from Theorem~\ref{theo.EL} with $(\rho^n,X^n)$ in place of $(\rho^0,X^0)$, we obtain 
\[
\int_0^{T_\tau}\lt\|\rho'_\tau(t)\rt\|_{L^2(\R_+)}^2\,dt
\le\sum_{n=0}^{N\tau -1}\|\rho^{n+1}\|_{L^\oo(\R_+)}\,\dfrac{\Wass^2(\rho^n,\rho^{n+1})}\tau\,\le\rot\sum_{n=0}^{N\tau -1}\dfrac{\Wass^2(\rho^n,\rho^{n+1})}\tau.
\]
Using the estimate of Proposition~\ref{prop.Ttau}(f), we see that the right hand side is uniformly bounded, which proves~\eqref{rho_H1}. 
\medskip

\noindent
\textit{Step 2. Proof of~\eqref{time_translate}.} Let $0\le n\le N_\tau-1$ be fixed and $\phi\in H^1(\R_+)$. We consider the quantity
\[
I^n(\phi):=\int_{\R_+}(\rho^{n+1}(x)-\rho^n(x))\,\phi(x)\, dx.
\]
Let us assume that $M^{n+1}\ge M^n$ (the proof in the case $M^{n+1}<M^n$ is similar). We have
\[
I^n(\phi) =\int_0^{\ell^{n+1}_+}\rho^{n+1}(x)\phi(x)\, dx +\int_{\ell^{n+1}_+}^{+\oo}(\phi(x)-\phi(T^{n+1}_+(x)))\rho^{n+1}(x)\, dx=:I^n_1(\phi)+I^n_2(\phi),
\]
where the optimal transport map $T^{n+1}_+$ is defined as in~\eqref{T+} of Subsection~\ref{sec.not}. 

Let us estimate the first term. We have, 
\[
|I^n_1(\phi)|\le\|\phi\|_{L^\oo(\R_+)}\,\lt(\int_0^{\ell^{n+1}_+}\rho^{n+1}(x)\, dx\rt)\le\|\phi\|_{H^1(\R_+)}m^n,
\]
where we have used the Sobolev embedding $H^1(\R_+)\hookrightarrow L^\infty(\R_+)$ (whose operator norm is 1).

For the second term, using Fubini, $\rho^{n+1}\le\rot$ and  $\rho^n(y)\ge\rob$ for $y\in [0,X^n]$, and then the Cauchy-Schwarz inequality, we compute, 
\begin{align*}
|I^n_2(\phi)|  =\lt |\int_{\ell^{n+1}_+}^{X^{n+1}}\rho^{n+1}(x)\int_{T^{n+1}_+(x)}^ x\phi'(s)\,ds\, dx\rt |  &=\lt |\int_{0}^{X^n}\phi'(s)\int_{s}^{\lt(T_+^{n+1}\rt)^{-1}(s)}\rho^{n+1}(x)\,dx\, ds\rt |\\
 &\le\dfrac{\rot}{\sqrt{\rob}}\int_{0}^{X^n}\sqrt{\rho^n(s)}\lt |s-\lt(T_+^{n+1}\rt)^{-1}(s)\rt|\,\lt|\phi'(s)\rt|\,ds\\
  &\le\dfrac{\rot}{\sqrt{\rob}}\,\|\phi'\|_{L^2(\R_+)}\Wass(\rho^{n+1},\rho^n).
\end{align*}
With the previous estimate, we obtain 
\[
|I^n(\phi)|\le\|\phi\|_{H^1(\R_+)}m^n+\dfrac{\rot}{\sqrt{\rob}}\,\|\phi'\|_{L^2(\R_+)}\Wass(\rho^{n+1},\rho^n).
\]
Since this holds for every $\phi\in H^1(\R)$, we get for $0\le n\le N_\tau-1$,
\[
\|\rho^{n+1}-\rho^n\|_{H^*}\le |m^n|+\dfrac{\rot}{\sqrt{\rob}}\Wass(\rho^{n+1},\rho^n).
\]
Summing over $n$ and using  the Cauchy-Schwarz inequality and $N_\tau\tau\le\TT$, we get
\begin{align*}
\dfrac1\tau\int_\tau^{T_\tau}\lt\|\rho_\tau(t)-\sigma_{-\tau}\rho_\tau(t)\rt\|_{H^*}\, dt  
&\le\sum_{n=0}^{N_\tau-1}|m^n|+\dfrac{\rot}{\sqrt{\rob}}\sum_{n=0}^{N_\tau-1}\Wass(\rho^{n+1},\rho^n)\\
&\le\sum_{n=0}^{N_\tau-1}|m^n|+\dfrac{\rot}{\sqrt{\rob}}\sqrt{N_\tau}\left(\sum_{n=0}^{N_\tau-1}\Wass^2(\rho^{n+1},\rho^n)\right)^{1/2}\\
&\le\sum_{n=0}^{N_\tau-1}|m^n|+\dfrac{\rot}{\sqrt{\rob}}\sqrt{\TT}\left(\sum_{n=0}^{N_\tau-1}\dfrac{\Wass^2(\rho^{n+1},\rho^n)}\tau\right)^{1/2}.
\end{align*}
By Proposition~\ref{prop.Ttau}(f) the right hand side is uniformly bounded by a constant independent from $\tau$. This proves~\eqref{time_translate}.\\
Eventually, by proposition~\ref{prop.rhonXn}(e) and with the corresponding notation, we have for $\tau>0$ small enough and $0\le t\le T_\tau$,
\[
\Lip(\rho_\tau(t))\le \|\rho_\tau(t)\|_{L^\oo(\R_+)} \Lip\lt(\log({\rho_\tau(t)}_{|[0,X_\tau(t)]})\rt) \le \rot \max \lt(A\tau^{-\om},B\rt)\le C\tau^{-\om}.
\]
This proves~\eqref{bound.rho'oo} and the proposition.
\end{proof}

\subsection{Compactness of $\{(\rho_\tau,X_\tau)\}$}
\label{Sec:Compacite}

In this section we establish the existence of some time $T>0$ and some function $t\in[0,T]\mapsto (\rho(t),X(t))$ which is the limit up to extraction of $(\rho_\tau,\tilde X_\tau)$. Indeed, the previous bounds on the time-discrete solutions lead to the following compactness results.

\begin{prop}\label{Prop_compacite1}~
There exist ${\ov T},\, X,\, \rho $ such that
\begin{enumerate}[(a)]
\item ${\ov T}>0$ with $\TT\ge {\ov T}\ge\min(X_0/C_b,\TT)$;\smallskip
\item $X\in\Lip((0,\ov T),(0,+\oo))$ such that 
\[
X(0)=X_0,\qquad\Lip(X)\le C_b\quad\text{ with }
\begin{cases}
 X_{|[0,\ov T[} >0 \text{ and } X(\ov T)=0, & \text{ if }\ov T<\TT \\
 X_{|[0,\ov T]} >0, & \text{ else; }
\end{cases}
\]
\item $\rho\in L^2(\R_+\t(0,\ov T))$  such that with $D_{\ov T}$ defined as in~\eqref{def_D_T}, there hold,
\be\label{rho_lim}
\supp\rho=\ov{D_{\ov T}},\qquad\rob\le\rho\le\rot\text{ on }D_{\ov T},\qquad\rho'=(\pt_x\rho)_{|D_{\ov T}}\in L^2(D_{\ov T}).
\ee
\item Writing $M(t):=\MM(\rho(t))$, we have   $M\in BV(0,{\ov T})$.\smallskip
\item Moreover, up to extraction of a subsequence, there hold for  $0<T<\ov T$, as $\tau\dw0$,
\begin{align}
\label{XtautoX}
\tilde{X}_\tau\to X\quad\text{in }\,\Co([0,T])\quad&\text{with }\sup_\tau\Lip(\tilde{X}_\tau)\le C_b,\\
\label{MtautoM}
M_\tau\to M\quad\text{in }\, L^1(0,T)\quad&\text{with }\sup_\tau\|M_\tau\|_{BV(0,T)}\le C,\\
\label{rhotau_lim}
\rho_\tau\to\rho\quad\text{in }\, L^1(\R_+\t(0,T))\quad&\text{with }\|\rho_\tau\|_{L^\oo(\R_+)}\le\rot,
\end{align}
\be
\label{rhotau'_lim}
\rho_\tau'\to\rho'\quad\mbox{weakly in }L^2(\R_+\t (0,T)).
\ee
\end{enumerate}
\end{prop}

\medskip

\begin{proof}~

\noindent  
\textit{Step 1. Convergence of $\tilde X_\tau$ and definition of $X$ and $\ov T$.} By Proposition~\ref{prop.Ttau}(d) there exists $\hat T>0$ such that $T_\tau\ge\hat T$ for $\tau$ small enough.
Hence ,
\[
\TT\ge T_0:=\limsup_{\tau\dw0}T_\tau\ge\hat T>0.
\]
By~\eqref{bounds.XtauMtau} of Proposition~\ref{prop_apriori}, $ \tilde X_\tau$ is $C_b$-Lipschitz continuous for every $\tau>0$ and $\tilde X_\tau(0)=X^0$. Therefore, there exists a Lipschitz function $X:[0,T_0]\to\R$ such that $X(0)=X^0$, $\Lip(X)\le C_b$ and, up to extraction, 
\[
T_\tau\to T_0\qquad\text{and}\qquad \tilde X_\tau\to X\quad\text{ in }\Co([0,T])\quad\text{ for every }0<T<T_0.
\] 
We set
\[
\ov T:=\sup\{0<T<T_0:\text{such that }X>0\text{ on }[0,T]\}.
\]
Let us assume by contradiction that $\ov T<\TT$ and that $X(\ov T)>0$, equivalently $\min_{[0,\ov T]}X>0$. This implies $\ov T=T_0$ and by definition of $T_0$ we have, as $\tau\dw0$,
\[
T_\tau=N_\tau\tau\to T_0\le\TT.
\]
By Proposition~\ref{prop.Ttau}(c), we get that $X_\tau^{N_\tau}\le C_c\sqrt\tau$ for $\tau$ small enough. Sending $\tau$ to $0$ and taking into account $\Lip(\tilde X_\tau)\le C_b$, this leads to $X(T_0)=\lim_{\tau\dw0}X_\tau^{N_\tau}=0$. This contradicts the assumption. We conclude that $\ov T<\TT\implies X(\ov T)=0$. Now, using $\Lip(X)\le C_b$, we deduce, 
\[
\ov T<\TT\ \implies\ X^0-0=X(0)-X(\ov T)\le C_b\ov T\ \implies\ \ov T\ge X^0/C_b.
\] 
At this point we have established that $\ov T$ and $X$ satisfy points (a)\&(b) of the proposition and the convergence property~\eqref{XtautoX} of point~(e).\medskip

\noindent 
\textit{Step 2. Convergence of $M_\tau$ and $\rho_\tau$ on a time interval $[0,T]$, $0<T<\ov T$.} Let $T\in(0,\ov T)$. Recalling that $M_\tau(0)=\MM(\rho^0)$, we deduce from the last estimate of~\eqref{bounds.XtauMtau} that the sequence $(M_\tau)$ is uniformly bounded in $W^{1,1}(0,T)$. Using a classical compactness criterion (see \textit{e.g.}~\cite[Theorem 5.5]{EG15}), we obtain that there exists $M\in BV(0,T)$ such that, up to further extraction,
\[
M^\tau\rightarrow M\quad\mbox{in }L^1(0,T).
\]
This proves~\eqref{MtautoM}.\\
Now, thanks to the estimates~\eqref{bounds.rhotau}--\eqref{time_translate} of Proposition~\ref{prop_apriori}, we see that  $(\rho_\tau)$ complies to the compactness criteria of~\cite[Theorem 1]{DJ12} with $\mathcal{X}=BV(\R_+)$, $\mathcal{B}=L^1(\R_+)$,  $\mathcal{Y}=H^\ast$ and $p=r=1$ (see Theorem~\ref{thmabstarctcompacness} below). This leads to the existence of $\rho\in L^1(\R_+\times (0,T))$ such that, up to further extraction, we have
\be
\label{conv_forte_Lp}
\rho_\tau\rightarrow\rho\quad\mbox{in }L^1(\R_+\times(0,T)).
\ee
We first deduce from this convergence that
\[
M(t) =\MM(\rho(t)) =\int_{\R_+}\rho(x,t)\, dx\qquad\mbox{for almost every }t\in (0,T),
\]
that is point~(d) of the proposition. Next, up to further extraction, the sequence $(\rho_\tau)$ converges to $\rho$ almost everywhere in $\R_+\times (0,T)$. Besides $\tilde X_\tau$ is uniformly-Lipschitz and $\supp\rho_\tau(t)=[0,X_\tau(t)]$, thereby
%The convergence property~\eqref{conv_forte_Lp} and the fact that by construction $\supp\rho_\tau(t)=[0,\tilde X(t)]$ with $\rob\le\rho_\tau(t)\le\rot$ imply that
\[
\hat{\rho}_-\le\rho\le\hat{\rho}_+\quad\mbox{in }D_T\quad\mbox{ and }\quad\supp\rho =\overline{D_T}.
\]
Eventually, the uniform estimate~\eqref{rho_H1} yields the weak convergence in $L^2(\R_+\times(0,T))$ of a subsequence of $(\rho'_\tau)$ towards a function $w\in L^2(\R_+\times(0,T))$. Using classical arguments we identify this function with $\rho'$ on $\R_+\times(0,T)$ in the sense of distributions.

Eventually, taking $T=T_k:=\ov T(1-1/k)$ for $k\ge 1$, and using a diagonal extraction, we see that we can assume that the subsequence is the same for $0<T<\ov T$.
\end{proof}

For the reader's convenience we state the compactness result used in the preceding proof. It is a time-discrete counterpart of the Aubin--Lions--Simon lemma (see~\cite{S86}).

\begin{theo}{\cite[Theorem 1]{DJ12}}\label{thmabstarctcompacness}
Let $\mathcal{X}\subset \mathcal{B}\subset \mathcal{Y}$ be Banach spaces such that the embedding $\mathcal{X} \hookrightarrow \mathcal{B}$ is compact and the embedding $\mathcal{B} \hookrightarrow \mathcal{Y}$ is continuous. Let $p,r\ge1$, let $\tau_j>0$ be sequence decreasing and let $(v_j)\subset L^p(0,T;\mathcal{X})$ such that for $j\ge1$, $v_j$ is constant on each subinterval $((k-1)\tau_j, k \tau_j)$.\\
We assume that  there exists $C>0$  such that for $j\ge1$,
\[
\tau_j^{-1} \, \|v_j - \sigma_{-{\tau_j}} \, v_j \|_{L^r({\tau_j},T;\mathcal{Y})} + \|v_j\|_{L^p(0,T;\mathcal{X})} \leq C \qquad\text{for every }j\ge1,
\]
where $\sigma_{-\tau}$ is the shift operator defined by~\eqref{sigmatau}. Then $(v_j)$ is relatively compact in $L^p(0,T; \mathcal{B})$.
\end{theo}
%\begin{theo}{\cite[Theorem 1]{DJ12}}\label{thmabstarctcompacness}
%Let $X\subset B\subset Y$ be Banach spaces such that the embedding $X \hookrightarrow B$ is compact and the embedding $B \hookrightarrow Y$ is continuous. Let $p\in[1,\oo]$ and $r\ge1$ with $r>1$ if $p=\oo$, let $\tau_j>0$ be sequence decreasing to 0 and  for $j\ge1$, let $u_j:\R_+\to X$ such that $u_j$ is constant on each subinterval $((k-1)\tau_j, k \tau_j)$.\\
%We assume that  there exists $C>0$ independent of $j$ such that for $j\ge1$,
%\[
%\tau^{-1} \, \|u_j - \sigma_{-{\tau_j}} \, u_j \|_{L^r({\tau_j},T;Y)} + \|u_j\|_{L^p(0,T;X)} \leq C \qquad \forall {\tau_j} > 0,
%\]
%where and $\sigma_{-\tau}$ is the shift operator defined by~\eqref{sigmatau}. Then if $p<\oo$, the sequence $(u_j)$ is relatively compact in $L^p(0,T; B)$ and if $p=\oo$ (and $r>1$), up to extraction, $(u_j)$ converges in each $L^q(0,T; B)$ for $1\leq q <\oo$ towards a function $u\in \Co^0([0,T],B)$.
%\end{theo}

In the rest of this subsection we study the convergence of the traces of $(\rho_\tau)$ at $x=0$ and at $x= X_\tau(t)$.

\begin{prop}\label{prop_traces}
Let $(\rho,X)$ be the limit obtained by Proposition~\ref{Prop_compacite1} and defined on $[0,\ov T)$. Up to the extraction of a subsequence, there hold for every $p\in [1,\oo)$, 
\begin{align}
\label{trace_inter}
\int_0^{\ov T}\lt|\rho_\tau(X_\tau(t),t) -\rho(X(t),t)\rt|^p\, dt &\st{\tau\dw0}\longto 0,\\
\label{trace_zero}
\int_0^{\ov T}|\rho_\tau(0,t) -\rho(0,t)|^p\, dt &\st{\tau\dw0}\longto 0.
\end{align}
\end{prop}

\begin{proof}
The proofs of the two limits are similar, we only establish~\eqref{trace_inter}. Thanks to the $L^\oo$ estimates of~\eqref{bounds.rhotau} and~\eqref{rho_lim} it is sufficient to establish~\eqref{trace_inter} in the case $p=1$, besides, for the same reason, we can replace $\ov T$ by $T\in(0,\ov T)$. Let $0<T<\ov T$ and let $\eps>0$ such that $X>3\eps$ on $[0,T]$. For $\tau$ small enough, there holds $T\le T_\tau$ and 
\[
 0<X(t)-3\eps \le X_\tau(t)-2\eps< X_\tau(t)-\eps\le X(t)\quad \text{for }t\in[0,T].
\]
Using the triangle inequality, we write
\[
\int_0^T\lt|\rho_\tau(X_\tau(t),t)-\rho(X(t),t)\rt|\,dt\le Q^1_\tau(\eps) + Q^2_\tau(\eps) + Q^3_\tau(\eps),
\]
with
\begin{align*}
Q^1_\tau(\eps)&:=\frac1\eps\int_0^T\int_\eps^{2\eps}\lt|\rho_\tau(X_\tau(t))-\rho_\tau(X(t)-y)\rt|\, dy\, dt,\\
Q^2_\tau(\eps)&:=\frac1\eps\int_0^T\int_\eps^{2\eps}\lt|\rho_\tau(X(t)-y) -\rho(X(t)-y)\rt|\, dy\, dt,\\
Q^3_\tau(\eps)&:=\frac1\eps\int_0^T\int_\eps^{2\eps}\lt|\rho(X(t)-y) -\rho(X(t)\rt|\, dy\, dt.
\end{align*}
For the second term, we have by~\eqref{rhotau_lim} of Proposition~\eqref{Prop_compacite1},
\be\label{Q2rhoX}
Q^2_\tau(\eps)\le\dfrac1\eps\|\rho_\tau -\rho\|_{L^1(\R_+\times (0,T))} \ \st{\tau\dw0}\longto\ 0.
\ee
For the first term, we first write for $y\in(\eps,2\eps)$,
\[
\lt|\rho_\tau(X_\tau(t))-\rho_\tau(X(t)-y)\rt|
=\lt|\ds\int_{X(t)-y}^{X_\tau(t)} \rho_\tau'(r)\,dr\rt|
\le\lt(\int_0^{X_\tau(t)} |\rho_\tau'(r)|^2\, dr\rt)^{1/2}\sqrt{3\eps},
\]
where we have used the Cauchy-Schwarz inequality and $|X_\tau(t)-(X(t)-y)|\le3\eps$. Integrating over $y\in[\eps,2\eps]$ and $t\in[0,T]$ and dividing by $\eps$, we obtain 
\[
Q^1_\tau(\eps)
\le\sqrt{3\,\eps\,}\int_0^T\lt(\int_{\R_+} |\rho_\tau'(r)|^2\, dr\rt)^{1/2}dt
\le\sqrt{3\,T\eps\,}\,\|\rho_\tau'\|_{L^2(D_T)}.
\]
According to~\eqref{rhotau'_lim}, there holds, $\|\rho_\tau'\|_{L^2(D_T)}\le C$ for $\tau>0$ small enough and some $C\ge 0$ that only depends on $(\rho^0,X^0)$, $T$ and on the parameters of the model. We have the same estimate for the term $Q^3_\tau(\eps)$ and we get
\be\label{Q13rhoX}
Q^1_\tau(\eps)+Q^3_\tau(\eps)\le C\sqrt\eps.
\ee
Eventually, sending $\eps$ to 0, we deduce from~\eqref{Q2rhoX},~\eqref{Q13rhoX} and a diagonal argument that~\eqref{trace_inter} holds true. The proof of~\eqref{trace_zero} is similar.
\end{proof}

We can now establish that the boundary constraint $\rom\le\rho(0,t)\le\rop$ holds true for almost every $t\in(0,{\ov T})$.

\begin{prop}
\label{prop_0_rho}
The function $\rho$ obtained in Proposition~\ref{Prop_compacite1} satisfies
\[%\be\label{5.esttrace}
\rom\le\rho(0,t)\le\rop\qquad\mbox{for almost every }t\in (0,{\ov T}).
\]%\ee
\end{prop}

\begin{proof} 
Let $\tau>0$ and $0\le n< N_\tau$ and let us consider the case $\rho^n>\rop$. Recalling the definition~\eqref{2.def.penalty}  of the penalty function $p_\tau$, we  have, thanks to the boundary conditions of point~(b) of Proposition~\ref{prop.rhonXn},
\[
\lt(\log(\rho^{n+1}(0))-\log(\rop)\rt)^2_+ =\lt(p'_\tau(-m^n)\rt)^2 = 2\,\tau^{\omega-1}\, p_\tau(m^n)\quad\text{ for }0\le n\le N_\tau-1.
\]
Multiplying by $\tau$, summing over $n$ and using the estimate~\eqref{e.a.p} we deduce,
\[
\int_0^T\lt(\log\rho_\tau(0,t)-\log\rop\rt)^2_+\, dt\le 2\,\tau^\omega\sum_{n=0}^{N_\tau-1}p_\tau(m^n)\le 2\,\tau^\omega\lt({\bf\Psi}(\rho^0,X^0)-C_e\rt),
\]
for any $0<T<\ov T$ and $\tau$ small enough (so that $T_\tau\ge T$).\\
Finally, sending $\tau$ to 0 and recalling that $\omega\in (0,1)$, the above estimate and~\eqref{trace_zero} of Proposition~\ref{prop_traces} lead to 
\[
\int_0^T\lt(\log\rho(0,t)-\log\rop\rt)^2_+\, dt\le0.
\]
This implies $\rho(0,t)\le \rop$ for almost every $t\in (0,\ov T)$. The proof of the lower bound is similar. This concludes the proof of Proposition~\ref{prop_0_rho}.
\end{proof}

\subsection{Proof that the limit $(\rho,X)$ is a weak solution of~\eqref{P}-\eqref{P.bord1}}\label{subsec.proof.wealsol}~

From now on, $(\rho,X)$ is the limit obtained by Proposition~\ref{Prop_compacite1} and defined on $(0,\ov T)$. In this subsection, we show that the pair $(\rho, X)$ is a weak solution of~\eqref{P}-\eqref{P.bord1} in the sense that the variational formulation~\eqref{wf_0}-\eqref{def:dotM} holds true. 

We first show that, up to a controlled error term, the pair $(\rho_\tau,X_\tau)$ is a weak solution.

\begin{prop}
\label{prop_inega_EL}
Assume that the parameter $\om$ of the definition~\eqref{2.def.penalty} of the penalty function $p_\tau$ satisfies $0<\om < 1/2$ and let $0<T<\ov T$. Then the following inequality holds for $\tau>0$ small enough and $\vhi\in\Co^\oo_c(\R_+\t [0,T))$. 
\begin{multline}
\label{inega_EL1}
 \lt| -\int_0^{T-\tau}\int_{\R_+}\rho_\tau(x,t)\,\dfrac{\sigma_{\tau}\vhi(x,t)-\vhi(x,t)}{\tau}dx\, dt 
 +\int_0^T\int_{\R_+}\rho'_\tau(x,t)\,\pt_x\vhi(x,t)\, dx\, dt\rt.\\
\lt.-\int_0^T\dot{M}_\tau(t)\,\vhi(0,t)\, dt
-\dfrac{1}{\tau}\int_0^\tau\int_{\R_+}\rho^0(x)\,\vhi(x,t)\, dx\, dt
+\dfrac{1}{\tau}\int_{T-\tau}^{T}\int_{\R_+}\rho_\tau(x,t)\,\vhi(x,t)\, dx\, dt\rt|\\
\le R_\tau(\vhi),
\end{multline}
% \sum_{n=0}^{N_\tau-1}\tau\, |Q_\tau^n(\vhi)|
where $\sigma_\tau$ denotes the shift operator defined by~\eqref{sigmatau} and with \\ 
\be
\label{cont_EL_reste2}
 R_\tau(\vhi) \le C\lt(\|\pt_x^2\vhi\|_{L^\oo(\R_+\t [0,T))}+\|\pt_x\vhi\|_{L^\oo(\R_+\t [0,T))}\rt)\tau^{1-2\omega}\ \st{\tau\dw0
}\longto\ 0,
\ee
with $C>0$ a constant only depending on $(\rho^0,X^0)$, $\TT$ and on the parameters of the model. Eventually, for all $\xi\in\Co(0,T)$, we have
\be
\label{FF_EDO_inter}
\lambda\int_0^T\overset{\cdot }{\tilde{X}_\tau}(t)\,\xi(t)\, dt =\int_0^T\rho_\tau(X_\tau(t),t)\,\xi(t)\, dt -\ros\int_0^T\xi(t)\, dt.
\ee
\end{prop}

\smallskip

\begin{proof}
We first notice that~\eqref{inega_EL1} is a direct consequence of point~(a) of Proposition~\ref{prop.rhonXn} and of a discrete integration by parts in time for the term that involves $\rho^{n+1}-\rho^n$. Then~\eqref{FF_EDO_inter} is a consequence of point~(c) of Proposition~\ref{prop.rhonXn}. We are left to establish~\eqref{cont_EL_reste2}. Thanks to point~(a) of Proposition~\ref{prop.rhonXn} we have the following estimate for every $\vhi\in\Co^\oo_c(\R_+\t [0,T))$,
 \[
R_\tau(\vhi)\le\tau\,\|\pt_x^2\vhi\|_{L^\oo}\,\sum_{n=0}^{N_\tau-1}\dfrac{\Wass^2\lt(\rho^{n+1},\rho^n\rt)}{2\tau}
+\|\pt_x\vhi\|_{L^\oo}\,\sum_{n=0}^{N_\tau-1}\int_0^{\ell_-^n}y\,\rho^n(y)\, dy,
\]
where we recall definition~\eqref{ell-n} of $\ell^n_-$. According to point~(e) of Proposition~\ref{prop.rhonXn}, there holds $0\le\ell^n_-\le A\tau^{1-\om}$. This leads to
\[
R_\tau(\vhi)\le\tau\,\|\pt_x^2\vhi\|_{L^\oo}\sum_{n=0}^{N_\tau-1}\dfrac{\Wass^2\lt(\rho^{n+1},\rho^n\rt)}{2\tau}
+\tau^{1-2\omega}\,\dfrac{\rot\, A^2\,\TT}{2}\|\pt_x\vhi\|_{L^\oo}.
\]
Using, the estimate~\eqref{e.a.p} of Proposition~\ref{prop.Ttau}(f), we obtain
\[
R_\tau(\vhi)\le C\,\lt(\|\pt_x^2\vhi\|_{L^\oo}
+\|\pt_x\vhi\|_{L^\oo}\rt)\,\tau^{1-2\omega},
\]
where $C\ge0$ only depends on $(\rho^0,X^0)$, $\TT$ and on the parameters of the model. This completes the proof of Proposition~\ref{prop_inega_EL}.
\end{proof}

\begin{coro}\label{coro_VarForm}
Assume that $0<\om<1/2$, then $(\rho,X)$ is a weak solution of~\eqref{P} on $(0,\ov T)$ in the sense that it satisfies conditions~(a) to~(f) of Definition~\ref{def.sol.faible}.
\end{coro}
\begin{proof}
We have already established in the points~(b),~(c),~(d) of Proposition~\ref{Prop_compacite1} that the pair $(\rho,X)$ satisfies points~(a),~(b),~(c) of Definition~\ref{def.sol.faible}. Then, sending $\tau$ to 0 in~\eqref{inega_EL1} and using the convergences~\eqref{MtautoM},~\eqref{rhotau_lim},~\eqref{rhotau'_lim} of Proposition~\ref{Prop_compacite1} yield identity~\eqref{1.eq_limit_rho}, that is point~(d) of the definition. Next, sending $\tau$ to 0 in~\eqref{FF_EDO_inter} and using~\eqref{trace_inter} of Proposition~\ref{prop_traces} leads to~\eqref{1.eq_limit_X}, that is point~(e) of the definition. Eventually, point~(f) is just Proposition~\ref{prop_0_rho}.
\end{proof}

\subsection{Proof that the limit $(\rho,X)$ satisfies the variational inequality~\eqref{1.ineg.var}}\label{subsec.var}~

In this subsection, we prove that $\rho$ satisfies %the variational inequality~\eqref{1.ineg.var}, that is the
 weak formulation of the boundary condition~\eqref{P.bord2} at $x=0$. We first establish a time-discrete counterpart of~\eqref{1.ineg.var} and then pass to the limit $\tau\dw 0$.

We fix $0 < T < {\ov T}$ and set $X_* =(\min_{[0,T]}X)/2>0$. Notice that for $\tau>0$ small enough, there holds, $X_\tau>X_*$ in $[0,T]$. Let $\chi\in\Co^\oo_c(\R_+)$ satisfying~\eqref{X*chi}, {\it i.e.}
\[
0\le\chi\le1,\qquad\chi\equiv1\ \text{ on }[0,X_*/2)\quad\text{ and }\quad\supp\chi\sub [0,3X_*/4).
\]
In this subsection, we first fix $\tau>0$ and with the notation of~Subsection~\ref{sec.JKOscheme}, we set
\[
w^n(x):=\chi(x)\,\rho^n(x)\quad\text{ for }x\in\R_+\text{ and }0\le n\le N_\tau.
\]

\begin{prop}\label{prop.prel.inegvard}
For $\tau>0$ small enough, there holds for every $\vhi\in\Co^\oo_c(\R_+)$ and $0\le n\le N_\tau-1$,
\begin{multline}\label{prel.inegvard}
\int_0^{X_*}\dfrac{w^{n+1}(x) - w^n(x)}{\tau}\,\vhi(x)\, dx -\dfrac{M^{n+1}-M^n}{\tau}\vhi(0)\\
+\int_0^{X_*}\lt(w^{n+1}\rt)'(x)\,\vhi'(x)\, dx =\int_0^{X_*}g^{n+1}(x)\,\vhi(x)\, dx + Q^n_{\tau}(\chi\vhi),
\end{multline}
with
\[
g^{n+1}(x) :=- 2\lt(\rho^{n+1}\rt)'(x)\,\chi'(x)-\rho^{n+1}(x)\,\chi''(x)\quad\text{ for }\,0\le n\le N_\tau-1,
\]
and
\be\label{reste.prel.inegvard}
|Q^n_\tau(\chi\vhi)|\le\dfrac{\lt\|\lt(\chi\,\vhi\rt)''\rt\|_{L^\oo(\R_+)}}{2\tau}\,\Wass^2\lt(\rho^{n+1},\rho^n\rt)
+\dfrac{\lt\|\lt(\chi\,\vhi\rt)'\rt\|_{L^\oo(\R_+)}}{2\tau}\rot(\ell^n_-)^2.
\ee
\end{prop}

\smallskip

\begin{proof}
Let $\vhi$ and $n$ as in the statement of the proposition. We set $\psi:=\chi\vhi$ and use $\psi$ as a test function in the identity of Proposition~\ref{prop.rhonXn}(a). We obtain
\begin{multline}\label{Qnwn+1}
\int_0^{X_*}\dfrac{w^{n+1}(x) - w^n(x)}{\tau}\,\vhi(x)\, dx -\dfrac{M^{n+1}-M^n}{\tau}\vhi(0) +\int_0^{X_*}\lt(\rho^{n+1}\rt)'(x)\,\lt(\chi\vhi\rt)'(x)\, dx \\
= Q^n_{\tau}(\chi\vhi).
\end{multline}
The bound~\eqref{reste.prel.inegvard} is a direct application of Proposition~\eqref{prop.rhonXn}(a). Now, denoting $\sigma:=\rho^{n+1}$ we have the elementary identities,
\[
\sigma'(\chi\vhi)'=\sigma'\chi\vhi'+\sigma'\chi'\vhi
=(\sigma\chi)'\vhi' -\sigma\chi'\vhi' +\sigma'\chi'\vhi
=(\sigma\chi)'\vhi'  - (\sigma\vhi)'\chi' +2\sigma'\chi'\vhi.
\]
 Hence,
\begin{align*}
\int_0^{X_*}\lt(\rho^{n+1}\rt)'(\chi\vhi)'\,& =\int_0^{X_*}\lt(w^{n+1}\rt)'\vhi'\, + 2\int_0^{X_*}\lt(\rho^{n+1}\rt)'\vhi \chi'\, -\int_0^{X_*}(\rho^{n+1}\vhi)'\chi'\,\\
&=\int_0^{X_*}\lt(w^{n+1}\rt)'\vhi'\, + 2\int_0^{X_*}\lt(\rho^{n+1}\rt)'\vhi\chi'\, +\int_0^{X_*}(\rho^{n+1}\vhi) \chi''\,,
\end{align*}
where we have used an integration by parts for the last integral together with $\supp\chi'\sub(0,X^*)$. Putting this identity in~\eqref{Qnwn+1} we obtain~\eqref{prel.inegvard}. This completes the proof of the proposition.
\end{proof}
\medskip

Now, we wish to consider the variational formulation~\eqref{prel.inegvard}  with test functions of the form $\vhi=(w^{n+1}-\eta)$ where  $\eta$ is smooth and satisfies $\eta(0)\in [\rom,\rop]$. In this case, thanks to the boundary condition of Proposition~\ref{prop.rhonXn}(b), the second term of~\eqref{prel.inegvard} has a sign. Namely,
\[
-(M^{n+1}-M^n)\lt(w^{n+1}(0)-\eta(0)\rt)\ge0,\quad\text{ for }0\le n\le N_\tau-1.
\]
Hence~\eqref{prel.inegvard} leads to the inequality
\begin{multline}\label{w(n+1)discussion}
\int_0^{X^*}\dfrac{w^{n+1}(x) - w^n(x)}{\tau}\, (w^{n+1}-\eta)(x)\, dx +\int_0^{X_*}\lt(w^{n+1}\rt)'(x)\,\lt(w^{n+1}-\eta\rt)'(x)\, dx\\
\le\int_0^{X^0}g^{n+1}(x)\,\lt(w^{n+1}-\eta\rt)(x)\, dx+  Q^{n}_{\tau}\lt(\chi (w^{n+1}-\eta)\rt).
\end{multline}
This requires some justifications.  Indeed, as far as we know, the function $\vhi=(w^{n+1}-\eta)$ is merely of class $H^1$ and it is not clear whether we can use it as the test function in~\eqref{prel.inegvard}.\\
In order to recover the expected variational inequality~\eqref{1.ineg.var} in the limit $\tau\dw 0$, we have to establish that the right hand side of~\eqref{w(n+1)discussion} is negligible as $\tau\dw0$. If we directly estimate $Q^{n}_{\tau}\lt(\chi (w^{n+1}-\eta)\rt)$ by~\eqref{reste.prel.inegvard} we see that we need controlling the second order derivative of $w^{n+1}$. However, we only have a $H^1$-control on $w^{n+1}$. In the sequel, we deal  with these difficulties by convolving $w^{n+1}$ with a smoothing kernel $\zeta_\delta$ with characteristic length $\delta=\delta(\tau)$. The choice of $\delta$ is driven by a trade-off between the smoothness of $w^{n+1}\ast\zeta_\delta$ and the smallness of $w^{n+1}\ast\zeta_\delta-w^{n+1}$.\medskip

For this we need to introduce further notation. First we define a piecewise constant in time function:
\be\label{def_wtau}
w_\tau(t) := w^{n+1}\quad\mbox{for }t\in (n\tau, (n+1)\tau]\quad\text{ and }-1\le n\le N_\tau-1.
\ee
We extend the functions $w^n$ on $\R $ by
\be\label{deftildew}
\tilde{w}^n(x):=
\begin{cases}
 \qquad\qquad 0 &\mbox{for }x\ge X_*,\\
  \qquad\quad w^n(x) &\mbox{for }x\in [0,X_*],\\
        2\, w^n(0) - w^n(-x) &\mbox{for }x\in [-\delta,0],\\
        2\, w^n(0) - w^n(\delta) &\mbox{for }x\le-\delta.
\end{cases}
\ee
We also define $\tilde w_\tau: \R \t[0,T_\tau]\to\R_+$ by
\[
\tilde{w}_\tau(0):=\tilde{w}^0\qquad \text{and}\qquad\tilde{w}_\tau(t):=\tilde{w}^{n+1}\quad\text{ for }t\in(n\tau,(n+1)\tau],\ 0\le n\le N_\tau-1.
\]
Let $\xi\in\Co^\oo_c\ds\lt((-1/2,1/2),\R_+\rt)$ with $\int\xi\, dx =1$ and let us set $\zeta:=\xi\ast\xi$. We get that  $\zeta$ is an even nonnegative function of $\Co^\oo_c(-1,1)$ with $\int\zeta\, dx = 1$. For $\delta>0$, we denote by $\zeta_\delta$ the mollifier  defined as 
\[
\zeta_\delta(y):=\zeta(y/\delta)/\delta\quad\text{ for }y\in\R,
\] 
and we define 
\[
w_{\tau,\delta}(x,t):=\lt(\tilde{w}_\tau\ast\zeta_\delta\rt)(x,t) =\int_\R\tilde{w}_\tau(x-y,t)\,\zeta_\delta(y)\, dy.
\]
Finally, we introduce a new parameter $\omega'\in (0,1)$ which will be fixed at the end of the proof  and we set 
\be\label{defdelta}
\delta=\delta(\tau,\om'):=\tau^{\omega'},
\ee
\medskip

Let us  establish some estimates, used in the proof of the variational inequality.

\begin{lem}\label{lem_bounds_pre_inegvar}
Assume that $\tau>0$ is small enough so that $\delta=\tau^{\omega'} < X_*/2$. Then the functions $\tilde{w}_\tau$ and $w_{\tau,\delta}$ satisfy the following properties.
\begin{enumerate}[(a)]
\item There exist constants $C_1, C_2>0$ depending on $(\rho^0,X^0)$, the parameters of the model and on $\chi$ such that 
\begin{align}
\label{bound1.ext}
\lt\|\tilde{w}_\tau\rt\|_{L^\infty(\R\t (0,T))}&\le2\rot -\rob,\\
\label{bound3.ext}
\lt\|\tilde{w}_\tau'\rt\|_{L^2(\R\t (0,T))}&\le C_1,\\
\label{bound2.ext}
\lt\|\tilde{w}_\tau'\rt\|_{L^\infty(\R\t (0,T))}&\le C_2\tau^{-\omega},
\end{align}
where we recall that the parameter $\omega\in(0,1)$ appears in the definition~\eqref{2.def.penalty} of $p_\tau$.\smallskip
\item The function $w_{\tau,\delta}$ also satisfies the estimates~\eqref{bound1.ext}--\eqref{bound2.ext}. 
\item We have $w_{\tau,\delta}(0,t) = w_\tau(0,t)$ for almost every  $t\in (0,T)$.
\item There exists $C_3 > 0$, only depending on $(\rho^0,X^0)$, $\chi$, $\zeta$ and on the parameters of the model such that
\be\label{moll.bound2}
\lt\| w_{\tau,\delta}''\rt\|_{L^\infty(\R\t (0,T))}\le C_3\,\tau^{-\omega-\omega'}.
\ee
\item There exists a constant $C_4> 0$, only depending on $(\rho^0,X^0)$, $\chi$, $\zeta$ and on the parameters of the model such that
\be\label{moll.bound1}
\int_0^T\int_0^{X_*}\dfrac{(w_{\tau,\delta} - w_\tau)^2}{\tau}\,dx\,dt\le C_4\,\tau^{2\omega'-2\omega-1}.
\ee
\end{enumerate}
\end{lem}

\medskip

\begin{proof}~
The estimate~\eqref{bound1.ext} follows from the definition~\eqref{deftildew} of $\tilde w^n$ and the inequalities $0\le\chi\le1$ and $\rob\le\rho\le\rot$ on $[0,X^*]\t[0,\ov T)$.  Estimate~\eqref{bound3.ext} follows from the estimate~\eqref{rho_H1} of Proposition~\ref{prop_apriori} and ~\eqref{bound2.ext} follows from the estimate~\eqref{bound.rho'oo} in the same proposition.

Point~(b) is a consequence of~(a) and of the standard properties of convolution.

The identity of point~(c) follows from the definition of $\tilde{w}_\tau$ in the neighborhood of $0$ and from the symmetry of the kernel~$\zeta$. 

For point~(d), we compute
%\begin{align*}
%\|(w_{\tau,\delta})''(t)\|_\oo
%=\|\lt[ \tilde{w}_\tau(t)\ast\zeta_\delta\rt]''\|_\oo 
%&\le\|\zeta_\delta'\|_\oo \|\tilde{w}_{\tau}'(t)\|_\oo\\
%&=\|\zeta'\|_\oo\,\delta^{-1}\|\tilde{w}_\tau'(t)\|_\oo
%=\|\zeta'\|_\oo \tau^{-\om'}\|\tilde{w}_\tau'(t)\|_\oo
%\le C_2\|\zeta'\|_\oo \tau^{-\om-\om'}.%\medskip
%\end{align*}
\begin{align*}
\|(w_{\tau,\delta})''(t)\|_{L^\oo}=\|\lt[ \tilde{w}_\tau(t)\ast\zeta_\delta\rt]''\|_{L^\oo} 
&\le\|\zeta_\delta'\|_{L^\oo}  \|\tilde{w}_{\tau}'(t)\|_{L^\oo}\\
&=\|\zeta'\|_{L^\oo}  \delta^{-1} \|\tilde{w}_\tau'(t)\|_{L^\oo}\\
&=\|\zeta'\|_{L^\oo}  \tau^{-\om'} \|\tilde{w}_\tau'(t)\|_{L^\oo}
\le C_2\|\zeta'\|_{L^\oo}  \tau^{-\om-\om'}.
\end{align*}
This proves the estimate~\eqref{moll.bound2} with $C_3=C_2\|\zeta'\|_{L^\oo}$. For point~(e) we have the classical inequality 
\[
\int_0^T\int_0^{X_*}\dfrac{(w_{\tau,\delta} - w_\tau)^2}{\tau}
=\sum_{n=0}^{N_\tau-1}\lt\| (\tilde w^{n+1})*\zeta_\delta-\tilde w^{n+1}\rt\|_{L^2}^2
\le C^2\delta^2 X_* \sum_{n=0}^{N_\tau-1}\lt\| (\tilde w^{n+1})'\rt\|_{L^\oo}^2,
\]
with $C=\int|y|\zeta(y)\,dy$. This leads to 
\[
\int_0^T\int_0^{X_*}\dfrac{(w_{\tau,\delta} - w_\tau)^2}{\tau}
\le C^2\delta^2 X_* N_\tau \sup_n\lt\| (\tilde w^{n+1})'\rt\|^2_{L^\oo}
\le C^2 \delta^2 X_*\dfrac\TT\tau( C_2\tau^{-\om})^2
=C^2(C_2)^2 X_*\TT\tau^{2\om'-2\om-1}.
\]
This establishes inequality~\eqref{moll.bound1} with $C_4= C^2(C_2)^2 X_*\TT$ and ends the proof of the lemma.
\end{proof}
%\begin{align*}
%\int (\zeta*u-u)^2(x)\,dx
%&=\int \lt(\int\zeta(y)[u(x-y)-u(x)]\,dy\rt)^2\,dx\\
%&=\int \lt(\int\zeta(y)\int^y_0u'(x-z)\,dz\,dy\rt)^2\,dx\\
%&\le \int\int\zeta(y)|y|\int_0^\delta|u'(x-z)|^2\,dz\,dy\,dx\\
%&\le \delta^2\int\zeta(y/\delta)|y|/\delta\,dy/\delta \int|u'(x)|^2\,dx\\
%&=\int \lt(\int u'(x-z) \int \zeta(y)\,dy\,dz\rt)^2\,dx
%\end{align*}

We now establish a semidiscrete version of the variational inequality~\eqref{1.ineg.var} as planed. Let $\phi\in\Co^\oo_c([0,T),\R_+)$ and $\eta\in\Co^\oo_c(\R_+\t [0,T))$ with $\eta(0,t)\in [\rom,\rop]$ on $[0,T)$. Taking $\vhi=w_{\tau,\delta} -\eta$ as a test function in~\eqref{prel.inegvard}, multiplying by $\phi$, integrating over $(n \tau, (n+1) \tau ) $ and summing over $0\le n\le N_\tau-1$ we obtain
\begin{multline}\label{VFw-eta}
\int_0^{T}\phi(t)\int_0^{X_*}\dfrac{w_\tau(x,t)-\sigma_{-\tau}w_\tau(x,t)}{\tau}\,\lt(w_{\tau,\delta}-\eta\rt)(x,t)\, dx\,dt\\
-\int_0^T\phi(t)\,\dot{M^\tau }(t)\,\lt(w_{\tau,\delta}-\eta\rt)(0,t)\, dt
+\int_0^T\phi(t)\int_0^{X_*}w_\tau'(x,t)\lt( w_{\tau,\delta}-\eta\rt)'(x,t)\, dx\,dt\\
 =\int_0^T\phi(t)\int_0^{X_*}\, g_\tau(x,t)\,\lt(w_{\tau,\delta}-\eta\rt)(x,t)\, dx\,dt\\
 +\sum_{n=0}^{N_\tau-1}\int_{n\tau}^{(n+1)\tau}\phi(t)\lt[ Q^{n}_\tau(\chi w_{\tau,\delta}(t)) - Q^{n}_\tau(\chi \eta(t))\rt]\, dt,
\end{multline}
where we have used the linearity of the remaining term $Q^{n}_\tau$ and have introduced the function $g_\tau:\R_+\t[0,T_\tau]\to\R$ defined in the obvious way:
\[
g_\tau(0):=g^0\qquad \text{and}\qquad g_\tau(t):=g^{n+1}\quad\text{ for }t\in(n\tau,(n+1)\tau],\ 0\le n<N_\tau.
\]
Since $w_{\tau,\delta}(0,t)=w_\tau(0,t)=\rho_\tau(0,t)$ for almost every $t\in [0,T]$ we obtain from the boundary condition of Proposition~\ref{prop.rhonXn}(b),
\begin{multline*}
-\int_0^T\phi(t)\,\dot{M^\tau }(t)\,\lt(w_{\tau,\delta}-\eta\rt)(0,t)\, dt\\
=-\sum_{n=0}^{N_\tau-1}\int_{n\tau}^{(n+1)\tau}\phi(t)\,\dfrac{\MM\lt(\rho^{n+1}\rt)-\MM\lt(\rho^n\rt)}{\tau}\,\lt(\rho^{n+1}(0)-\eta(0,t)\rt)\, dt
\ge0.
\end{multline*}
Putting this inequality in~\eqref{VFw-eta}, we obtain the following variational inequality:
\be\label{semi.discreteVI}
A_\tau^1 - A_\tau^2 + A_\tau^3 - A_\tau^4\le A_\tau^5 + A_\tau^6 - A_\tau^7,
\ee
with
\begin{align*}
A_\tau^1 &:=\int_0^{T}\phi(t)\int_0^{X_*}\dfrac{ w_\tau(x,t)-\sigma_{-\tau}w_\tau(x,t)}{\tau}\, w_{\tau,\delta}(x,t)\, dx\,dt,\\
A_\tau^2 &:=\int_0^{T}\phi(t)\int_0^{X_*}\dfrac{ w_\tau(x,t)-\sigma_{-\tau}w_\tau(x,t)}{\tau}\,\eta(x,t)\, dx\,dt,\\
A_\tau^3 &:=\int_0^T\phi(t)\int_0^{X_*}w_\tau'(x,t) w'_{\tau,\delta}(x,t)\,dx\,dt,\\
A_\tau^4 &:=\int_0^T\phi(t)\int_0^{X_*}w_\tau'(x,t)\eta'(x,t)\,dx\,dt,\\
A_\tau^5 &:=\int_0^T\phi(t)\int_0^{X_*}g_\tau(x,t)\,\lt(w_{\tau,\delta}-\eta\rt)(x,t)\, dx\,dt,\\
A_\tau^6 &:=\sum_{n=0}^{N_\tau-1}\int_{n\tau}^{(n+1)\tau}\phi(t)\, Q^{n}_\tau(\chi w^{n+1}_\delta)\,dt,\\
A_\tau^7 &:=\sum_{n=0}^{N_\tau-1}\int_{n\tau}^{(n+1)\tau}\phi(t)\, Q^{n}_\tau(\chi\eta(t))\, dt.
\end{align*}
The expressions of $A_\tau^1$ and $A_\tau^2$ involve $\sigma_{-\tau}w_\tau(t)=w_\tau(t-\tau)$ for $0\le t<\tau$. In this range of values of $t$ we have by definition~\eqref{def_wtau},   $w_\tau(t-\tau)=w^0=\chi\rho^0$.

In order to establish the variational inequality~\eqref{1.ineg.var} it remains to pass to the limit $\tau\downarrow 0$ in the above terms. This is done in the next result.

\begin{prop}\label{prop.inegvar}
Assume that $0 <\tau < 1$ is small enough such that $\delta < X_*/2$ and that $\omega,\omega'\in(0,1)$ comply to the conditions
\be\label{cond.vartheta}
3\omega +\omega' < 1,\quad\omega'-\omega > 1/2,\quad\omega' > 2\omega,\quad\omega <1/2.
\ee
We can take for instance $\om=1/16$, $\om'=3/4$. Then, there hold
\begin{align}\label{lim.A1}
\liminf_{\tau\downarrow 0}A_\tau^1&\ge-\int_0^T\dot{\phi}\int_0^{X_*}\dfrac{(w)^2}{2}\,dx\,dt-\phi(0)\int_0^{X_*}\dfrac{(w)^2}{2}(x,0)\,dx,\\
\label{lim.A2}
\lim_{\tau\downarrow 0}A_\tau^2 &=-\int_0^T\int_0^{X_*}w(\dot{\phi}\,\eta+\phi\dot\eta)\,dx\,dt-\phi(0)\int_0^{X_*}w^0(x)\eta(x,0)\,dx,\\  
\label{lim.A3}
\liminf_{\tau\downarrow 0}A_\tau^3 &\ge\int_0^T\int_0^{X_*}\phi\,\lt(w'\rt)^2\, dx\,dt,\\
\label{lim.A4}
\lim_{\tau\downarrow 0}A_\tau^4 &=\int_0^T\int_0^{X_*}\phi w'\,\eta'\, dx\,dt,\\
\label{lim.A5}
\lim_{\tau\downarrow 0}A_\tau^5 &=\int_0^T\int_0^{X_*}\phi\, g\, (w-\eta)\, dx\,dt,
\end{align}
and
\[%\be\label{lim.A6A7}
\lim_{\tau\downarrow 0}A_\tau^6 =\lim_{\tau\downarrow 0}A_\tau^7 = 0.
\]%\ee
\end{prop}

\begin{Remark}\label{Remomom'} For instance, the choice $\om'=3/4$, $\om=1/{16}$ complies to conditon~\eqref{cond.vartheta}.
\end{Remark}
\medskip

\begin{proof}[Proof of Proposition~\ref{prop.inegvar}]~

\noindent
\textit{Step 1. Limits of the terms $A_\tau^2$, $A_\tau^4$, $A_\tau^5$ and $A_\tau^7$.}  Let us first treat the easiest terms. The limit~\eqref{lim.A2} follows from a discrete integration by parts with respect to the time-shift $\sigma_{-\tau}$, the convergence $\rho_\tau\to\rho$ in $L^2((0,X^*)\t(0,T))$ (see~\eqref{rhotau_lim} in Proposition~\ref{Prop_compacite1}(e)) and the convergence $\rho_\tau(0,t)\to\rho(0,t)$ in $L^2(0,T)$ (see~\eqref{trace_zero} in Proposition~\ref{prop_traces}).\\
The limit~\eqref{lim.A4} follows from the weak convergence $\rho_\tau'\to\rho'$  in $L^2((0,X^*)\t(0,T))$ (see~\eqref{rhotau'_lim} in Proposition~\ref{Prop_compacite1}(e)).\\
 The limit~\eqref{lim.A5} also  follows from the convergence  $\rho_\tau\to\rho$ of Proposition~\ref{Prop_compacite1}(e).
 
 Next, using the estimate~\eqref{reste.prel.inegvard} of Proposition~\ref{prop.prel.inegvard}, we compute
 \[
\lt|A_\tau^7\rt|\le \dfrac{\|\phi\|_{L^\oo}}{2}\lt(\lt\| \lt(\chi\,\eta\rt)''\rt\|_{L^\oo} \sum_{n=0}^{N_\tau-1}\Wass^2\lt(\rho^{n+1},\rho^n\rt)+
\lt\| \lt(\chi\,\eta\rt)'\rt\|_{L^\oo}\rot\sum_{n=0}^{N_\tau-1}\tau(\ell_-^n)^2\rt).
\]
Then, the estimates~\eqref{e.a.p} of Proposition~\ref{prop.Ttau}(f) and~\eqref{borneell+ell-} of Corollary~\ref{coroborneM} lead to 
\[
\lt|A_\tau^7\rt|\le C\|\phi\|_{L^\oo} \|\chi\|_{\Co^2} \|\eta\|_{\Co^2}(\tau+\tau^{2-2\om}) \ \st{\tau\dw0}\longto\ 0.\medskip
\]

\noindent
\textit{Step 2. Lower bound for $A_\tau^1$.} We rewrite $A_\tau^1$ as 
\[
A_\tau^1 =\int_0^{T}\int_0^{X_*}\phi\,\dfrac{w_\tau -\sigma_{-\tau}w_\tau}{\tau}w_\tau\, dx\,dt 
+\int_0^{T}\int_0^{X_*}\phi\,\dfrac{ w_\tau -\sigma_{-\tau}w_\tau}{\tau}\, (w_{\tau,\delta}-w_\tau)\, dx\,dt
=: A_\tau^{11}+A_\tau^{12}.
\]
Using the identity $(a-b)a=(a^2-b^2)/2+(a-b)^2/2$ in the integrand of $A_\tau^{11}$, we obtain
\[
A_\tau^{11}=\dfrac{1}{2}\int_0^{T}\phi\int_0^{X_*}\dfrac{\lt( w_\tau\rt)^2-\lt(\sigma_{-\tau}w_\tau\rt)^2}{\tau}\, dx\,dt
+\dfrac{1}{2}\int_0^{T}\phi\int_0^{X_*}\dfrac{\lt( w_\tau -\sigma_{-\tau}w_\tau\rt)^2}{\tau}\,dx\,dt
=:A_\tau^{11a}+A_\tau^{11b}.
\]
Considering the first term, we get after a discrete integration by parts in time,
\[
A_\tau^{11a}= -\dfrac{1}{2}\int_0^{T-\tau}\int_0^{X_*}\dfrac{\sigma_\tau\phi-\phi}{\tau}\,\lt(w_\tau\rt)^2\,dx\,dt 
-\dfrac{1}{2\tau}\int_{-\tau}^0\phi\int_0^{X_*}\lt(w_\tau\rt)^2\,dx\,dt,
\]
where we assume that $\tau$ is small enough so that $\phi\equiv0$ on $(T-\tau,T)$. By~\eqref{def_wtau}, in the last integral, we have $w_\tau(t)=w^0=\chi\rho^0$ for $-\tau<t\le0$. The $L^2$-convergence $\rho_\tau\to\rho$ then leads to
\be\label{lim.A11a}
\lim_{\tau\downarrow 0}A_\tau^{11a}= -\int_0^T\dot{\phi}\int_0^{X_*}\dfrac{(w)^2}{2}\, dx\,dt -\phi(0)\int_0^{X_*}\dfrac{(w)^2}{2}(x,0)\, dt.
\ee
Now, rearranging the remaining terms, we see that 
\begin{align*}
A_\tau^{11b}+A_\tau^{12}
&=-\dfrac12\int_0^T\phi\int_0^{X_*}\dfrac{(w_{\tau,\delta}-w_\tau)^2}{\tau}\,dx\,dt
+\dfrac12\int_0^T\phi\int_0^{X_*}\dfrac{(w_{\tau,\delta}-\sigma_{-\tau}w_\tau)^2}{\tau}\,dx\,dt\\
&\ge-\dfrac12\int_0^T\phi\int_0^{X_*}\dfrac{(w_{\tau,\delta}-w_\tau)^2}{\tau}\,dx\,dt\\
&\st{\eqref{moll.bound1}}\ge\ -\dfrac{C_4\,\|\phi\|_{L^\infty}}{2}\,\tau^{2\omega'-2\omega-1}\to 0,\quad\mbox{as }\tau\downarrow 0.
\end{align*}
With~\eqref{lim.A11a}, this proves~\eqref{lim.A1}.\medskip

\noindent
\textit{Step 3. Lower bound for $A_\tau^3$.} We write
\[
A_\tau^3 =\int_0^T\phi\int_0^{X_*}\lt(w_{\tau,\delta}'\rt)^2\, dx\,dt +\int_0^T\phi\int_0^{X_*}w'_{\tau,\delta}\lt(w'_\tau -w'_{\tau,\delta}\rt)\, dx\,dt =: A_\tau^{31}+A_\tau^{32}.
\]
By Lemma~\ref{lem_bounds_pre_inegvar} there holds on the one hand $\|w_{\tau,\delta}'\|_{L^2}\le C_1$. On the other hand, we deduce from the convergence $\tilde w_\tau\to w$ in $L^2(\R_+\t(0,T))$ and~\eqref{moll.bound1} (with $\om'>\om$) that  $w_{\tau,\delta}\to w$ in $L^2(\R_+\t(0,T))$. This implies that $w'_{\tau,\delta}\to w'=\lt(\chi\rho\rt)'$ weakly in $L^2(\R_+\t (0,T))$ and we get by the lower semicontinuity of the functional $v\mapsto\int_0^T\int_0^{X_*}v^2(x,t)\phi(t)\, dx\,dt$ with respect to weak convergence,  
\be\label{lim.A31}
\liminf_{\tau\downarrow 0}A_\tau^{31}\ge\int_0^T\int_0^{X_*}\phi\lt(w'\rt)^2\, dx\,dt.
\ee
Next, we rewrite the term $A_\tau^{32}$ as
\[
A_\tau^{32}=\int_0^T\phi\int_\R w'_{\tau,\delta} \lt(\tilde{w}'_\tau- w'_{\tau,\delta}\rt)\, dx\,dt -\int_0^T\phi\int_{\R\setminus[0,X_*]}w'_{\tau,\delta}\lt(\tilde{w}'_\tau-w'_{\tau,\delta}\rt)\,dx\,dt =: A_\tau^{32a}+A_\tau^{32b}.
\]
Let us treat the term $A_\tau^{32a}$. We recall that we have $\zeta=\xi\ast\xi$ with $\xi\in\Co^\oo_c(\R)$. Setting $\xi(y):=\xi(y/\delta)/\delta$  we have $\zeta_\delta=\xi_\delta\ast\xi_\delta$ and since $\xi$ is even, there holds for $t\in(0,T)$,
\[
\int_\R w'_{\tau,\delta}(t)\tilde{w}'_\tau(t)\, dx  = \int_\R \tilde{w}'_{\tau}\ast\zeta_\delta(t)\tilde{w}'_\tau(t)\, dx  = \int_\R \lt[(\tilde{w}'_{\tau}(t)\ast\xi_\delta)\ast\xi_\delta\rt]\tilde{w}'_\tau(t)\, dx = \int_\R \lt(\tilde{w}'_{\tau}(t)\ast\xi_\delta\rt)^2\,dx.
\]
Thereby, 
\[
\int_\R w'_{\tau,\delta} \lt(\tilde{w}'_\tau- w'_{\tau,\delta}\rt)\, dx=
\int_\R \lt(\tilde w'_{\tau}(t)\ast\xi_\delta\rt)^2\, dx
-\int_\R \lt(\lt[\tilde w'_{\tau}(t)\ast\xi_\delta\rt]\ast\xi_\delta\rt)^2\, dx\ge0.
\]
We get the inequality from  $\|z\ast\xi_\delta\|_{L^2}\le\|\xi_\delta\|_{L^1}\|z\|_{L^2}=\|z\|_{L^2}$ with $z=\tilde w'_{\tau}(t)\ast\xi_\delta\in L^2(\R)$. Now, since $\phi\ge0$ we deduce that $A_\tau^{32a} \ge 0$.\\
For $A_\tau^{32b}$ we observe that $w'_{\tau,\delta}(x,t) = 0$ for almost every $x \in\R\setminus (-2\delta, X_*+\delta)$. Therefore, applying the $L^\oo$ bounds~\eqref{bound2.ext} of Lemma~\ref{lem_bounds_pre_inegvar} on $\tilde{w}'_\tau$ and $w'_{\tau,\delta}$ we get
\[
\lt|A_\tau^{32b}\rt|\le 6\, C_2^2\,\|\phi\|_{L^1}\delta\tau^{-2\om}=6\, C_2^2\,\|\phi\|_{L^1}\tau^{\om'-2\om} \ \st{\tau\dw0}\longto\ 0.
\]
With~\eqref{lim.A31} and $A_\tau^{32a} \ge 0$ this proves~\eqref{lim.A3}.
\medskip

\noindent
\textit{Step 4. Control of $A_\tau^6$.}
The estimate~\eqref{cont_EL_reste2} of Proposition~\ref{prop_inega_EL} leads to
\begin{align*}
\lt|A_\tau^6\rt|
&\le C\,\|\phi\|_{L^\oo}\,\lt(\lt\|\lt(\chi w_{\tau,\delta}\rt)''\rt\|_{L^\oo}+\lt\|\lt(\chi w_{\tau,\delta}\rt)'\rt\|_{L^\oo}\rt)\tau^{1-2\omega}\\
&\le C'\lt(\|\tilde w_\tau\|_{L^\oo}+\|\tilde w_\tau'\|_{L^\oo}+\|w_{\tau,\delta}''\|_{L^\oo}\rt)\tau^{1-2\omega}\\
&\st{\eqref{moll.bound2}}\le\ C''\tau^{1-3\omega-\omega'}\  \st{\tau\dw0}\longto\ 0,
\end{align*}
where the constants $C',C''\ge0$ only depend on $\phi$, $\chi$, $(\rho^0,X^0)$ and on the parameters of the problem. This completes the proof of Proposition~\ref{prop.inegvar}.
\end{proof}

\subsection{Proof of the main result, Theorem~\ref{th.main}}~
%\label{subsecEND}~

We already have established in Corollary~\ref{coro_VarForm} that, assuming $0<\om<1/2$, the pair $(\rho,X)$ obtained in Proposition~\ref{Prop_compacite1} satisfies conditions~(a) to~(f) of Definition~\ref{def.sol.faible}.\\
Choosing $\om,\om'\in(0,1)$ satisfying~\eqref{cond.vartheta} (for instance $\om=1/16$, $\om'=3/4$ as in Remark~\ref{Remomom'}), and putting the limits of Proposition~\ref{prop.inegvar} in~\eqref{semi.discreteVI} we get that $(\rho,X)$ also satisfies the variational inequality~\eqref{1.ineg.var} which is condition~(g) of Definition~\ref{def.sol.faible}. We conclude that $(\rho,X)$ is a weak solution of~\eqref{P}--\eqref{P.bord} in the sense of Definition~\ref{def.sol.faible}. 

We now send $\TT$ to infinity. Let $\TT_k:=2^k$ for $k\ge0$. With obvious notation we get a sequence of solutions $(\rho_k,X_k)$ defined on $[0,\ov T_k]$ with $\ov T_k$ nondecreasing and such that $(\rho_{k+1},X_{k+1})=(\rho_k,X_k)$ on $[0,\ov T_k]$. We set $\ov T:=\lim_{k\up\oo}\ov T_k$ and define $(\rho,X)$ by  
\[
(\rho(t),X(t)):=(\rho_k(t),X_k(t))\quad\text{ for }0\le t<\ov T_k.
\] 
We have either $\ov T_k\to+\oo$ or the sequence $\ov T_k$ is stationary after some index $k_0\ge 1$. In this case $\ov T=\ov T_{k_0}$ and by Proposition~\ref{Prop_compacite1}(b) we have $X(\ov T)=0$ so that  $(\rho,X)$ satisfies the conclusion (a) of Theorem~\ref{th.main}.\\
We also have $\rob\le\rho\le\rot$ on $\supp\rho=\ov{D_{\ov T}}\,$ so that the conclusion~(b) of the theorem holds true. Finally, the conclusion~(c) follows from Proposition~\ref{prop.Ttau}(a) and from the lower semicontinuity of $\PSI$. This concludes the proof of Theorem~\ref{th.main}.

%%%%%%%%%%%%%%%%%%%%%%%%%%%%%%%%%%%%%%%%%%%%%%%%%%%%%%%%%%%%%%%%%%%%%%%%%%%%%%%%%%%%%%%%%

\nocite{TheseZ19}

\bibliographystyle{plain}
\bibliography{Biblio}

\end{document}